\RequirePackage{fix-cm}
\documentclass[smallextended]{svjour3}       
\smartqed  
\usepackage{graphicx}

\usepackage{bm}

\RequirePackage[OT1]{fontenc}
\RequirePackage{amsmath}

\RequirePackage[colorlinks,citecolor=blue,urlcolor=blue]{hyperref}

\usepackage{graphicx}
\usepackage{amsmath}
\usepackage{graphics}
\usepackage{epsfig}
\usepackage{amssymb}
\usepackage{color}
\usepackage{float}
\usepackage[section]{placeins}
\usepackage{subcaption,graphicx}

\newcommand{\eps}{\varepsilon}
\newcommand{\be}{\begin{eqnarray}}
\newcommand{\ee}{\end{eqnarray}}
\newcommand{\bea}{\begin{eqnarray*}}
\newcommand{\eea}{\end{eqnarray*}}

\newcommand{\ve}{\varepsilon}


\begin{document}

\title{Covering of high-dimensional cubes and quantization  
}


\author{Anatoly Zhigljavsky        \and
        Jack Noonan 
}


\institute{A. Zhigljavsky \at
              School of Mathematics, Cardiff University, Cardiff, CF244AG, UK \\
              \email{ZhigljavskyAA@cardiff.ac.uk}           
           \and
           J. Noonan \at
               School of Mathematics, Cardiff University, Cardiff, CF244AG, UK \\
              \email{Noonanj1@cardiff.ac.uk}
}

\date{Received: date / Accepted: date}

\maketitle
\begin{abstract}
As the main problem, we consider covering of a $d$-dimensional cube by $n$ balls with reasonably large   $d$ (10 or more) and reasonably small $n$, like $n=100$ or $n=1000$. We do not require the full coverage but only 90\% or 95\% coverage.  We establish  that efficient covering schemes have several important properties which are not seen in small dimensions and in asymptotical considerations, for very large $n$. One of these properties can be termed  `do not try to cover the  vertices' as the vertices of the cube and their close neighbourhoods are very hard to cover and for large $d$  there are far too many of them. We clearly demonstrate that, contrary to a common belief, placing balls at points which form a low-discrepancy sequence in the cube, makes for a very inefficient covering scheme. For a family of random coverings, we are able to provide very accurate approximations to the coverage probability. We then extend our results to the problems of coverage of a cube by smaller cubes and  quantization, the latter being also referred to as  facility location. Along with  theoretical considerations and derivation of approximations, we discuss results of a large-scale numerical investigation.
 MSC 2010: 	90C26,	65K99, 	65B99

\keywords{covering \and quantization  \and facility location \and space-filling \and computer experiments \and high dimension}
\end{abstract}

\section{Introduction}

In this paper, we develop and study efficient schemes for covering and quantization  in high-dimensional cubes.  In particular, we will demonstrate that the proposed schemes are much superior to the so-called `low-discrepancy sequences'.
The paper starts with introducing the main notation, then we formulate the main problem of covering a $d$-dimensional cube by $n$ Euclidean balls. This is followed by a discussion on the main principles we have adopted for construction of our algorithms.
Then we briefly formulate problems   of covering a cube by smaller cubes (which are balls in the $L_\infty$-norm) and the problem of quantization. Both problems have many similarities with the main problem of covering a cube by $n$ balls.
At the end of this section, we describe the structure of the remaining sections of the paper and summarize our main findings.

\subsection{Main notation}

\begin{itemize}

  \item $\mathbb{R}^d$: $d$-dimensional space;
   \item  $\|\cdot \|$ and $\|\cdot \|_\infty$: Euclidean and $L_\infty$-norms in $\mathbb{R}^d$;
  \item ${\cal B}_d(Z,{ r })= \{ Y \in \mathbb{R}^d: \| Y-Z \| \leq { r } \}$: $d$-dimensional ball of radius $r$ centered at $Z \in \mathbb{R}^d$;
        \item  ${\cal B}_d({ r })= {\cal B}_d(0,{ r })=\{ Y \in \mathbb{R}^d: \| Y \| \leq { r }\}$;
        \item ${\cal S}_{d}(Z,r)=  \{ Y \in \mathbb{R}^d: \| Y-Z \| = { r } \}$: $d$-dimensional sphere of radius $r$ centered at $Z \in \mathbb{R}^d$;
            \item
        ${\cal C}_{d}(Z,\delta)=\{ Y \in \mathbb{R}^d: \| Y -Z\|_\infty \leq { \delta } \}$: $d$-dimensional cube of side length $2 \delta$ centered at $Z$ (it is also the $d$-dimensional ball in the $L_\infty$-norm with radius $\delta$ and center  $Z$);
        \item
        ${\cal C}_{d}(\delta)=[-\delta,\delta]^d={\cal C}_{d}(0,\delta) $;
        \item ${\cal C}_{d}=[-1,1]^d={\cal C}_{d}(1) $.
        \end{itemize}

\subsection{Main problem of interest}
\label{sec:main_intro}

The main problem discussed in the paper is the following problem of covering a cube by $n$ balls.
Let ${\cal C}_{d}=[-1,1]^d$ be a $d$-dimensional cube, $Z_1, \ldots, Z_n$ be some points in $\mathbb{R}^d$ and ${\cal B}_d(Z_j,r)$ be the corresponding balls of radius $r$ centered at $Z_j$ $(j=1, \ldots, n)$.
The dimension $d$, the number of balls $n$ and their radius $r$  could be arbitrary.

We are interested in the problem of choosing the locations of the centers of the balls $Z_1, \ldots, Z_n$ so that the union of the balls $\cup_{j}{\cal B}_d(Z_j,r)$ covers the largest possible proportion of the cube ${\cal C}_{d}$. That is, we are interested in choosing a scheme (a collection of points) $\mathbb{Z}_n=\{Z_1, \ldots, Z_n\}$  so that
\be
\label{eq:cover1}
\mbox{$C_d(\mathbb{Z}_n,r):= $vol$({\cal C}_{d} \cap {\cal B}_d(\mathbb{Z}_n,r))/2^d$}
\ee
is as large as possible (given $n$, $r$ and the freedom we are able to use in choosing $Z_1, \ldots, Z_n$). Here
\be
\label{eq:cover1a}
{\cal B}_d(\mathbb{Z}_n,r)= \bigcup_{j=1}^n {\cal B}_d(Z_j,r)
\ee
and $C_d(\mathbb{Z}_n,r)$ is the proportion of the cube ${\cal C}_{d}$ covered by the balls ${\cal B}_d(Z_j,r)$ $(j=1, \ldots, n)$.

For a  scheme $\mathbb{Z}_n$, its covering radius is defined by CR$(\mathbb{Z}_n)= \max_{X \in {\cal  C}_{d}} \min_{Z_j \in \mathbb{Z}_n}\|X-Z_j\|$. In computer experiments, covering radius is called minimax-distance criterion, see \cite{johnson1990minimax} and \cite{pronzato2012design}; in the theory of low-discrepancy sequences, covering radius is called dispersion, see \cite[Ch. 6]{niederreiter1992random}.
The problem of optimal covering of a cube by $n$ balls has very high importance  for the theory of global optimization and many branches of numerical mathematics.
In particular, the celebrated results of A.G.Sukharev  imply that an $n$-point design $\mathbb{Z}_n$ with smallest CR provides the following:
(a)
 min-max $n$-point global optimization method in the set of all adaptive $n$-point optimization strategies, see \cite[Ch.4,Th.2.1]{sukharev2012minimax}, and
  (b) the $n$-point min-max optimal quadrature, see \cite[Ch.3,Th.1.1]{sukharev2012minimax}.
  In both cases, the class of (objective) functions is the class of Liptshitz functions with known  Liptshitz constant.

If $d$ is not small (say, $d>5$) then computation of the covering radius CR$(\mathbb{Z}_n) $ for any non-trivial design $\mathbb{Z}_n$ is a very difficult computational problem. This explains why the problem of construction of optimal $n$-point designs with smallest covering radius is notoriously difficult, see for example recent surveys \cite{toth20172,toth1993packing}.

If $r= $CR$(\mathbb{Z}_n)$, then
 $C_d(\mathbb{Z}_n,r)$ defined in \eqref{eq:cover1} is equal to 1, and  the whole cube ${\cal C}_{d}$ gets covered by the balls.
 However,  we are only interested in reaching the values like 0.9, when a large part of the ball is covered. There are two main reasons why we are not interested in reaching the value $C_d(\mathbb{Z}_n,r)= 1$:
(a) practical impossibility of making a numerical checking of the full coverage, if $d$ is large enough, and
  (b) our approximations lose accuracy when $C_d(\mathbb{Z}_n,r)$ closely approaches 1.

  If, for a given $\gamma \in [0,1)$, we have $C_d(\mathbb{Z}_n,r)\geq 1-\gamma$, then
  the corresponding coverage of ${\cal C}_{d}$ will be called $(1\!-\!\gamma)$-coverage; the corresponding value of $r$ can be called $(1\!-\!\gamma)$-covering radius.  If $\gamma=0$ then the $(1\!-\!\gamma)$-coverage becomes the full coverage {and 1-covering radius of $\mathbb{Z}_n$ becomes $C_d(\mathbb{Z}_n,r)$}.
Of course, for any $\mathbb{Z}_n=\{Z_1, \ldots, Z_n\}$ we can reach $C_d(\mathbb{Z}_n,r)= 1$ by means of increasing $r$. Likewise, for any given~$r$
we can
reach $V_d(\mathbb{Z}_n,r)= 1$ by sending $n \to \infty$. However, we are not interested in very large values of $n$ and try to get the coverage
of the most part of the cube ${\cal C}_{d}$ with the radius  $r$ as small as possible.  We will keep in mind the following typical values of $d$ and $n$:  $d=10,20,50$; $n=64,128,512,1024$.  Correspondingly, we will illustrate our results in
such scenarios.


\subsection{Two contradictory criteria and a   compromise}
\label{sec:two_criteria}

In choosing $\mathbb{Z}_n=\{Z_1, \ldots, Z_n\}$, the following two main criteria
must be followed:
\begin{itemize}
  \item[(i)] the volumes  of intersections of the cube ${\cal C}_{d}$ and each  individual ball ${\cal B}_d(Z_j,r)$ are
not very small;
  \item[(ii)] the volumes of intersections ${\cal B}_d(Z_j,r) \cap {\cal B}_d(Z_i,r)$ are  small  for all $i\neq j$ $(i,j=1, \ldots, n)$.
\end{itemize}

These two criteria do not agree with each other. Indeed, as shown in Section~\ref{sec:quantuty0}, see formulas
\eqref{eq:inters11}--\eqref{eq:inters2f}, the volume of intersection of the ball ${\cal B}_d(Z,r)$  and the cube  ${\cal C}_{d}$ is approximately inversely proportional to $\|Z\|$ and hence criterion (i) favours $Z_j$ with small norms. However, if at least some of  the points $Z_j$ get close to 0, then the distance between these points gets small and, in view of the formulas of Section~\ref{sec:inter}, the volumes of intersections ${\cal B}_d(Z_j,r) \cap {\cal B}_d(Z_i,r)$ get large.

This yields that the above two criteria require a compromise in the rule of choosing $\mathbb{Z}_n=\{Z_1, \ldots, Z_n\}$ as the points $Z_j$ should not be too far from 0 but at the same time, not too close. In particular, and this is clearly demonstrated in many examples throughout the paper, the so-called `uniformly distributed sequences of points' in ${\cal C}_{d}$, including  `low-discrepancy sequences' in ${\cal C}_{d}$, provide poor covering schemes.
This is in a sharp contrast with the asymptotic case  $n \to \infty$ (and hence $r \to 0$), when one of the  recommendations, see \cite[p.84]{janson1986random},  is to choose $Z_j$'s from a uniformly distributed sequence of points from a set which is slightly larger than ${\cal C}_{d}$; this is  to facilitate covering of the boundary of ${\cal C}_{d}$, as
it is much easier to cover the interior of the cube ${\cal C}_{d}$ than its boundary.

In our considerations, $n$ is not very large and hence the radius of balls $r$ cannot be small. One of our recommendations for choosing $\mathbb{Z}_n=\{Z_1, \ldots, Z_n\}$ is to choose $Z_j$'s at random in a cube
${\cal C}_{d}(\delta)=[-\delta,\delta]^d$ (with $0<\delta<1$) with components distributed according to a suitable Beta-distribution.
The optimal value of $\delta$ is always smaller than 1 and depends on $d$ and $n$. If $d$ is small or $n$ is astronomically large, then the optimal value of $\delta$ could be close to 1 but in most interesting instances this value is  significantly smaller than 1. This implies
that the choice  $\delta=1$
(for example, if $Z_j$'s form a uniformly distributed sequence of points in the whole cube ${\cal C}_{d}$)
often leads to very poor covering schemes, especially when the dimension $d$ is large  (see Tables~\ref{table_d_10}--\ref{Table_d_50} in  discussed in Section~\ref{sec:main_1}).
More generally, we show that for construction of efficient designs $\mathbb{Z}_n=\{Z_1, \ldots, Z_n\}$, either deterministic or randomized,
we have to restrict the norms of the design points $Z_j$. We will call this principle `$\delta$-effect'.


\subsection{Covering  a cube by smaller cubes and quantization}

In Section
\ref{sec:cubes}
we consider the problem of $(1-\gamma)$-coverage of  the cube ${\cal C}_{d}=[-1,1]^d$ by smaller cubes (which are $L_\infty$-balls). The problem of 1-covering of cube by cubes has attracted a reasonable attention in mathematical literature, see e.g. \cite{kuperberg1994line,januszewski1994line}. The problem of cube $(1-\gamma)$-covering by cubes  happened to be  simpler than the main
problem of $(1-\gamma)$-coverage of a cube by Euclidean balls and we have managed  to derive closed-form expressions for (a) the volume of intersection of two cubes, and (b) $(1-\gamma)$ coverage, the probability of covering a random point in ${\cal C}_{d}$ by $n$ cubes  ${\cal C}_{d}(Z_i,r)$
for a wide choice of randomized schemes of choosing designs $\mathbb{Z}_n=\{Z_1, \ldots, Z_n\}$. The results of Section
\ref{sec:cubes} show that the $\delta$-effect holds for the problem of coverage of  the cube by smaller cubes
in the same degree as for the main problem of Section~\ref{sec:main_1} of  covering with balls.

Section~\ref{sec:quantization}
is devoted to the following problem of quantization also known as the problem of facility location.
Let  $X=(x_1, \ldots, x_d)$ be uniform on ${\cal C}_{d}=[-1,1]^d$ and $\mathbb{Z}_n=\{Z_1, \ldots, Z_n\}$ be an $n$-point design.
The mean square quantization error is $Q(\mathbb{Z}_n)=\mathbb{E}_X\min_{i=1, \ldots, n} \|X-Z_i\|^2$. In the case where $Z_1, \ldots, Z_n$ are i.i.d. uniform on ${\cal C}_{d}(\delta)$, we  will derive a simple approximation for
the expected value of $Q(\mathbb{Z}_n)$ and clearly demonstrate the $\delta$-effect. Moreover, we will notice a strong similarity between efficient quantization designs and efficient designs constructed in Section~\ref{sec:main_1}.

\subsection{Structure of the paper and main results}

In Section~\ref{sec:quantuty0} we derive accurate approximations for the volume of intersection of an arbitrary $d$-dimensional cube with an arbitrary $d$-dimensional ball. These formulas will be heavily used in Section~\ref{sec:main_1}, which is the main section of the paper
dealing with  the  problem of  $(1-\gamma)$-coverage of a cube by $n$ balls.
In Section
\ref{sec:cubes}
we extend some considerations of Section~\ref{sec:main_1} to  the problem of $(1-\gamma)$-coverage of  the cube ${\cal C}_{d}$ by smaller cubes.
In Section~\ref{sec:quantization} we argue that there is a strong similarity between efficient quantization designs and efficient designs of Section~\ref{sec:main_1}.
In Appendix A, Section~\ref{sec:appA}, we  briefly mention several  facts, used in the main part of the paper, related to high-dimensional cubes and balls. In Appendix B, Section~\ref{sec:appB}, we prove two simple but very important lemmas about distribution and moments of certain random variables.\\

Our main contributions in this paper are:
\begin{itemize}
  \item  an accurate approximation \eqref{eq:inters2f_corrected_2} for the volume of intersection of an arbitrary $d$-dimensional cube with an arbitrary $d$-dimensional ball;
  \item  an accurate approximation \eqref{eq:accurate_app} for the expected volume of intersection of the cube ${\cal C}_{d}$ with
  $n$  balls with uniform random centers $Z_j \in {\cal C}_{d}(\delta)$;
 \item closed-form expression of Section~\ref{sec:closed_form_cubes}  for the expected volume of intersection the cube ${\cal C}_{d}$ with
  $n$  cubes with uniform random centers $Z_j \in {\cal C}_{d}(\delta)$;
  \item construction of efficient schemes of quantization and $(1-\gamma)$-coverage of the cube ${\cal C}_{d}$ by $n$  balls;
  \item large-scale numerical study.
\end{itemize}

We are preparing an accompanying paper \cite{second_paper} in which  we will further explore the topics of Sections~\ref{sec:main_1}-\ref{sec:quantization}
and also consider the problems of quantization and
$(1-\gamma)$-coverage in the whole space $\mathbb{R}^d$ and the problem of $(1-\gamma)$-coverage of simplices.

\section{Volume of  intersection of a cube and a ball }
\label{sec:quantuty0}

\subsection{The main quantity of interest}
\label{sec:quantuty1}

Consider the following problem. Let us take the cube ${\cal C}_{d}=[-1,1]^d$ of volume ${\rm vol}({\cal C}_{d})=2^d$ and a ball
$
{\cal B}_d(Z,{ r })= \{ Y \in \mathbb{R}^d: \| Y-Z \| \leq { r } \}
$  centered at a point  $Z=(z_1, \ldots, z_d)^\top \in \mathbb{R}^d$; this point $Z$ could be outside  $ {\cal C}_{d}$.
Denote the fraction of the cube ${\cal C}_{d}$ covered by the ball ${\cal B}_d(Z,{ r })$ by
\be
\label{eq:inters2a}
C_{d,Z,{ r }}={{\rm vol}({\cal C}_{d} \cap  {\cal B}_d(Z,{ r }))}/2^d\, .
\ee

Our aim is to approximate $C_{d,Z,{ r }}$ for arbitrary $d$, $Z$ and~${ r }$. We will derive a CLT-based normal  approximation in Section~\ref{sec:quantuty2} and then, using an asymptotic expansion in the CLT for non-identically distributed r.v., we will improve this normal approximation in Section~\ref{sec:quantuty2I}. In Section~\ref{sec:quantuty44}  we consider a more direct approach for  approximating $C_{d,Z,{ r }}$
based on the use of  characteristic functions and  the fact that  $ C_{d,Z,{ {r} }}$ is a c.d.f. of $\|U-Z\|$, where $U=(u_1, \ldots, u_d)^\top$ is  random vector with uniform distribution on ${\cal C}_{d}$. From this, $C_{d,Z,{ r }}$ can be expressed through the convolution of one-dimensional c.d.f's. Using this approach we can evaluate the quantity $C_{d,Z,{ r }}$ with high accuracy but the calculations are rather time-consuming. Moreover, entirely new computations have  to be made for different $Z$ and, therefore, we much prefer the approximation of Section~\ref{sec:quantuty2I}.

Note that in the special case $Z=0$, several approximations for the quantity $C_{d,0,{ r }}$ have been derived in \cite{SIAM} but their methods cannot be generalized to arbitrary $Z$.
Note also
that symmetry considerations imply the following relation between $C_{d,0,{ r }}$ and  $C_{d,V,{ r }}=C_{d,Z,{ r }}$ with $\|Z\|=\sqrt{d}$ (when $Z$ is a vertex of ${\cal C}_{d}$) and $r\leq1$:
$
C_{d,V,{ r }}=2^{-d}C_{d,0,{ r }}\, .
$

\subsection{A generalization of the quantity \eqref{eq:inters2a}}
In the next sections, we will need another quantity which slightly generalizes \eqref{eq:inters2a}. Assume
that we have the cube ${\cal C}_{d}(\delta)=[-\delta,\delta]^d$ of volume ${\rm vol}({\cal C}_{d}(\delta))=(2\delta)^d$, the ball
$
{\cal B}_d(Z^\prime,{ r^\prime })= \{ Y \in \mathbb{R}^d: \| Y-Z^\prime \| \leq { r^\prime } \}
$ with a center at a point  $Z^\prime=(z_1^\prime, \ldots, z_d^\prime)^\top $.
Denote the fraction of the cube ${\cal C}_{d}(\delta)$ covered by the ball ${\cal B}_d(Z^\prime,{ r^\prime })$ by
\be
\label{eq:inters2ab}
C^{(\delta)}_{d,Z^\prime,{ r^\prime }}={{\rm vol}({\cal C}_{d}(\delta) \cap  {\cal B}_d(Z^\prime,{ r^\prime }))}/(2\delta)^d\, .
\ee
Then the following change of the coordinates and the radius
\be
\label{eq:inters2ac}
Z=Z^\prime/\delta= (z_1^\prime/\delta, \ldots, z_d^\prime/\delta)^\top \,  \;\;{\rm and }\;\; r= r^\prime/\delta\, .
\ee
gives
\be
\label{eq:inters2ad}
C^{(\delta)}_{d,Z^\prime,{ r^\prime }}=C_{d,Z,{ r }}\, .
\ee

\subsection{Normal approximation for the quantity \eqref{eq:inters2a} }
\label{sec:quantuty2}

Let $U=(u_1, \ldots, u_d)^\top$ be a random vector with uniform distribution on ${\cal C}_{d}$ so that $u_1, \ldots, u_d$ are i.i.d.r.v. uniformly distributed on $[-1,1]$. Then for given $Z=(z_1, \ldots, z_d)^\top \in \mathbb{R}^d$ and any $r>0$,
\be
\label{eq:inters2c}
C_{d,Z,{ r }}\!= \! \mathbb{P} \left\{ \| U\!-\!Z \|\! \leq \! { r } \right\}\!= \! \mathbb{P} \left\{ \| U\!-\!Z \|^2 \leq  { r^2 } \right\}\!= \! \mathbb{P} \left\{\sum_{j=1}^d (u_j\!-\!z_j)^2 \leq  { r }^2 \right\}   .\;\;
\ee
That is,  $C_{d,Z,{ r }}$, as a function of ${ r }$,   is the c.d.f. of the r.v. $\| U-Z \| $.

Let  $u$ have a  uniform distribution on $[-1,1]$ and $|z| \leq 1 $.  In view of Lemma~1 of Section~\ref{sec:appB},  the density of the r.v. $\eta_z = (u-z)^2$ is
\be
\label{eq:inters2d}
\varphi_z(t) = \left\{\begin{array}{ll}
             1/(2\sqrt{t}) & {\rm for\;\;} 0<t\leq (1- |z|)^2\\
              1/(4\sqrt{t}) & {\rm for\;\;} (1- |z|)^2<t\leq (1+ |z|)^2 \\
              0 & {\rm otherwise}
            \end{array}
            \right.
\ee
and
\be
\label{eq:intersf}
\mathbb{E}\eta_z =z^2 +\frac1{3}, \;\;{\rm var}(\eta_z) = \frac43 \left(z^2 +\frac1{15} \right)\, , \;\;
\mu_{z}^{(3)}= \frac{16}{15} \left(z^2 +\frac{1}{63} \right)  \, ,
\ee
where $\mu_{z}^{(3)}$ is the third central moment: $\mu_{z}^{(3)}=E \left[\eta_{z} - E\eta_{z}\right]^3 $.

For   $|z|>1$,  the density of  $\eta_z = (u-z)^2$ is
\be
\label{eq:inters3d}
\varphi_z(t) = \left\{\begin{array}{ll}
              1/(4\sqrt{t}) & {\rm for\;\;} ( 1-|z|)^2<t\leq (1+|z|)^2 \\
              0 & {\rm otherwise}
            \end{array}
            \right.
\ee
with expressions \eqref{eq:intersf} for $\mathbb{E}\eta_z$, ${\rm var}(\eta_z)$ and $\mu_{z}^{(3)}$  not changing.


Consider the r.v.
\be
\label{eq:inters33d}
\| U-Z \|^2 =\sum_{i=1}^d \eta_{z_j}=\sum_{j=1}^d (u_j-z_j)^2\, .
\ee
From \eqref{eq:intersf}, its mean is
\be
\label{eq:inters11}
\mu_{d,Z}=\mathbb{E}\| U-Z \|^2  =\|Z\|^2 +\frac{d}{3}\, .
\ee
Using independence of $u_1, \ldots, u_d$, we also obtain from \eqref{eq:intersf}:
\be
\label{eq:inters12}
 {\sigma}_{d,Z}^2={\rm var}(\| U-Z \|^2 ) = \frac43 \left(\|Z\|^2 +\frac{d}{15}\right)\, \;\;\;\;\;\;
\ee
and
\be
\label{eq:inters123}
 {\mu}_{d,Z}^{(3)}= \mathbb{E}\left[\| U-Z \|^2- \mu_{d,Z}\right]^3  = \sum_{j=1}^d     \mu_{z_j}^{(3)} =
 \frac{16}{15} \left(\|Z\|^2 +\frac{d}{63}\right)\, .\;\;\;\;\;\;
\ee

If $d$ is large enough then the conditions of the CLT for $\| U-Z \|^2$ are approximately met and  the distribution of $\| U-Z \|^2 $
 is approximately normal with mean $\mu_{d,Z}$ and variance ${\sigma}_{d,Z}^2$. That is, we can approximate
$C_{d,Z,{ r }}$
by
\be
\label{eq:inters2f}
C_{d,Z,{ r }} \cong \Phi \left(\frac{{ r }^2-\mu_{d,Z}}{{\sigma}_{d,Z}} \right) \, ,
\ee
where $\Phi (\cdot)$   is the c.d.f. of the standard normal distribution:
$$
\Phi (t) = \int_{-\infty}^t \phi(v)dv\;\;{\rm with}\;\; \phi(v)=\frac{1}{\sqrt{2\pi}} e^{-v^2/2}\, .
$$

The approximation \eqref{eq:inters2f} has acceptable accuracy if $C_{d,Z,{ r }}$ is not very small; for example, it falls inside a $2\sigma$-confidence interval generated by the standard normal distribution; {see Figures~\ref{d_10_z_0}--\ref{d_50_z_0} as examples. Let $p_\beta$ be the quantile of the standard normal distribution defined by $\Phi(\beta)=1-p_\beta$; for example, $p_\beta \simeq 0.05$ for $\beta=2$.
 As follows from \eqref{eq:inters11}, \eqref{eq:inters12} and the approximation \eqref{eq:inters2f}, we expect the approximate inequality $C_{d,Z,{ r }} \gtrapprox p_\beta$ if
\be
\label{eq:cover2}
r \geq R_{d,\|Z\|,\beta}=\left[\|Z\|^2+d/3 - 2\beta   \sqrt{\|Z\|^2/3+d/45 }\right]^{1/2} \, .
\ee
In many cases discussed in Section~\ref{sec:main_1}, the radius $r$ does not satisfy the inequality \eqref{eq:cover2} with $\beta=2$ and even $\beta=3$ and hence the normal approximation \eqref{eq:inters2f} is not satisfactorily accurate; this can be evidenced from Figures~\ref{d_10_z_0} -- \ref{d_50_on_sphere} below.

In the next section, we improve the approximation \eqref{eq:inters2f} by
using an Edgeworth-type  expansion in the CLT for sums of independent  non-identically distributed r.v.

\subsection{Improved normal approximation  }
\label{sec:quantuty2I}

General expansion in the central limit theorem for sums of independent non-identical r.v. has been derived
by V.Petrov, see
Theorem 7 in Chapter 6 in  \cite{petrov2012sums}, see also Proposition 1.5.7 in \cite{rao1987asymptotic}.
 The first three terms of this expansion have been specialized by V.Petrov in
Section 5.6 in \cite{petrov}.
By using only the first term in this expansion,
we obtain the following approximation for the distribution function of $\| U-Z \|^2 $:
\bea
P\left(\frac{\| U-Z \|^2-\mu_{d,Z}}{\sigma_{d,Z}} \leq x \right) \cong \Phi(x) + \frac{ {\mu}_{d,Z}^{(3)}}{6  ({\sigma}_{d,Z}^2)^{3/2} }(1-x^2)\phi(x),
\eea
leading to the following improved  form of \eqref{eq:inters2f}:
\be
\label{eq:inters2f_corrected}
C_{d,Z,{ r }} \cong \Phi(t) + \frac{  \|Z\|^2+d/63}{5\sqrt{3} (\|Z\|^2+d/15)^{3/2} }(1-t^2)\phi(t) \, ,
\ee
where
\be
\label{eq:inters2f_corrected1}
t = t_{d,\|Z\|,{ r }}= \frac{{ r }^2-\mu_{d,Z}}{{\sigma}_{d,Z}}
=  \frac{\sqrt{3}(r^2- \|Z\|^2 -d/3)}{2\sqrt{ \|Z\|^2 +{d}/{15}} }\,.
\ee

From the viewpoint of Section~\ref{sec:main_1}, the range of most important values of $t$ from \eqref{eq:inters2f_corrected1} is $-3 \pm 1$. For such values of $t$,
the uncorrected normal approximation \eqref{eq:inters2f} significantly overestimates the values of $C_{d,Z,{ r }}$, see Figures~\ref{d_10_z_0} -- \ref{d_50_on_sphere} below.
The approximation \eqref{eq:inters2f_corrected}  brings the normal approximation down and makes it much more accurate.
The other terms
in Petrov's expansion of \cite{petrov2012sums} and \cite{petrov} continue to bring   the  approximation down (in a much slower fashion) so that the approximation \eqref{eq:inters2f_corrected} still slightly overestimates the true value of $C_{d,Z,{ r }}$ (at least, in the range of interesting values of $t$ from \eqref{eq:inters2f_corrected1}). However, if $d$ is large enough (say, $d \geq 20$) then the approximation  \eqref{eq:inters2f_corrected} is very accurate and no further correction is needed.

A very attractive feature of the approximations \eqref{eq:inters2f} and  \eqref{eq:inters2f_corrected1} is their dependence on $Z$ through $\|Z\|$ only. We could have specialized for our case the next terms in Petrov's approximation but these terms no longer depend on  $\|Z\|$ only (this fact can be verified from the formula \eqref{eq:inters1c4} for the fourth moment of the r.v. $\nu_z=(z-u)^2$) and hence the next terms are much more complicated. Moreover,  adding one or two extra terms from Petrov's expansion to the approximation \eqref{eq:inters2f_corrected}  does not fix the problem entirely for all $Z$ and $r$. Instead,
we propose a slight adjustment to the r.h.s of \eqref{eq:inters2f_corrected} to improve this approximation, especially for small dimensions. Specifically, we suggest  the approximation
\be
\label{eq:inters2f_corrected_2}
C_{d,Z,{ r }} \cong \Phi(t) + c_d\frac{  \|Z\|^2+d/63}{5\sqrt{3} (\|Z\|^2+d/15)^{3/2} }(1-t^2)\phi(t) \, ,
\ee
where $c_d= 1+ {3}/{d}$ if the point $Z$ lies on the diagonal of the cube ${\cal C}_{d}$ and $c_d= 1+ {4}/{d}$ for a typical (random) point $Z$.
For typical (random) points $Z \in {\cal C}_{d}$, the values of $C_{d,Z,{ r }}$ are marginally smaller than for the points on the diagonal  of ${\cal C}_{d}$ having the same norm, but the difference is very small. In addition to the points on the diagonal, there are other special points: the points whose components are all zero except for one. For such points, the values of $C_{d,Z,{ r }}$ are smaller than for typical points $Z$ with the same norm, especially for small $r$. Such points, however, are of no value for us as they are not typical and we have never observed
in simulations
random points that come close to these truly  exceptional points.

\subsection{Simulation study}

\label{sec:quantuty3}

In Figures~\ref{d_10_z_0} -- \ref{d_50_on_sphere}   we demonstrate  the accuracy of approximations \eqref{eq:inters2f}, \eqref{eq:inters2f_corrected} and \eqref{eq:inters2f_corrected_2} for $C_{d,Z,{ r }} $ in dimensions $d=10, 50$ for the following
 locations of $Z$:
\begin{itemize}
  \item[(i)] $Z=0$, the center of the cube ${\cal C}_{d}$;
  \item[(ii)] $\|Z\|=\sqrt{d}$, with $Z$ being a vertex of  the cube ${\cal C}_{d}$;
    \item[(iii)] $Z$ lies on a diagonal of ${\cal C}_{d}$ with $|z_j|=\lambda \geq 0$ for all $j=1,\ldots,d$ and
     $\|Z\|=\lambda \sqrt{d}$;
      \item[(iv)]  $Z$ is a  random vector uniformly distributed on the sphere
      ${\cal S}_{d}(0,v)$
        with some $v>0$.

\end{itemize}

There are  figures of two types. In the figures of the first type, we plot  $C_{d,Z,{ r }} $ over a wide range of $r$ ensuring that values of $C_{d,Z,{ r }} $  lie in the whole range $[0,1]$. In the figures of the second type, we plot $C_{d,Z,{ r }} $ over a much smaller range of $r$ with $C_{d,Z,{ r }}$ lying in the range $[0,\varepsilon]$ for some small positive $\varepsilon$ such as $\varepsilon = 0.015$. For the purpose of using the approximations of Section~\ref{sec:main_1}, we need to assess the accuracy of all approximations for smaller values of $C_{d,Z,{ r }} $ and hence the second type of plots are often more insightful. In Figures~\ref{d_10_z_0} -- \ref{d_50_on_sphere_non_zoom}, the solid black line depicts values of $C_{d,Z,{ r }} $ computed via Monte Carlo methods, the blue dashed, the red dot-dashed and  green long dashed lines display
 approximations \eqref{eq:inters2f},  \eqref{eq:inters2f_corrected} and
         \eqref{eq:inters2f_corrected_2}, respectively.

In the case where $Z$ is a  random vector uniformly distributed on a sphere
      ${\cal S}_{d}(0,v)$, the style of the figures of the second type is slightly changed to adapt for this choice of $Z$  and provide more information for $Z$ which do or do not belong to the cube ${\cal C}_{d}$. In Fig.~\ref{d_10_on_sphere} and Fig.~\ref{d_50_on_sphere}, the thick dashed red lines correspond to random points $Z \in {\cal S}_{d}(0,v) \cap {\cal C}_{d}$. The thick dot-dashed orange lines correspond to random points $Z \in {\cal S}_{d}(0,v)$ such that $Z\not\in {\cal C}_{d}$.
      Approximations \eqref{eq:inters2f} and \eqref{eq:inters2f_corrected} are depicted  in the same manner as previous figures but  the  approximation \eqref{eq:inters2f_corrected_2}  is now represented by a solid green line.
     The thick solid red line  displays values of $C_{d,Z,{ r }} $ for   $Z$ on the diagonal of ${\cal C}_{d}$ with $\|Z\|=v $ with $v=1.5$ for $d=10$ and
     $v=1.75$ for $d=50$.


\begin{figure}[!h]
\centering
\begin{minipage}{.5\textwidth}
  \includegraphics[width=1\linewidth]{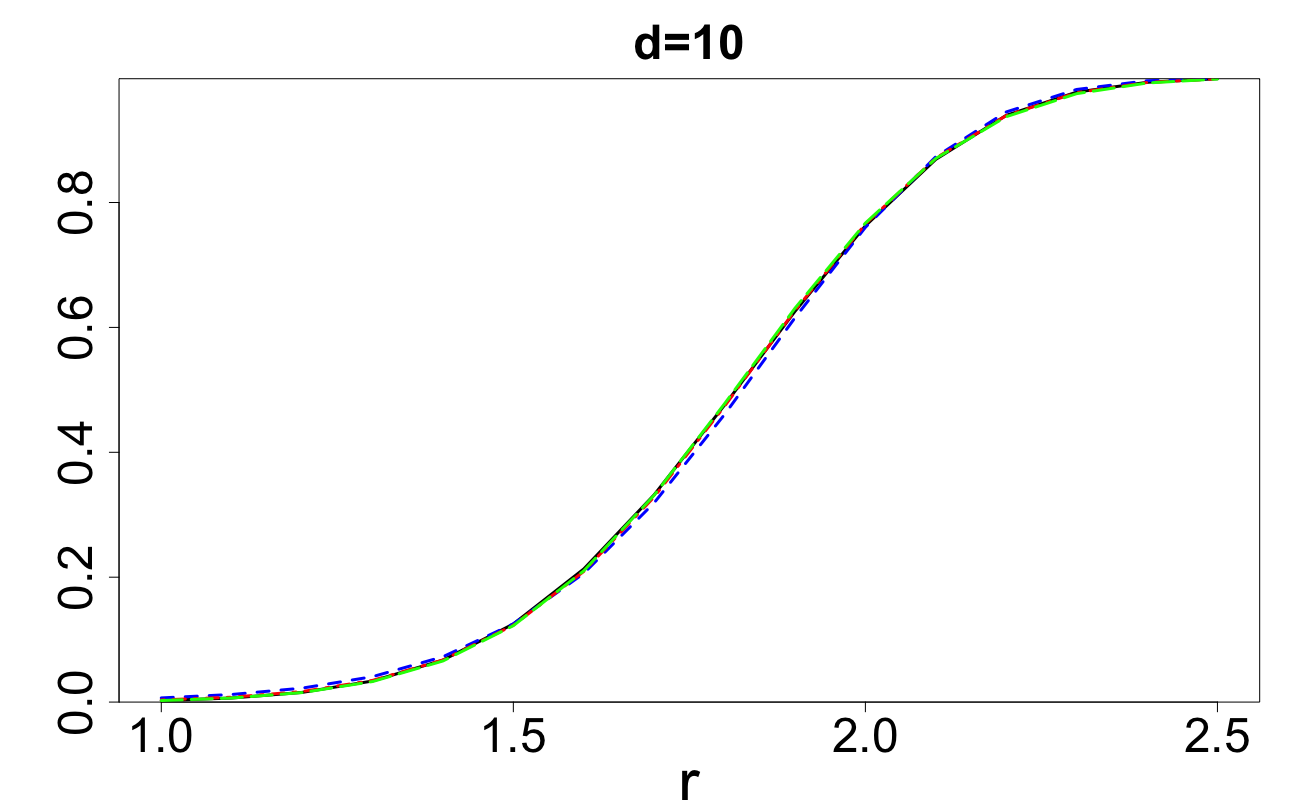}
  \caption{$d=10$, $Z=0$, $r \in [1,2.5]$. }
  \label{d_10_z_0}
\end{minipage}%
\begin{minipage}{.5\textwidth}
  \centering
  \includegraphics[width=1\linewidth]{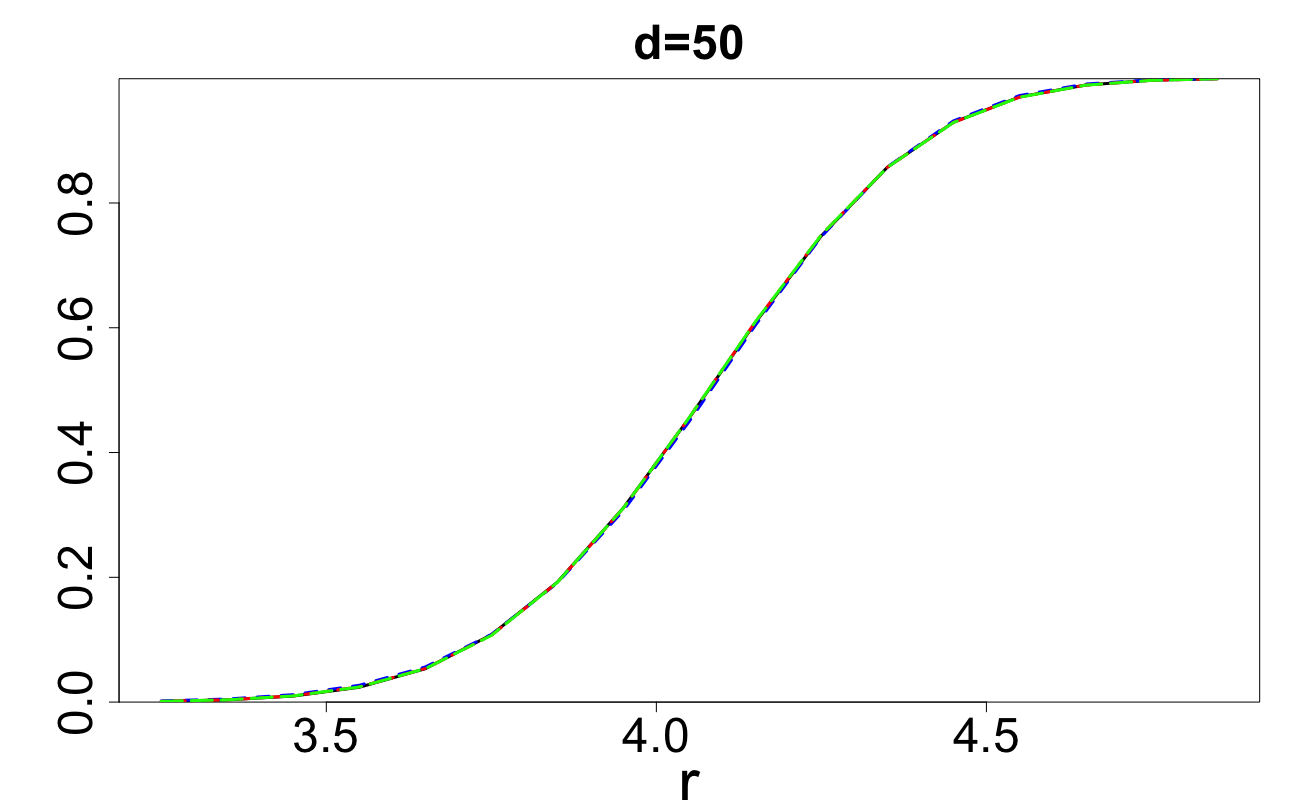}
  \caption{$d=50$, $Z=0$, $r \in [3.2,4.9]$. }
  \label{d_50_z_0}
\end{minipage}
\end{figure}

\begin{figure}[h]
\centering
\begin{minipage}{.5\textwidth}
  \centering
  \includegraphics[width=1\linewidth]{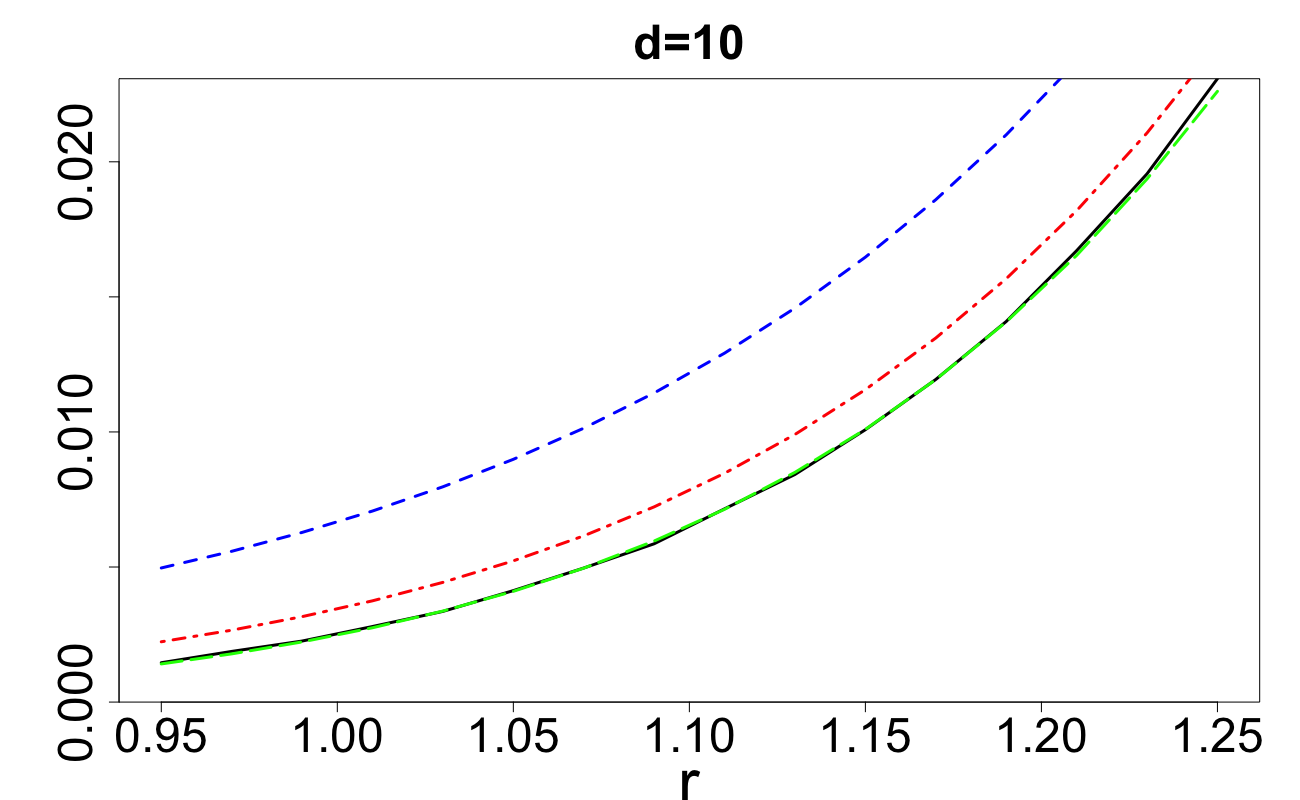}
  \caption{$d=10$, $Z=0$, $r \in [0.95,1.25]$. }
\end{minipage}%
\begin{minipage}{.5\textwidth}
  \centering
  \includegraphics[width=1\linewidth]{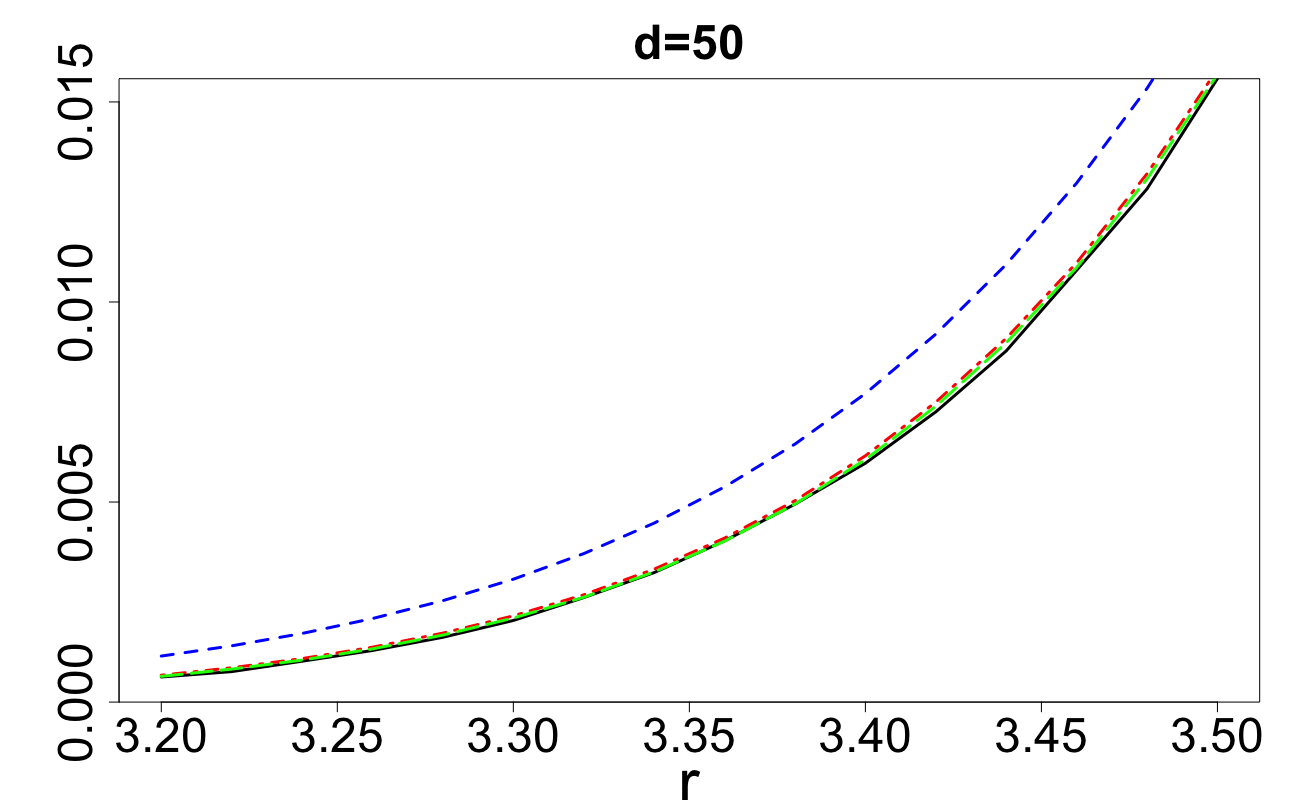}
  \caption{$d=50$, $Z=0$, $r \in [3.2,3.5]$. }
\end{minipage}
\end{figure}

\begin{figure}[!ht]
\centering
\begin{minipage}{.5\textwidth}
  \centering
  \includegraphics[width=1\linewidth]{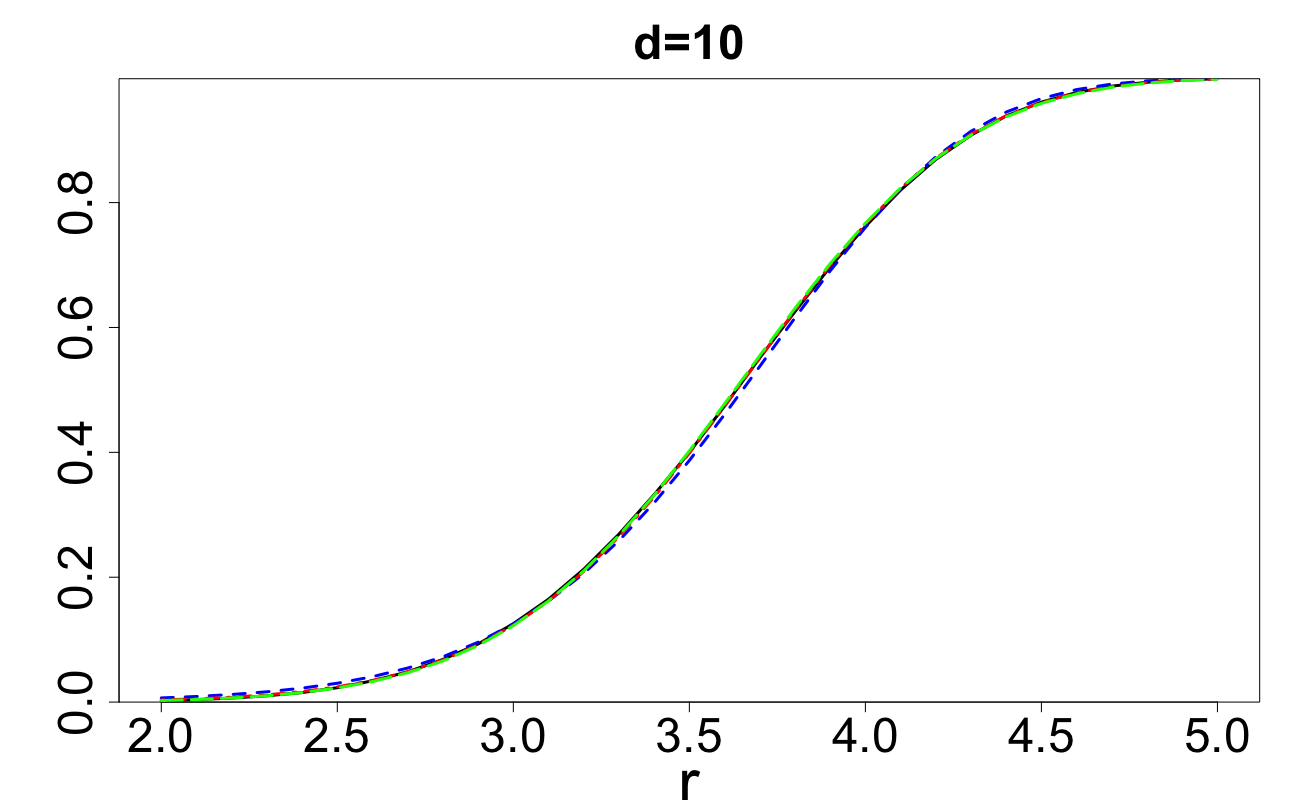}
  \caption{$d=10$, $Z$ is a vertex of ${\cal C}_{d}$, $r \in [2, 5]$. }
\end{minipage}%
\begin{minipage}{.5\textwidth}
  \centering
  \includegraphics[width=1\linewidth]{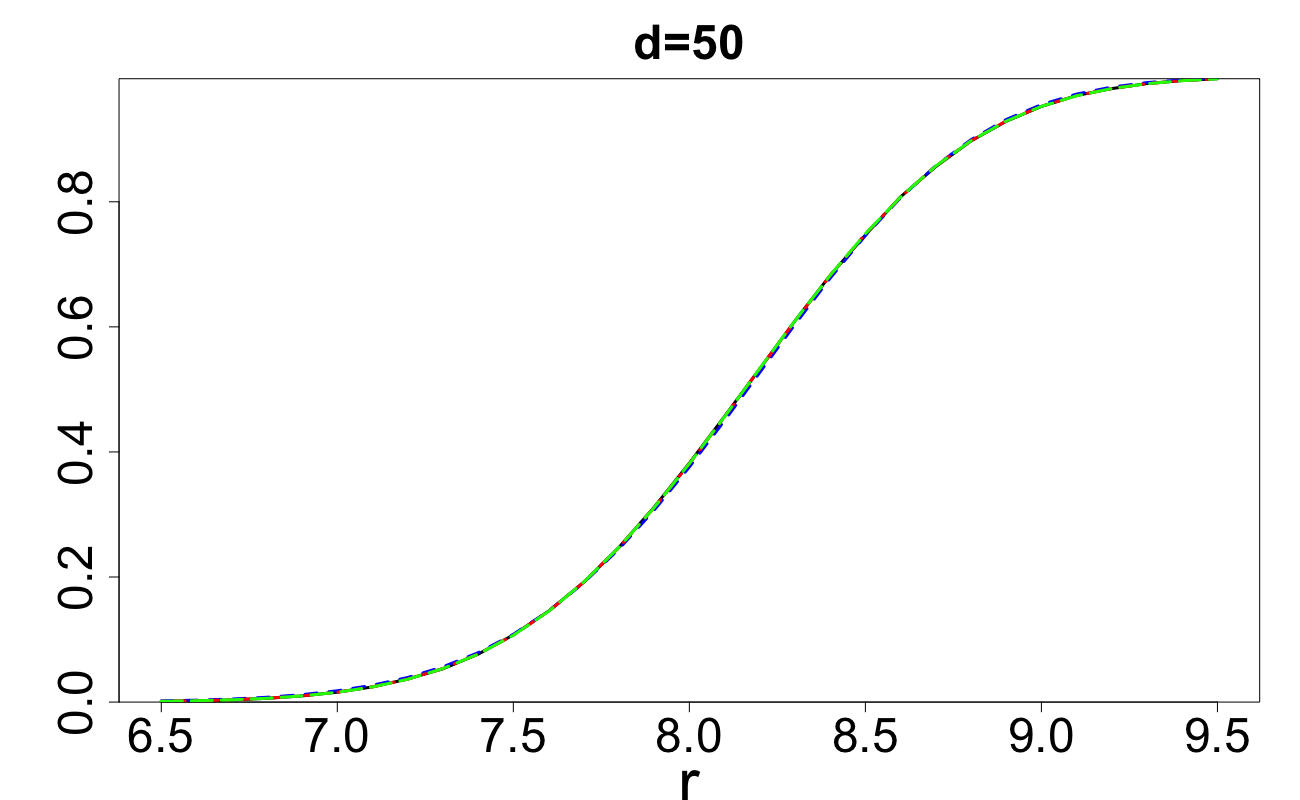}
\caption{$d=50$, $Z$ is a vertex of ${\cal C}_{d}$, $r \in [6.5, 9.5]$. }
\end{minipage}
\end{figure}

\begin{figure}[h]
\centering
\begin{minipage}{.5\textwidth}
  \centering
  \includegraphics[width=1\linewidth]{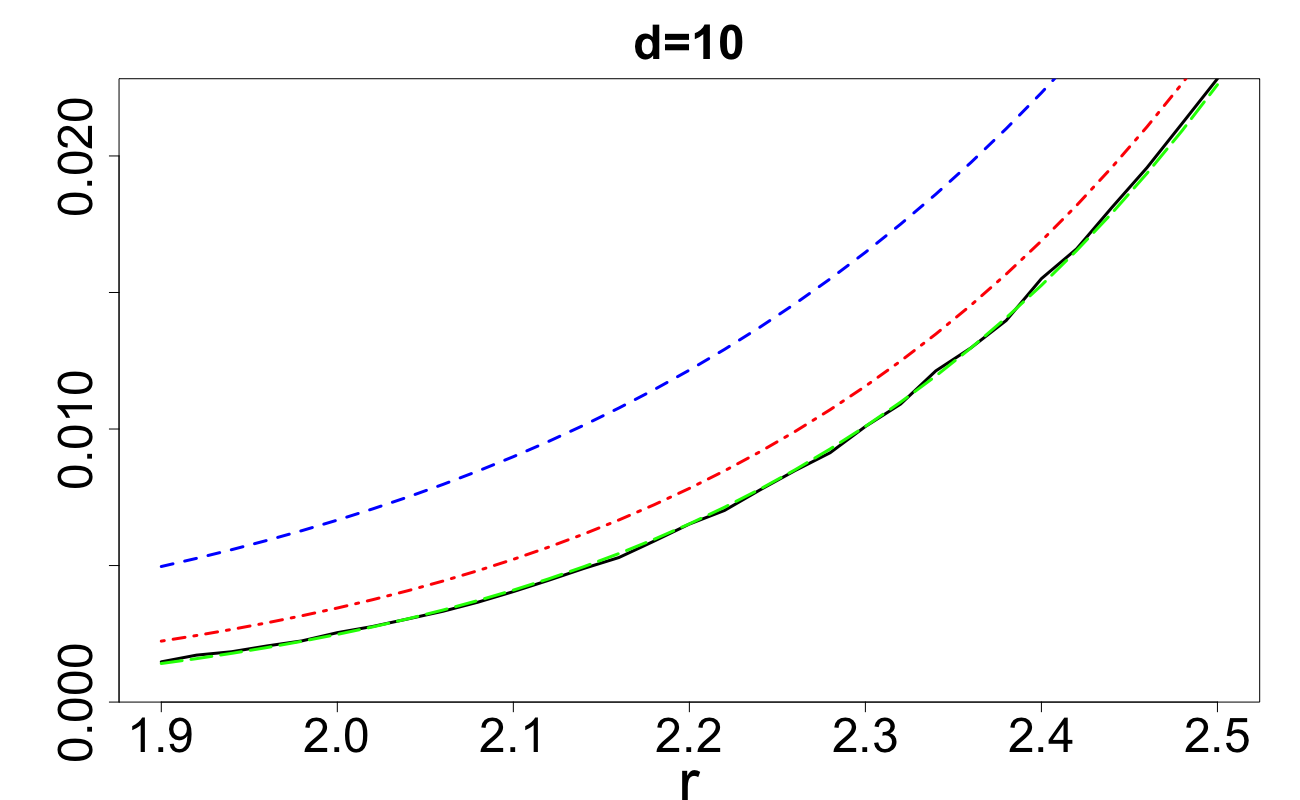}
  \caption{$d=10$, $Z$ is a vertex of ${\cal C}_{d}$, $r \in [1.9, 2.5]$. }
  \label{fig:test1}
\end{minipage}%
\begin{minipage}{.5\textwidth}
  \centering
  \includegraphics[width=1\linewidth]{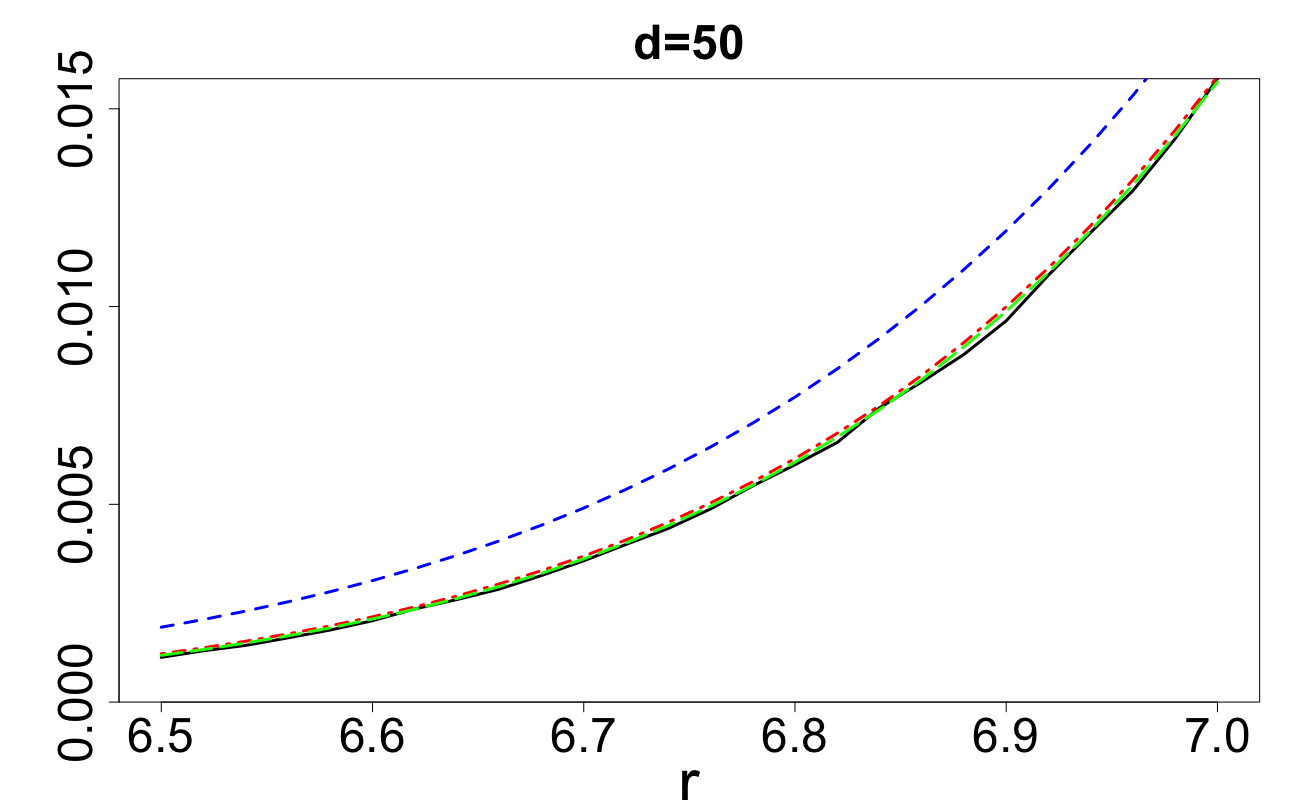}
  \caption{$d=50$, $Z$ is a vertex of ${\cal C}_{d}$, $r \in [6.5, 7]$.}
  \label{fig:test2}
\end{minipage}
\end{figure}

\begin{figure}[h]
\centering
\begin{minipage}{.5\textwidth}
  \centering
  \includegraphics[width=1\linewidth]{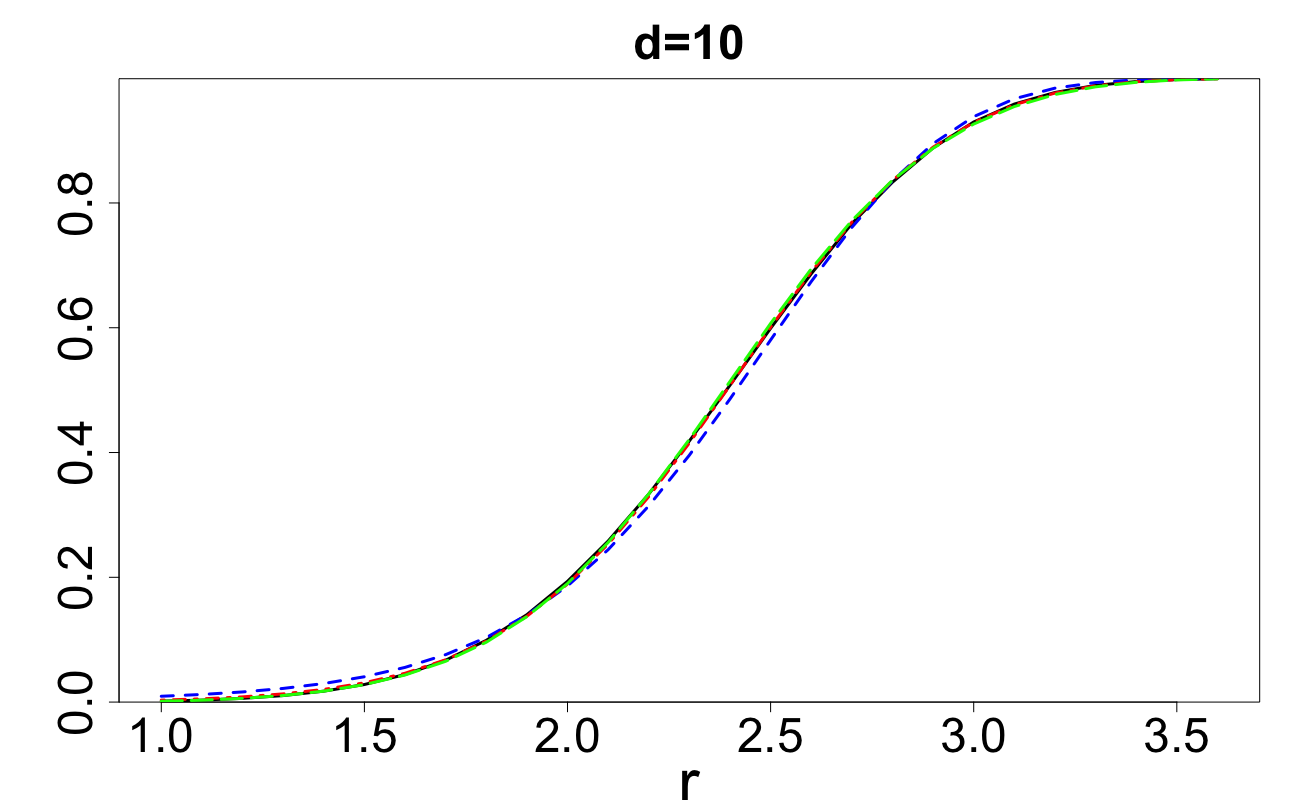}
  \caption{$Z$ is at half-diagonal with $\|Z\|=\frac12 \sqrt{10}$}
\end{minipage}%
\begin{minipage}{.5\textwidth}
  \centering
  \includegraphics[width=1\linewidth]{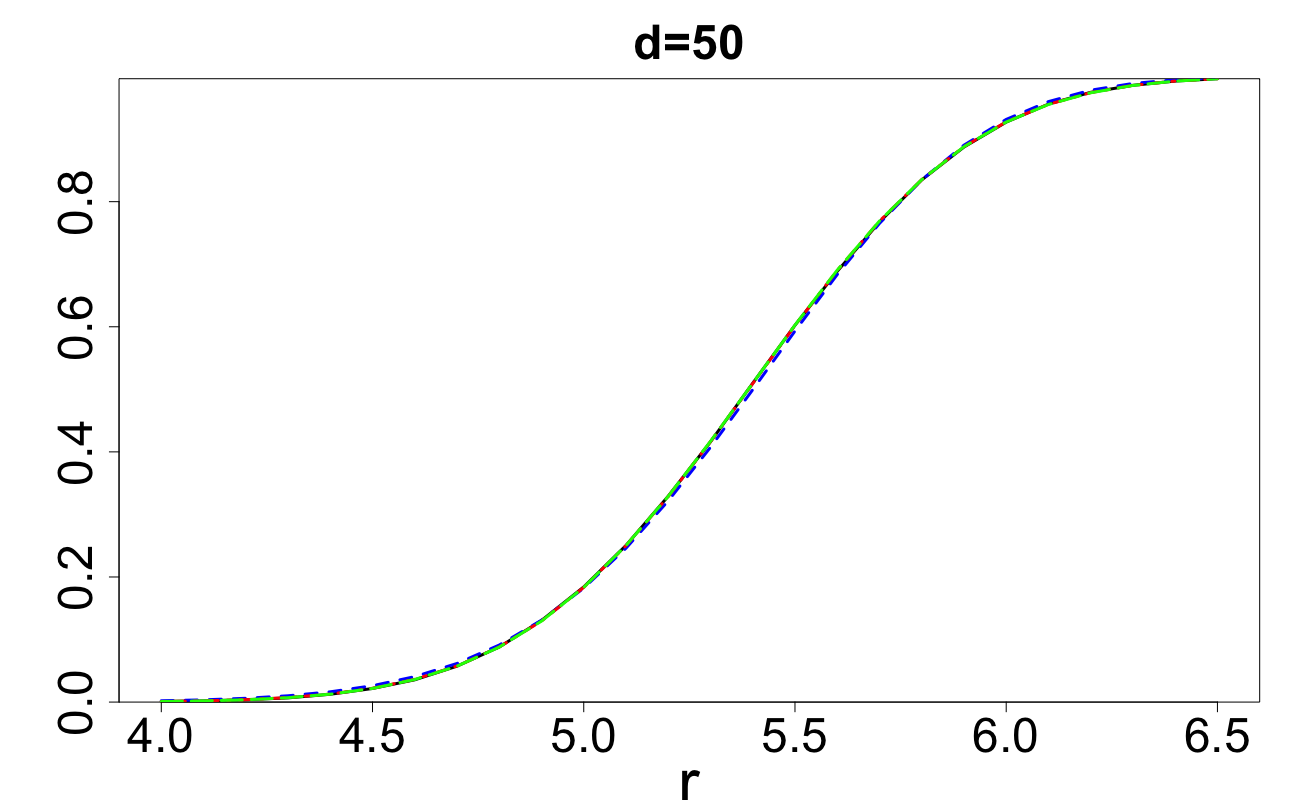}
  \caption{$Z$ is at half-diagonal, $\|Z\|=\frac12 \sqrt{50}$}
\end{minipage}
\end{figure}

\begin{figure}[!h]
\centering
\begin{minipage}{.5\textwidth}
  \centering
  \includegraphics[width=1\linewidth]{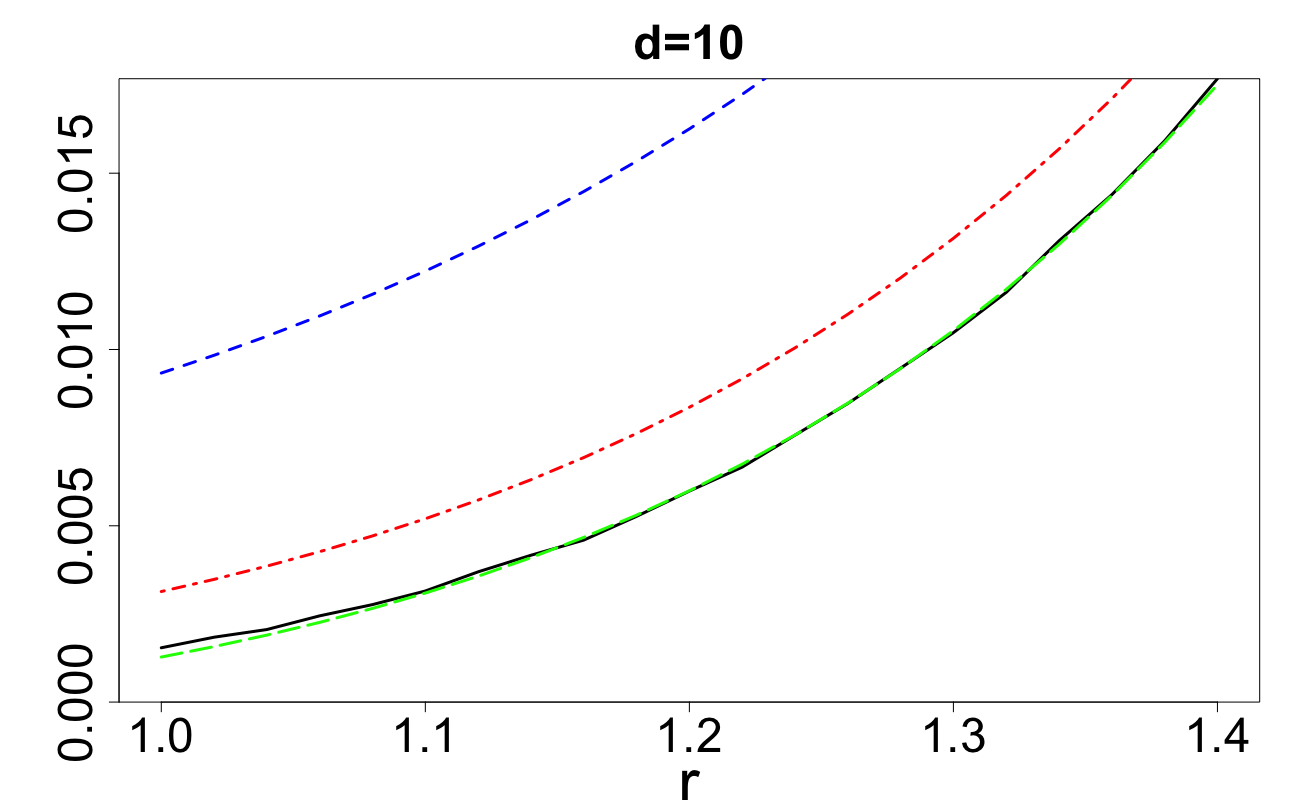}
  \caption{$Z$ is at half-diagonal, $\|Z\|=\frac12 \sqrt{10}$}
\end{minipage}%
\begin{minipage}{.5\textwidth}
  \centering
  \includegraphics[width=1\linewidth]{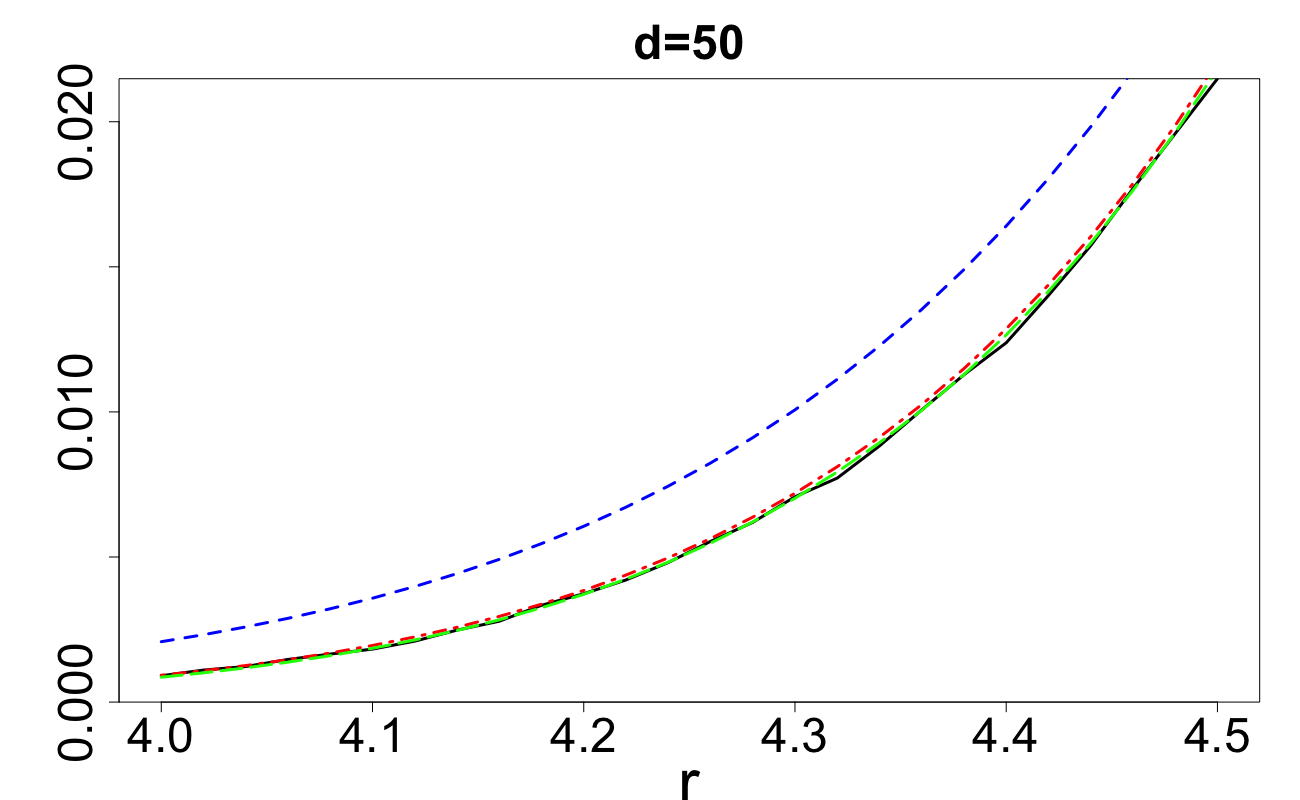}
  \caption{$Z$ is at half-diagonal, $\|Z\|=\frac12 \sqrt{50}$}
\label{z_half_diag_50}
\end{minipage}
\end{figure}

\begin{figure}[h]
\centering
\begin{minipage}{.5\textwidth}
  \centering
  \includegraphics[width=1\linewidth]{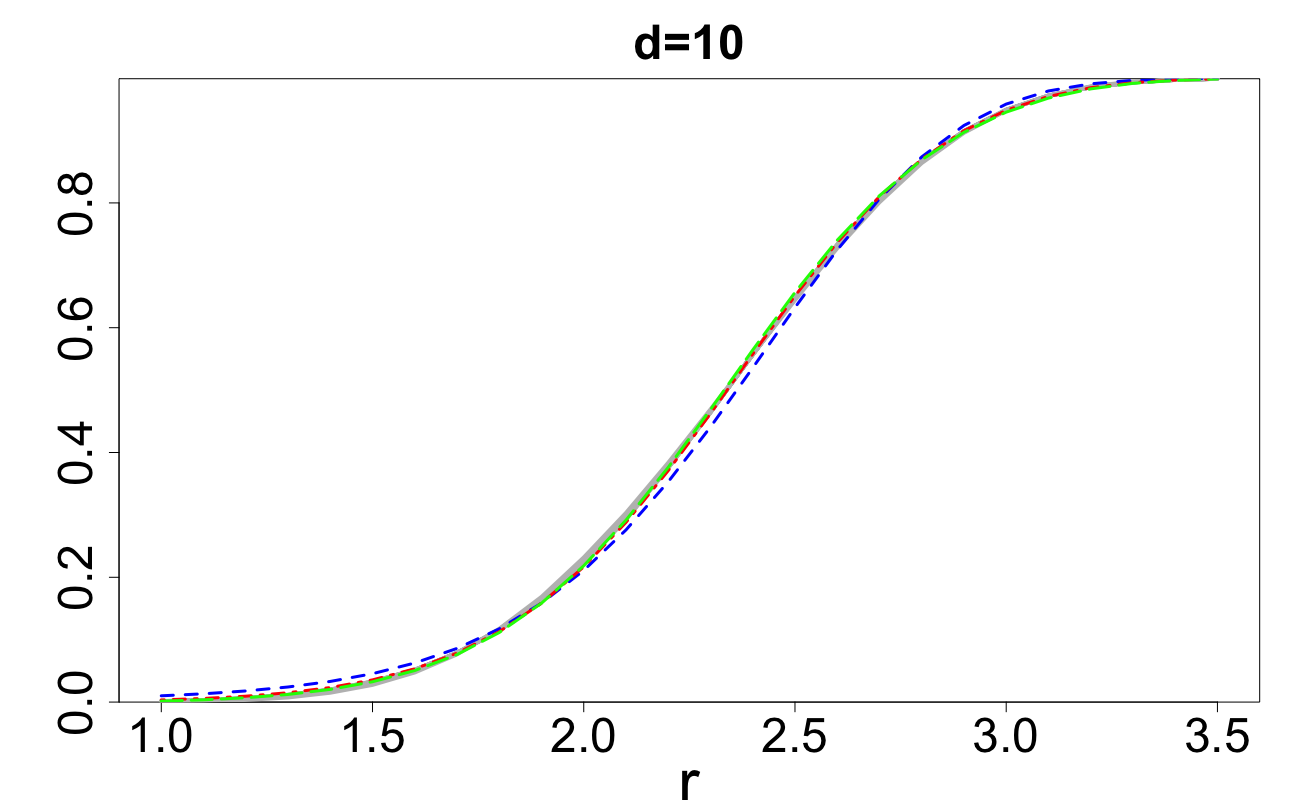}
  \caption{$d=10$, $Z\in{\cal S}_{10}(0,1.5)$, $r \in [1, 3.5]$ }

\end{minipage}%
\begin{minipage}{.5\textwidth}
  \centering
  \includegraphics[width=1\linewidth]{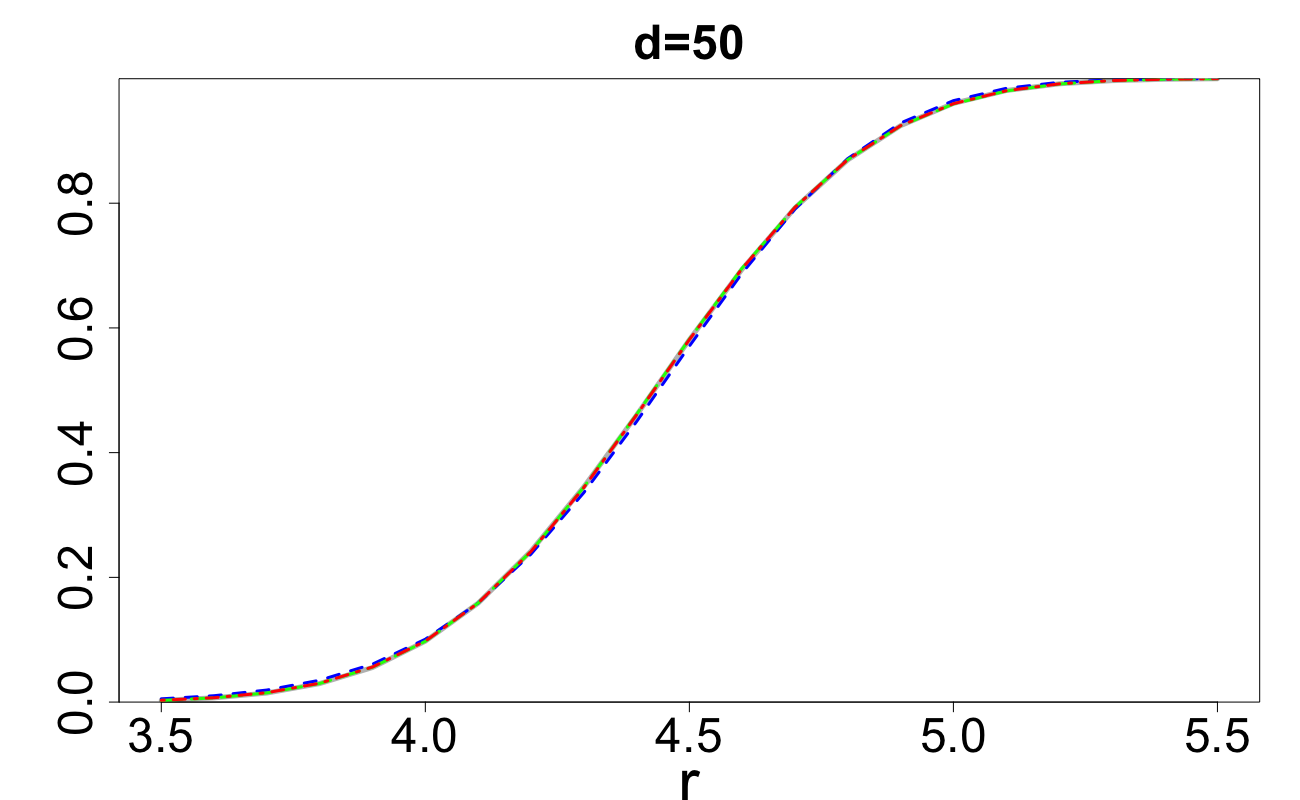}
  \caption{$d=50$, $Z\in{\cal S}_{50}(0,1.75)$, $r \in [3.5, 5.5]$ }
  \label{d_50_on_sphere_non_zoom}
\end{minipage}
\end{figure}

\clearpage
\begin{figure}[h]
\centering
\begin{minipage}{.5\textwidth}
  \centering
  \includegraphics[width=1\linewidth]{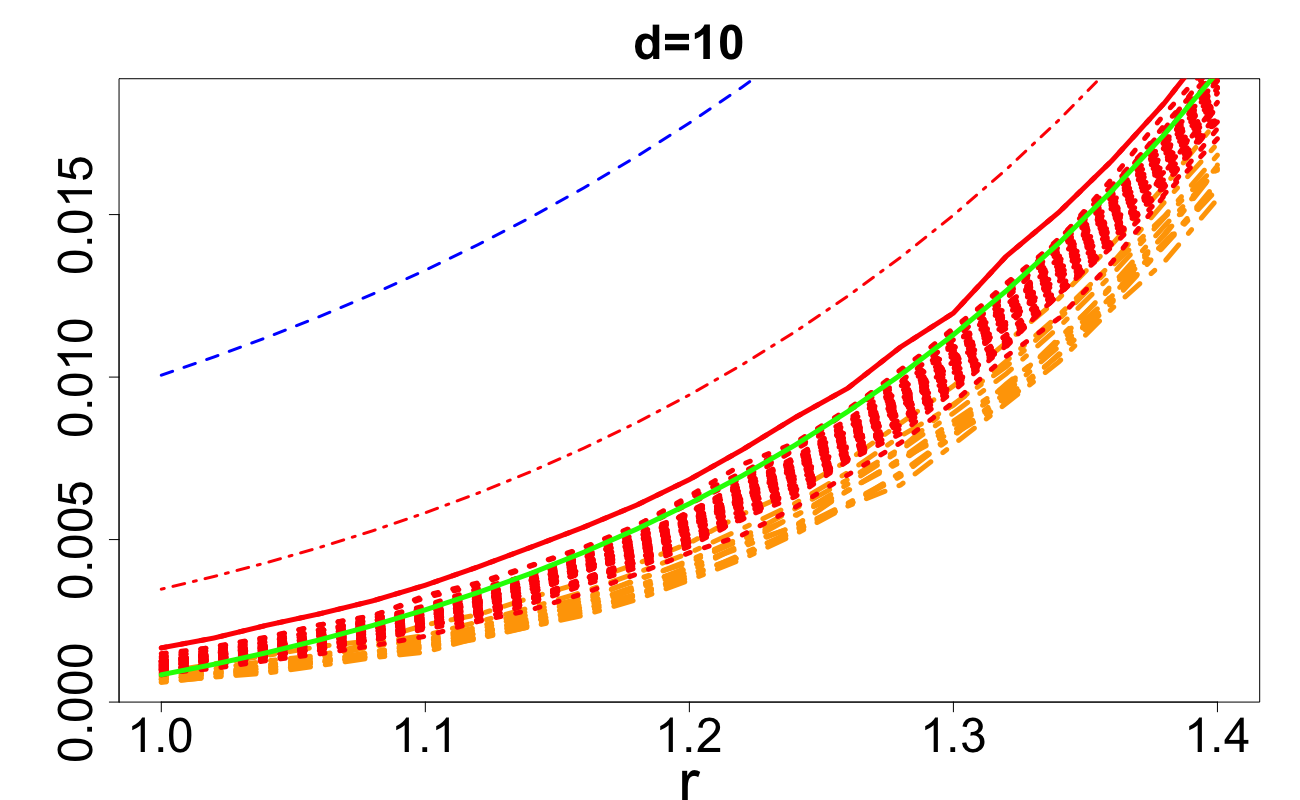}
  \caption{$d=10$, $Z\in{\cal S}_{10}(0,1.5)$, $r \in [1, 1.4]$ }
\label{d_10_on_sphere}
\end{minipage}%
\begin{minipage}{.5\textwidth}
  \centering
  \includegraphics[width=1\linewidth]{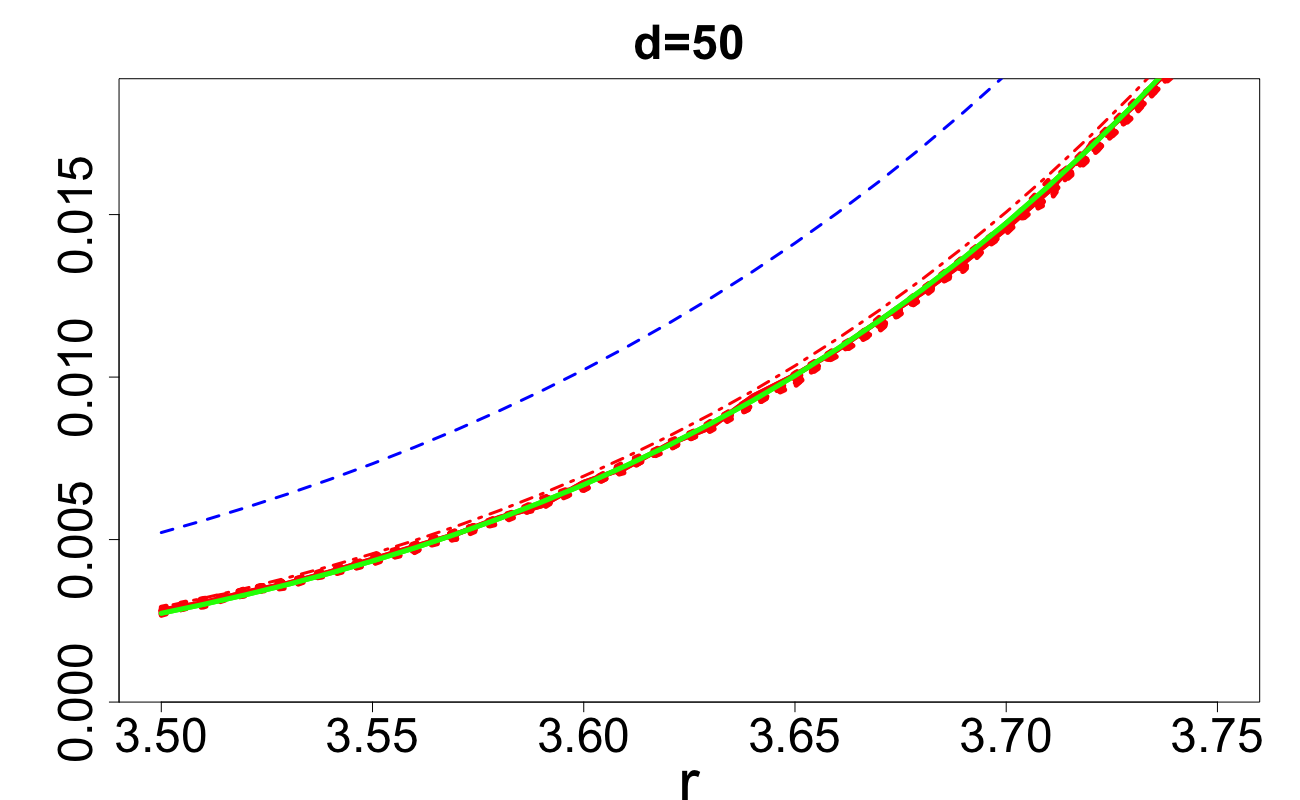}
  \caption{$d=50$, $Z\in{\cal S}_{50}(0,1.75)$, $r \in [3.5, 3.75]$ }
\label{d_50_on_sphere}
\end{minipage}
\end{figure}

%
%
%
%

From the simulations that led to Figures~\ref{d_10_z_0} -- \ref{d_50_on_sphere}  we can make the following conclusions.
\begin{itemize}
  \item The normal approximation \eqref{eq:inters2f} is quite satisfactory unless the value $C_{d,Z,{ r }} $ is small.
  \item The accuracy of all approximations improves as $d$ grows.
  \item The  approximation \eqref{eq:inters2f_corrected_2} is very accurate even if  the values $C_{d,Z,{ r }} $ are very small.
   \item If $d$ is large enough then the  approximations \eqref{eq:inters2f_corrected} and \eqref{eq:inters2f_corrected_2} are practically identical and are extremely accurate.
    \end{itemize}

\section{Covering  a cube by $n$ balls}

\label{sec:main_1}

In this section, we consider   the main problem of covering the cube ${\cal C}_{d}=[-1,1]^d$ by the union of $n$ balls ${\cal B}_d(Z_j,r)$ as formulated in Section~\ref{sec:main_intro}. We will discuss different schemes of choosing the set of ball centers
$\mathbb{Z}_n=\{Z_1, \ldots, Z_n\}$ for given $d$ and $n$. The radius $r$ will then be chosen to achieve the required probability of covering: $C_d(\mathbb{Z}_n,r)\geq 1-\gamma$.
Most of the schemes will involve one or several parameters which we will want to choose in an optimal way.

\subsection{The main covering scheme}

The following  will be our main  scheme for choosing $\mathbb{Z}_n=\{Z_1, \ldots, Z_n\}$.\\

{\bf Scheme 1.} {\it
$Z_1, \ldots, Z_n$ are i.i.d. random vectors uniformly distributed  in the cube  ${\cal C}_{d}(\delta)=[-\delta,\delta]^d$,
where $\delta \in [0,1] $ is a parameter.\\
}

We will formulate several other covering schemes and compare them with Scheme~1.
The reasons why we have chosen Scheme 1 as the main scheme are the following.

\begin{itemize}
  \item It is easier to theoretically investigate than all  other non-trivial schemes.
  \item It includes, as a special case when $\delta=1$, the scheme which is very popular in practice of Monte-Carlo
  \cite{niederreiter1992random} and global random search \cite{zhigljavsky2012theory,zhigljavsky2007stochastic}  and is believed to be rather efficient (this is not true).
  \item Numerical studies provided below show that Scheme 1 with optimal $\delta$ provides coverings which are rather efficient, especially for large $d$; see Section~\ref{Numerical_comparison_of_schemes} for a discussion regarding this issue.
\end{itemize}

\subsection{Theoretical investigation of Scheme 1}

Let $Z_1, \ldots, Z_n$ be i.i.d. random vectors uniformly distributed  in the cube  ${\cal C}_{d}(\delta)$ with $0<\delta\leq 1$.
Then, for given $U = (u_1, \ldots, u_d)^\top \in \mathbb{R}^d$,
\be
\mathbb{P} \left\{ U \in {\cal B}_d(\mathbb{Z}_n,r)  \right\}&=& 1-\prod_{j=1}^n \mathbb{P} \left\{ U \notin {\cal B}_d({Z}_j,r)  \right\} \nonumber \\
&=& 1-\prod_{j=1}^n\left(1-\mathbb{P} \left\{ U \in {\cal B}_d({Z}_j,r)  \right\} \right) \nonumber\\
&=& 1-\bigg(1-\mathbb{P}_{_Z} \left\{ \|U - {Z} \| \leq r \right\} \bigg)^n
\label{eq:prod}
\ee
where $
{\cal B}_d(\mathbb{Z}_n,r)$ is defined in \eqref{eq:cover1a}.
The main
 characteristic of interest $C_d(\mathbb{Z}_n,r)$, defined in \eqref{eq:cover1},
the proportion of the cube   covered by the union of balls ${\cal B}_d(\mathbb{Z}_n,r)$, is simply
\be \label{eq:prod5}
C_d(\mathbb{Z}_n,r) = \mathbb{E}_{_U} \mathbb{P} \left\{ U \in {\cal B}_d(\mathbb{Z}_n,r)  \right\}
\ee

Continuing \eqref{eq:prod}, note that
\be
\label{eq:sss}
\mathbb{P}_{_Z} \left\{ \|U - {Z} \| \leq r \right\} = \mathbb{P}_{_Z} \left\{ \sum_{j=1}^d (z_j-u_j)^2 \leq r^2 \right\} =
C^{(\delta)}_{d,U,{ r }}\, ,
\ee
where $C^{(\delta)}_{d,U,{ r }}$ is defined by the formula \eqref{eq:inters2ab}.
From \eqref{eq:inters2ac} and \eqref{eq:inters2ad} we have
$
\label{eq:inters2ae}
C^{(\delta)}_{d,U,{ r }}=C_{d,U/\delta,{  r/\delta }}
$
where $C_{d, U/\delta,{  r/\delta }}$ is the quantity defined by \eqref{eq:inters2a}. This quantity can be approximated in a number of different ways as shown in Section~\ref{sec:quantuty0}. We will compare \eqref{eq:inters2f}, the simplest of the approximations,  with the approximation given in \eqref{eq:inters2f_corrected_2}. Approximation  \eqref{eq:inters2f} gives
\be\label{eq:inters2fa}
C^{(\delta)}_{d,U,{ r }}=C_{d, U/\delta,{ r/ \delta }} \cong \Phi \left(\frac{{ (r/ \delta) }^2- \| U \|^2/\delta^2  -d/3 }{2 \sqrt{ \|U\|^2/(3\delta^2) + d/45} } \right) \,,
\ee
whereas approximation  \eqref{eq:inters2f_corrected_2} provides
\be\label{eq:inters2fa_corrected}
C^{(\delta)}_{d,U,{ r }} \cong \Phi(t_\delta) + c_d\frac{  \|U\|^2/\delta^2+d/63}{5\sqrt{3} (\|U\|^2/\delta^2+d/15)^{3/2} }(1-t_\delta^2)\phi(t_\delta)\, ,
\ee
with $c_d= 1+ {4}/{d}$  and
\bea
t_\delta =  \frac{(r/\delta)^2- \|U\|^2/\delta^2 -d/3}{2\sqrt{\|U\|^2/(3\delta^2) +d/45 } }  \,.
\eea

From \eqref{eq:moments_hd_u0},
$\mathbb{E} \| U \|^2= d/3$ and $ {\rm var}( \| U \|^2)= {4d}/{45}  .$ Moreover, if $d$ is large enough then $\| U \|^2= \sum_{j=1}^d u_j^2$
is approximately normal.

We shall simplify the expression \eqref{eq:prod} by using the approximation
\be\label{eq:inters2fb}
(1-t)^n \simeq e^{-nt}\, ,
\ee
 which is a good approximation for small values of $t$ and moderate values of $nt$; this agrees with the ranges of $d$, $n$ and $r$ we are interested in.

We can combine the expressions \eqref{eq:prod5} and \eqref{eq:prod} with approximations \eqref{eq:inters2fa},\eqref{eq:inters2fa_corrected} and \eqref{eq:inters2fb} as well as with the normal approximation for the distribution of $\| U \|^2$, to arrive at two final approximations for $C_d(\mathbb{Z}_n,r)$ that differ in complexity.
If the original normal approximation of \eqref{eq:inters2fa} is used then we obtain

\be
C_d(\mathbb{Z}_n,r) \simeq 1- \int_{-\infty}^{\infty} \psi_1(s) \phi(s)d s,   \label{eq:prod5a}\;\;
\ee
with
\bea
   \psi_1(s)=\exp \left\{-n
\Phi(c_s) \right\} \,, \,\,\,\ c_s= \frac{{3  (r/\delta) }^2-  s' -d }{2 \sqrt{ s' + d/5} }\,,  \,\,\,\,\, s'=(d+2s \sqrt{d/5})/\delta^2 \,.
\eea
 If approximation \eqref{eq:inters2fa_corrected} is used, we obtain:
\be
\label{eq:accurate_app}
C_d(\mathbb{Z}_n,r) \simeq 1- \int_{-\infty}^{\infty} \psi_2(s) \phi(s)d s, \;\; \label{eq:prod6a}\;\;
\ee
with
\bea
   \psi_2(s)=\exp \left\{-n \left(
\Phi(c_s) +\left(1+\frac{4}{d}\right) \frac{s'+d/21}{5[ s' +d/5 ]^{3/2}}(1-c_s^2)\phi(c_s)   \right) \right\} \,.
\eea

\subsection{Simulation study for assessing accuracy of approximations \eqref{eq:prod5a} and \eqref{eq:prod6a}}

In  Figures~\ref{d_10_n_128_approx}--\ref{d_80_n_512}, $C_d(\mathbb{Z}_n,r) $ is represented by a solid black line and has been obtained via Monte Carlo methods.
 Approximation~\eqref{eq:prod5a} is indicated by a dashed blue line and approximation~\eqref{eq:prod6a} is represented by long dashed green lines. All figures demonstrate that  approximation~\eqref{eq:prod6a} is extremely accurate across different dimensions and values of $n$. This  approximation is much superior to  approximation~\eqref{eq:prod5a}.

\begin{figure}[h]
\centering
\begin{minipage}{.5\textwidth}
  \centering
  \includegraphics[width=1\linewidth]{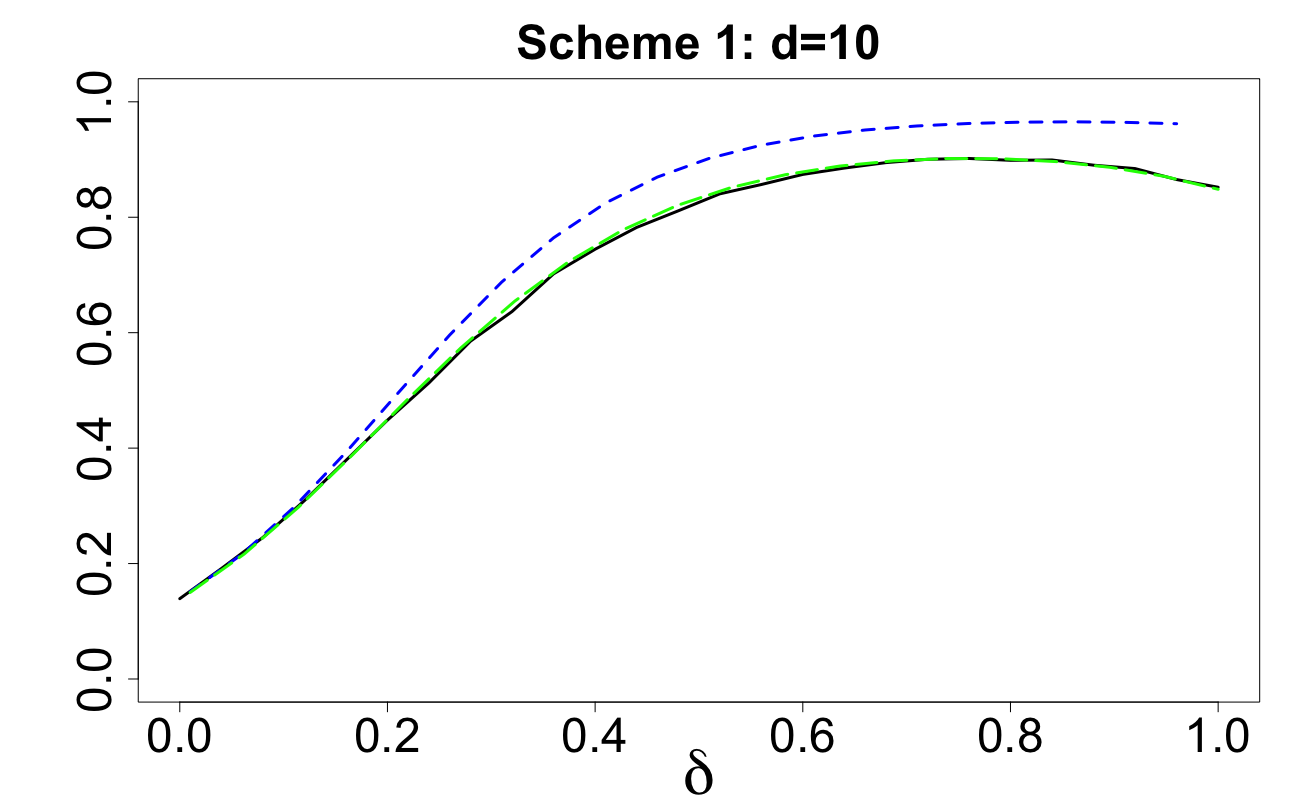}
  \caption{$C_d(\mathbb{Z}_n,r) $ and approximations: $n=128$.}
  \label{d_10_n_128_approx}

\end{minipage}%
\begin{minipage}{.5\textwidth}
  \centering
  \includegraphics[width=1\linewidth]{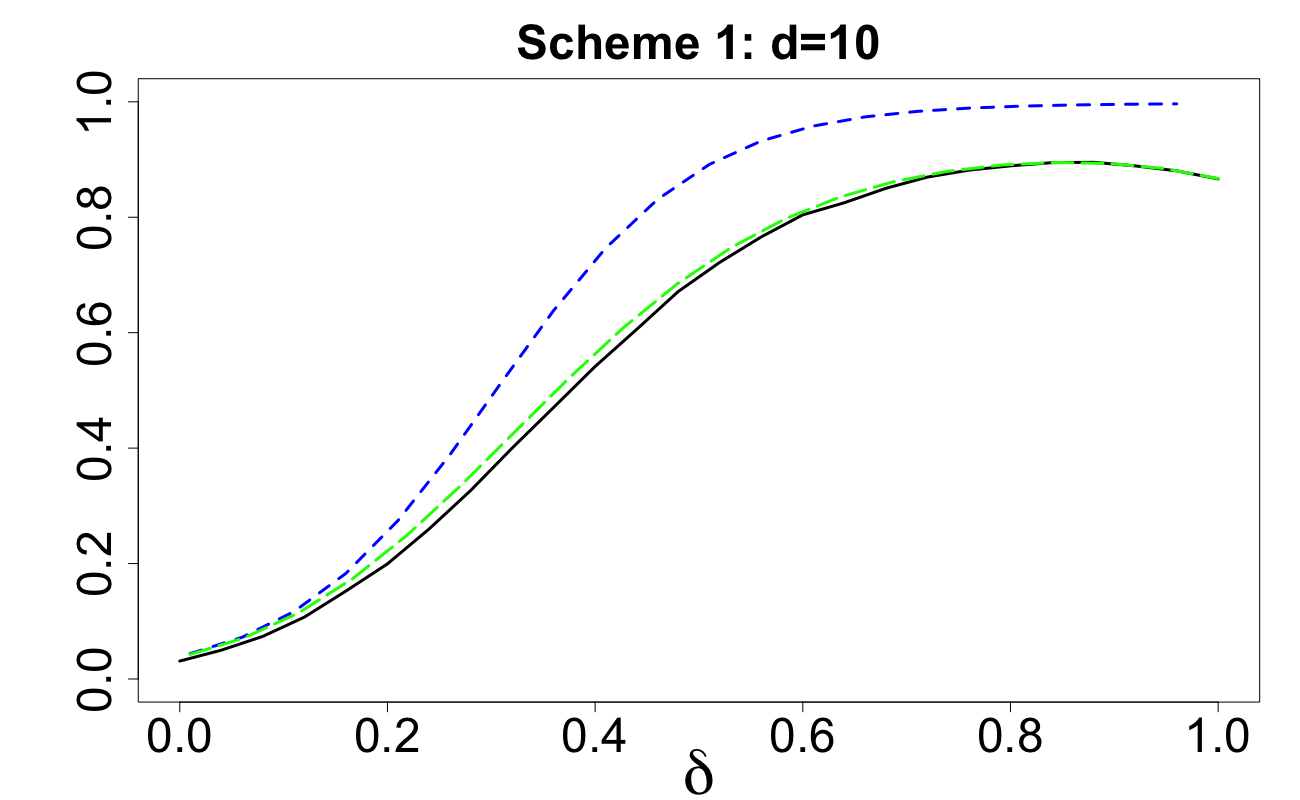}
  \caption{$C_d(\mathbb{Z}_n,r) $ and approximations: $n=512$.}

\end{minipage}
\end{figure}

\begin{figure}[h]
\centering
\begin{minipage}{.5\textwidth}
  \centering
  \includegraphics[width=1\linewidth]{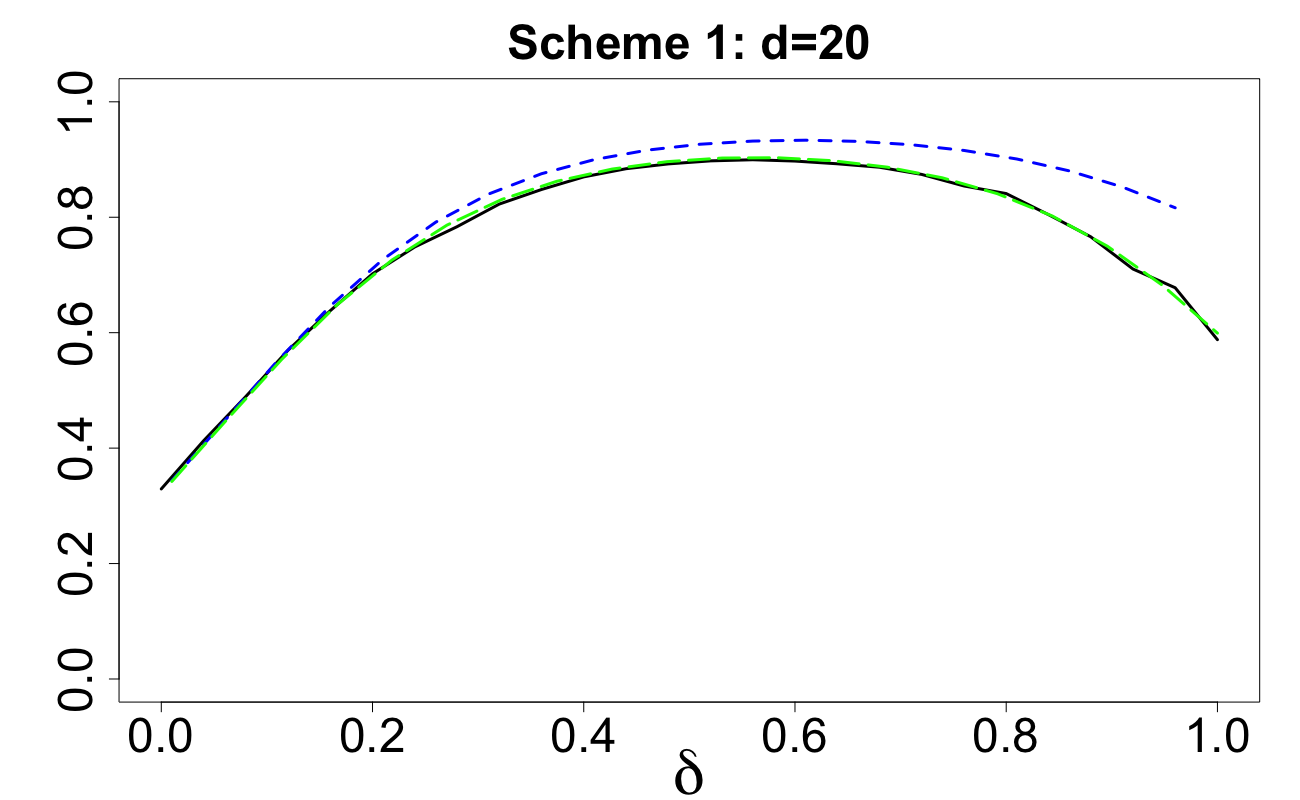}
  \caption{$C_d(\mathbb{Z}_n,r) $ and approximations: $n128$.}

\end{minipage}%
\begin{minipage}{.5\textwidth}
  \centering
  \includegraphics[width=1\linewidth]{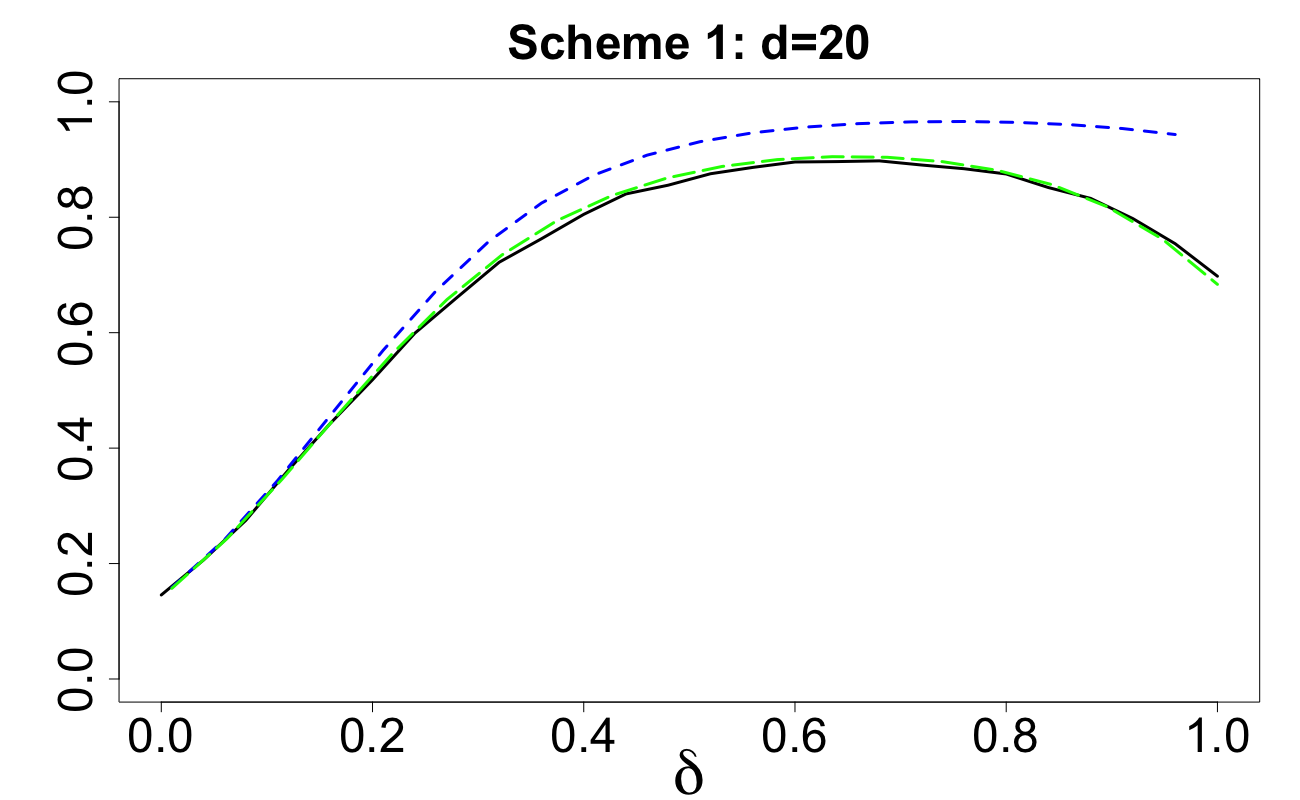}
  \caption{$C_d(\mathbb{Z}_n,r) $ andapproximations: $n=512$.}

\end{minipage}
\end{figure}

\begin{figure}[h]
\centering
\begin{minipage}{.5\textwidth}
  \centering
  \includegraphics[width=1\linewidth]{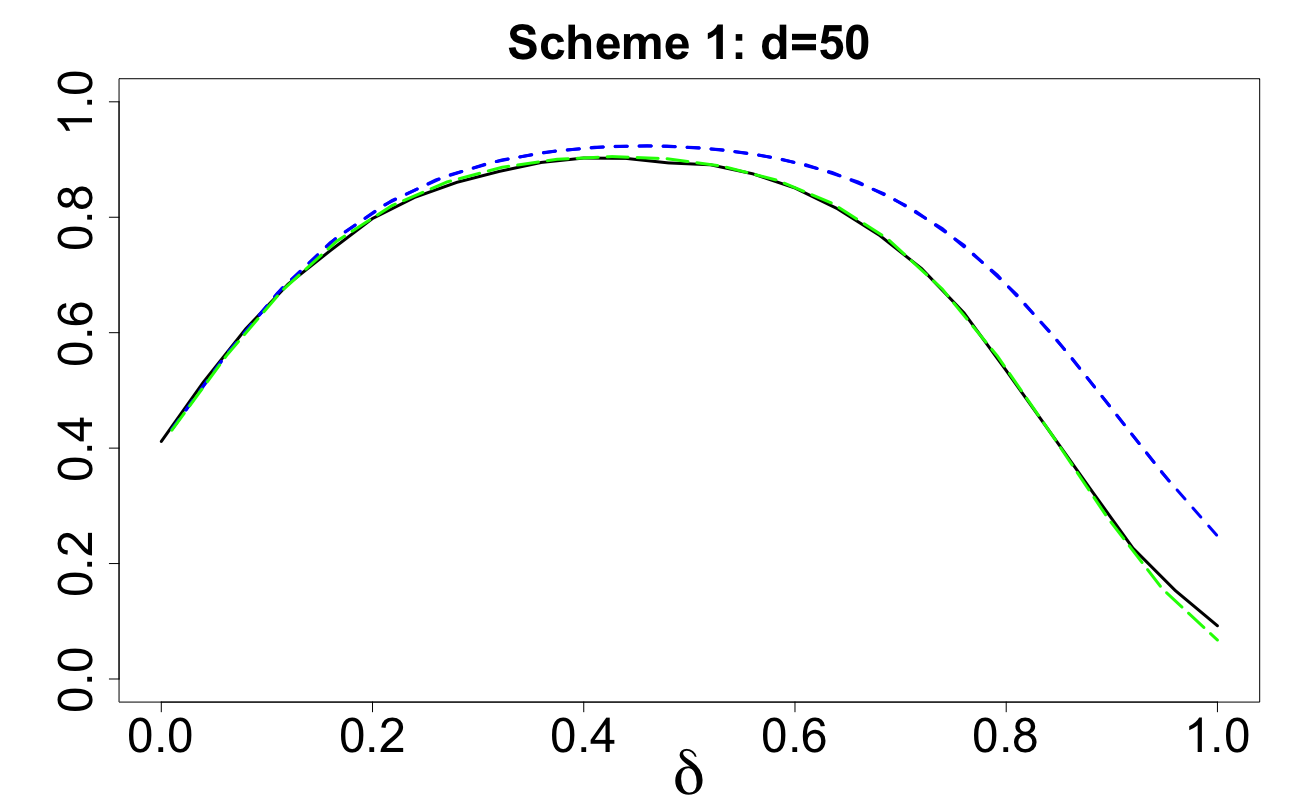}
  \caption{$C_d(\mathbb{Z}_n,r) $ and approximations: $n=512$.}

\end{minipage}%
\begin{minipage}{.5\textwidth}
  \centering
  \includegraphics[width=1\linewidth]{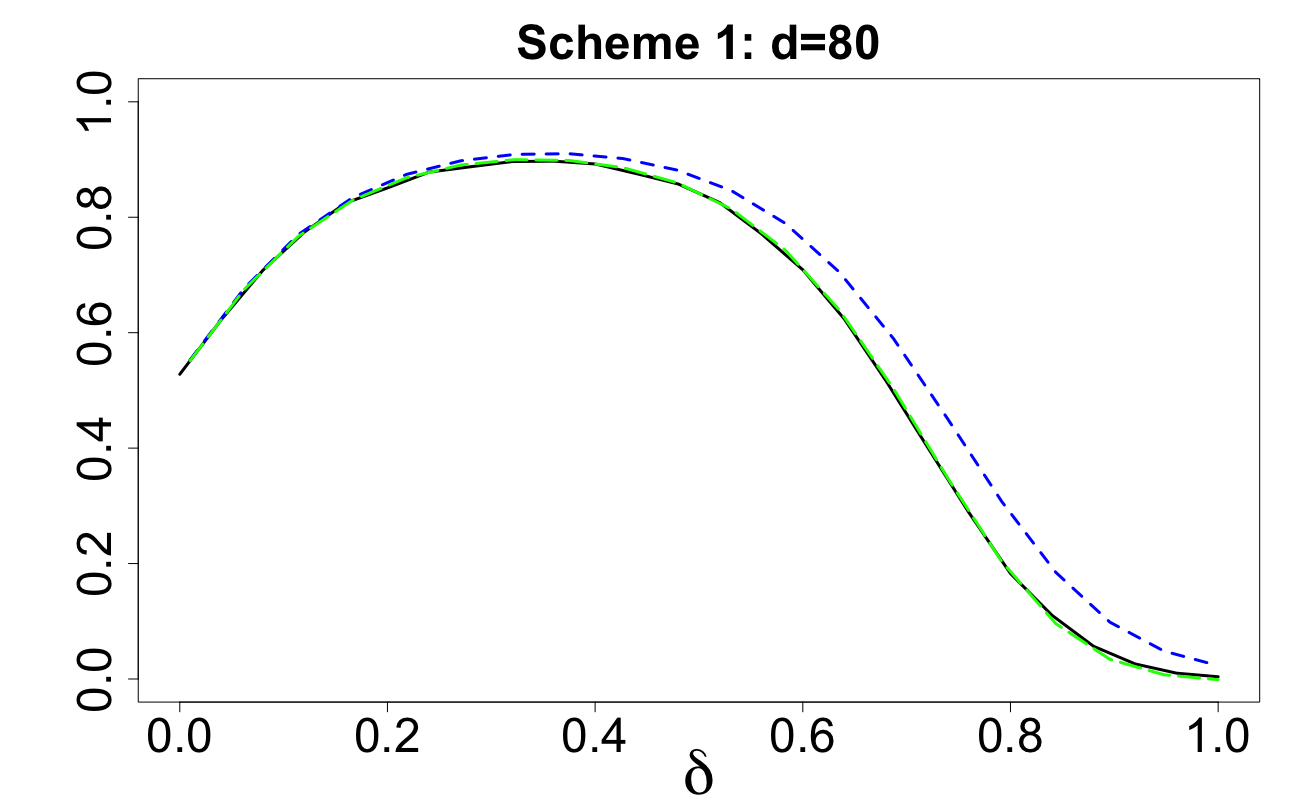}
  \caption{$C_d(\mathbb{Z}_n,r) $ and approximations: $n=512$.}
  \label{d_80_n_512}

\end{minipage}
\end{figure}

\subsection{Other schemes}
\label{sec:other}

In addition to Scheme 1, we have also considered the following schemes for choosing $\mathbb{Z}_n=\{Z_1, \ldots, Z_n\}$.\\

{\bf Scheme 2.} {\it $Z_1=0$;
$Z_2, \ldots, Z_n$ are i.i.d. random vectors uniformly distributed  in the cube  ${\cal C}_{d}(\delta)=[-\delta,\delta]^d$.
}\\

{\bf Scheme 3.} {\it
$Z_1, \ldots, Z_n$ are taken from the minimum-aberration fractional factorial design on vertices of the cube  ${\cal C}_{d}(\delta)=[-\delta,\delta]^d$.
}\\

{\bf Scheme 4.} {\it
$Z_1, \ldots, Z_n$ are i.i.d. random vectors on ${\cal C}_{d}(\delta)$ with independent components  distributed according to Beta-distribution with density \eqref{eq:beta1} with some $\alpha>0$.
}\\

{\bf Scheme 5.} {\it
$Z_1, \ldots, Z_n$ are i.i.d. random vectors uniformly distributed  in the ball  ${\cal B}_{d}(\delta)  $.
}\\

{\bf Scheme 6.} {\it
$Z_1, \ldots, Z_n$ are i.i.d. random vectors uniformly distributed  on the sphere ${\cal S}_{d}(\delta)  $.
}\\

{\bf Scheme 7.} {\it
$Z_1, \ldots, Z_n$ are taken from a low-discrepancy Sobol's sequence on  the cube  ${\cal C}_{d}(\delta)=[-\delta,\delta]^d$.
}\\

The rationale behind the choice of these schemes is as follows. By studying Scheme 2, we test the importance of inclusion of 0 into $\mathbb{Z}_n$. We propositioned that if we included 0 into $\mathbb{Z}_n$, the optimal value of $\delta$ may increase for some of the schemes making them more efficient; this effect has not been detected.

Scheme 3 with optimal $\delta$ is an obvious candidate for being the most efficient.  Unlike all  other  schemes considered, Scheme 3 is only defined for the values of $n$ of the form $n=2^k$ with $k \leq d$.

By using Scheme 4, we test the possibility of improving Scheme 1 by changing the distribution of points in the cube
${\cal C}_{d}(\delta)$. We have found that the effect of distribution is very strong and smaller values of $\alpha$ lead to more efficient covering schemes. By choosing $\alpha$  small enough, like $\alpha=0.1$, we can achieve the average efficiency of covering schemes very close to
the efficiency of Scheme 3. Tables~\ref{table_d_10}--\ref{Table_d_50} contain results obtained for Scheme 4 with $\alpha=0.5$ and $\alpha=1.5$; if $\alpha=1$ then Scheme 4 becomes Scheme 1.

From Section~\ref{sec:concentration}, we know that  for constructing efficient designs we have to  somehow restrict the norms of $Z_j$'s. In Schemes 5 and 6, we are trying to do this in an alternative way to Schemes 1 and 4.

Scheme 7 is a natural improvement of Scheme 1. As a particular case with $\delta=1$, it contains one of the best known low-discrepancy sequences and hence Scheme 7 with $\delta=1$ serves as the main benchmark with which we compare other schemes. For construction, we have used
the R-implementation of   the Sobol's sequences; it is based on \cite{joe2008constructing}.

For all the schemes  excluding Scheme 3,  the sequences $\mathbb{Z}_n=\{Z_1, \ldots, Z_n\}$ are nested so that $\mathbb{Z}_n   \subset \mathbb{Z}_m $ for all $n<m$; using  the terminology of \cite{kuperberg1994line}, these schemes provide on-line coverings
of the cube. Note  that for the chosen values of $n$, Scheme 7 also has some advantage over other schemes considered. Indeed, despite Sobol's sequences are nested, the values $n$ of the form $n=2^k$ are special for the Sobol's sequences and for such values of $n$ the Sobol's sequences possess extra uniformity properties that they do not possess for other values of $n$.

\subsection{Numerical comparison of schemes}
\label{Numerical_comparison_of_schemes}

In Tables~\ref{table_d_10}--\ref{Table_d_50}, for Schemes~1,2,4,5,6 we present the smallest values of $r$ required to achieve an 0.9-coverage on average. For these schemes, the value inside the brackets shows the average value of $\delta$ required to obtain 0.9-coverage. For Schemes 3 and 7, we give the smallest value of $r$ needed for a 0.9-coverage. For these two schemes, the value within the bracket corresponds to the (non random) value of $\delta$ with which we  attain such a coverage.

In Figures~\ref{Scheme_1_multiple_n}--\ref{Scheme_7_multiple_n} we plot $C_d(\mathbb{Z}_n,r) $ as a functions of $\delta \in [0,1]$ across a number schemes, $n$ and $d$. For these plots we have used the values of $r$ provided in Tables~\ref{table_d_10}--\ref{Table_d_50} such that for Figures~\ref{Scheme_1_multiple_n}--\ref{Scheme_2_multiple_n} which correspond to Scheme~1 and Scheme~2, the maximum coverage is very close to $0.9$ and the optimal $\delta$ is very close to the values presented in Tables~\ref{table_d_10}--\ref{Table_d_50}. For Figures~\ref{Scheme_3_multiple_n}--\ref{Scheme_7_multiple_n} the maximum coverage  0.9 is attained with $\delta$ provided in Tables~\ref{table_d_10}--\ref{Table_d_50}. In Figures~\ref{Scheme_1_multiple_n}--\ref{Scheme_7_multiple_n} the solid green line, long dashed red line, dashed blue line and dot dashed orange line correspond to $n=64,128,512,1024$ respectively. The vertical lines on these plots indicate the value of $\delta$ where the maximum coverage is obtained.

\begin{table}[h]
\centering
\begin{tabular}{ |p{2.5cm}||p{2cm}|p{2cm}|p{2cm}|p{2cm}|   }
 \hline
 \multicolumn{5}{|c|}{$d=10$} \\
 \hline
  & $n=64$ &$n=128$& $n=512$ & $n=1024$\\
 \hline
  Scheme 1&  1.632 (0.70)   & 1.520 (0.78) & 1.291 (0.86) &  1.195 (0.90) \\
{   Scheme 1, $\delta=1$}&  1.720 (1.00)   & 1.577 (1.00) & 1.319 (1.00) &1.215 (1.00) \\
Scheme 2&   1.634 (0.70)   &   1.520 (0.78)  & 1.291 (0.86) & 1.195 (0.90) \\
Scheme 3 & 1.530 (0.44) & 1.395 (0.48)  & 1.115 (0.50)   & 1.075 (0.50)  \\
Scheme 4,  $\alpha=0.5$&  1.629 (0.58) &  1.505 (0.65)  & 1.270 (0.72)& 1.165 (0.75)  \\
Scheme 4,  $\alpha=1.5$ & 1.635 (0.80)  & 1.525 (0.88)  &1.310 (1.00) & 1.210 (1.00) \\
Scheme 5& 1.645 (1.40)   & 1.530 (1.50) & 1.330 (1.75) &1.250 (1.75) \\
Scheme 6&  1.642 (1.25)   & 1.532 (1.35) & 1.330 (1.50) & 1.250 (1.70)  \\
Scheme 7&   1.595 (0.72)  & 1.485 (0.80) & 1.280 (0.85) & 1.170 (0.88)   \\
{   Scheme 7, $\delta=1$}& 1.678 (1.00)   &  1.534 (1.00) & 1.305  (1.00) & 1.187  (1.00) \\
 \hline
\end{tabular}
\caption{Values of $r$ and $\delta$ (in brackets) to achieve 0.9 coverage for $d=10$.}
\label{table_d_10}
\end{table}

\begin{table}[h]
\centering
\begin{tabular}{ |p{2.5cm}||p{2cm}|p{2cm}|p{2cm}|p{2cm}|   }
 \hline
 \multicolumn{5}{|c|}{$d=20$} \\
 \hline
  & $n=64$ &$n=128$& $n=512$ & $n=1024$\\
 \hline
  Scheme 1& 2.545 (0.50)    & 2.460 (0.55)  & 2.290 (0.68) &2.205 (0.70)   \\
  {   Scheme 1, $\delta=1$}& 2.840  (1.00)   & 2.702 (1.00) & 2.444 (1.00) & 2.330 (1.00) \\
Scheme 2&   2.545 (0.50)   &  2.460 (0.55)   & 2.290 (0.68)  & 2.205 (0.70)   \\
Scheme 3 & 2.490 (0.32) & 2.410 (0.35) &  2.220 (0.40) &2.125 (0.44)   \\
Scheme 4,  $\alpha=0.5$ & 2.540 (0.44)  & 2.455 (0.48)  &  2.285 (0.55) &2.220 (0.60)  \\
Scheme 4,  $\alpha=1.5$ & 2.545 (0.60) & 2.460 (0.65)  & 2.290 (0.76) &2.215 (0.78)  \\
Scheme 5& 2.550 (1.40)   & 2.467 (1.60) &2.305 (1.75)& 2.235 (1.90)\\
Scheme 6& 2.550 (1.40)     &2.467 (1.58)  &2.305 (1.75) & 2.235 (1.90) \\
Scheme 7&  2.520 (0.50)   & 2.445 (0.60)  &  2.285 (0.68) & 2.196 (0.72)  \\
{   Scheme 7, $\delta=1$}&  2.750 (1.00)   &  2.656 (1.00) &  2.435 (1.00) &  2.325 (1.00) \\
 \hline
\end{tabular}
\caption{Values of $r$ and $\delta$ (in brackets) to achieve 0.9 coverage for $d=20$.}
\end{table}

\begin{table}[h]
\centering
\begin{tabular}{ |p{3cm}||p{2cm}|p{2cm}|p{2cm}|   }
 \hline
 \multicolumn{4}{|c|}{$d=50$} \\
 \hline
  & $n=128$& $n=512$ & $n=1024$\\
 \hline
  Scheme 1     & 4.130 (0.38) & 4.020 (0.45)&3.970 (0.46)   \\
  {   Scheme 1, $\delta=1$}&  4.855  (1.00) & 4.625 (1.00) & 4.520 (1.00) \\
Scheme 2      &  4.130 (0.38)   & 4.020 (0.45)  & 3.970 (0.46) \\
Scheme 3    &  4.110 (0.21)&  4.000  (0.25) &  3.950 (0.28)  \\
Scheme 4  $\alpha=0.5$    & 4.130 (0.30) &4.020 (0.36)  & 3.970 (0.40) \\
Scheme 4 $\alpha=1.5$   & 4.130 (0.42) & 4.020 (0.48)& 3.970 (0.52)  \\
Scheme 5&  4.130 (1.50)   &4.020 (1.75) &3.970 (2.00) \\
Scheme 6&   4.130 (1.50)  & 4.020 (1.75)&3.970 (2.00)   \\
Scheme 7&  4.115 (0.40)   & 4.015 (0.45)  &  3.965 (0.47)   \\
{   Scheme 7, $\delta=1$} & 4.395 (1.00) &  4.379 (1.00) & 4.366  (1.00) \\
 \hline
\end{tabular}
\caption{Values of $r$ and $\delta$ (in brackets) to achieve 0.9 coverage for $d=50$.}
\label{Table_d_50}
\end{table}

\clearpage

\begin{figure}[h]
\centering
\begin{minipage}{.5\textwidth}
  \centering
  \includegraphics[width=1\textwidth]{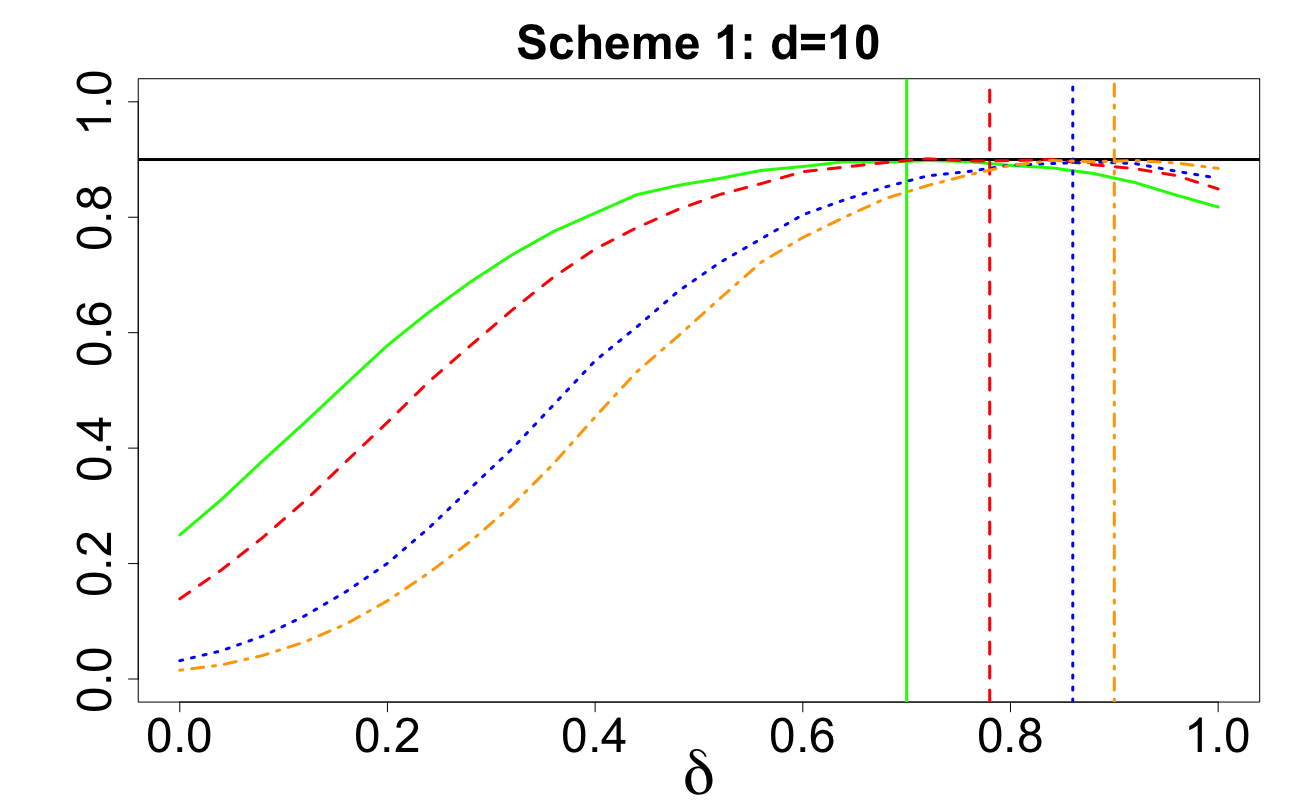}
\caption{Scheme~1: $C_d(\mathbb{Z}_n,r)$ across $\delta$ for $d=10$}
\label{Scheme_1_multiple_n}
\end{minipage}%
\begin{minipage}{.5\textwidth}
  \centering
 \includegraphics[width=1\textwidth]{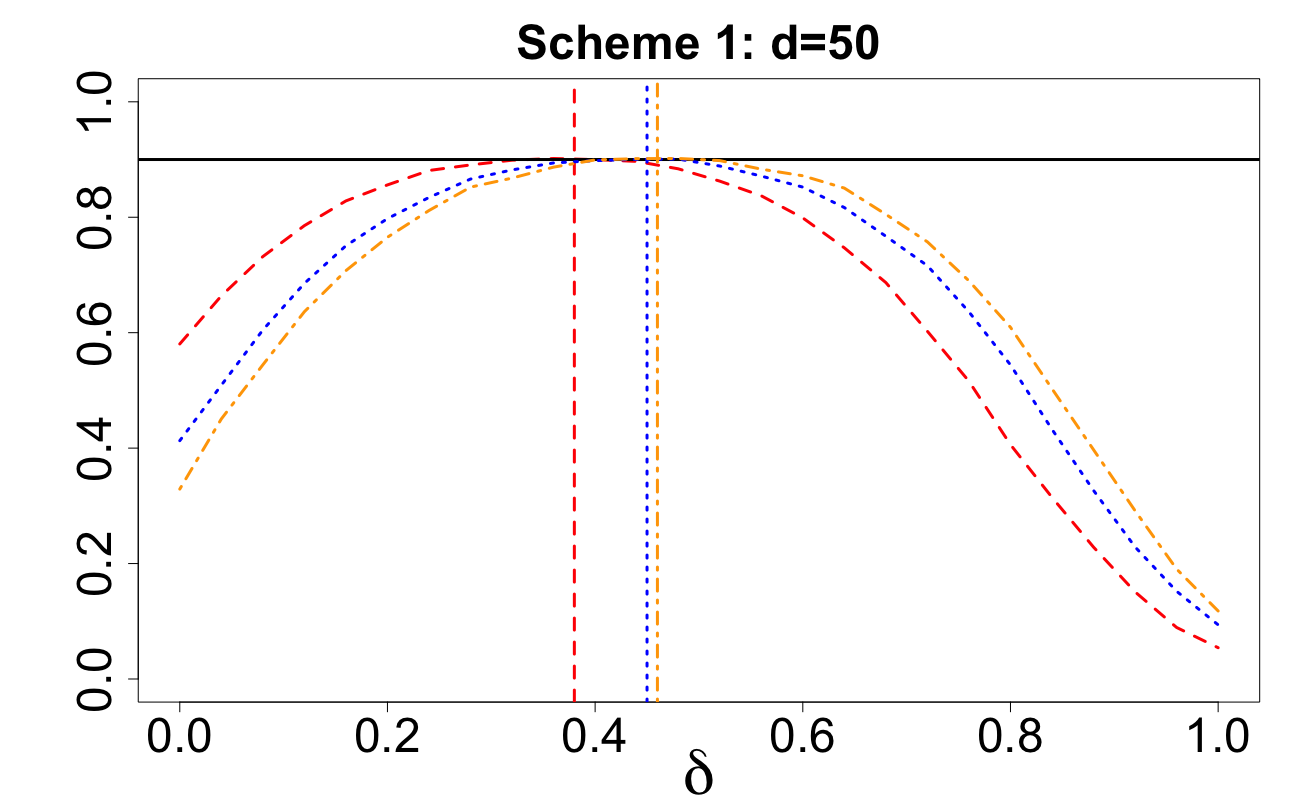}
\caption{Scheme~1: $C_d(\mathbb{Z}_n,r)$ across $\delta$ for $d=50$}

\end{minipage}
\end{figure}

\begin{figure}[h]
\centering
\begin{minipage}{.5\textwidth}
  \centering
  \includegraphics[width=1\textwidth]{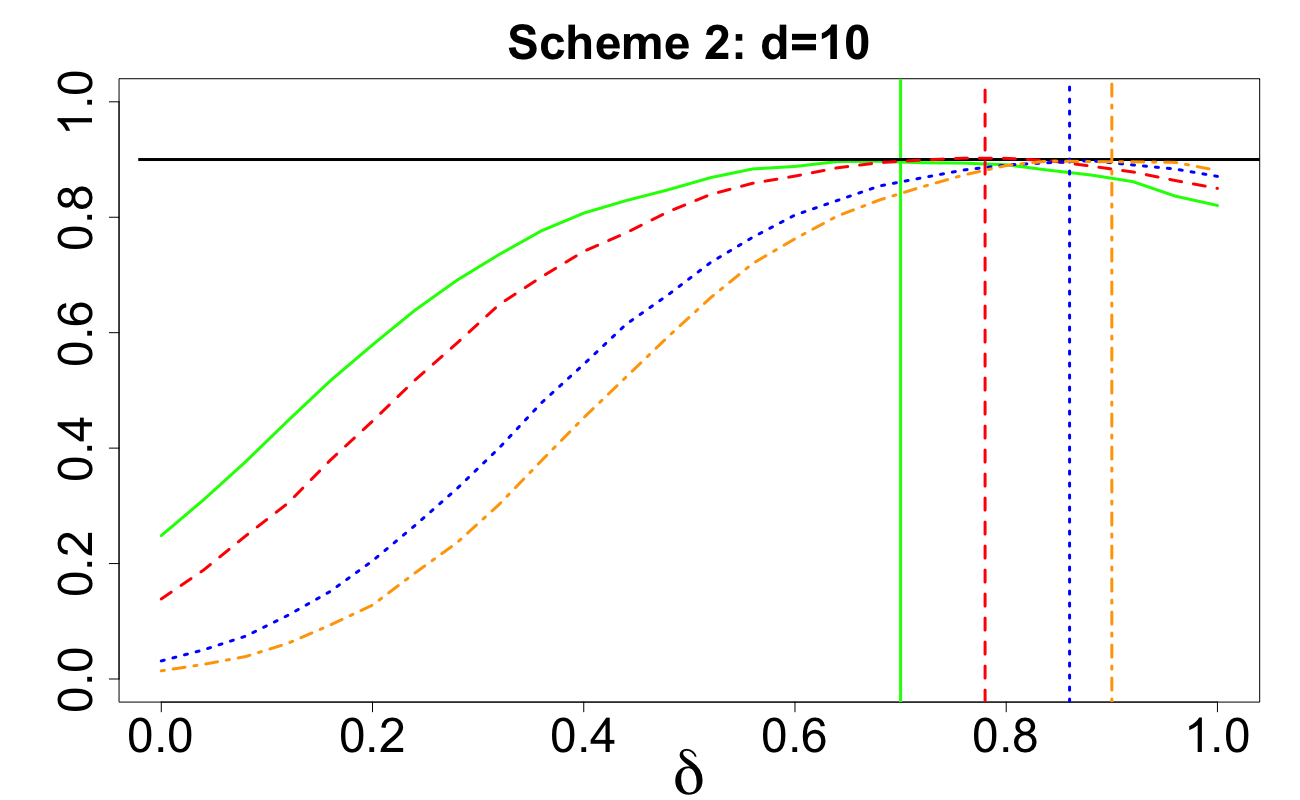}
\caption{Scheme~2: $C_d(\mathbb{Z}_n,r)$ across $\delta$ for $d=10$}
\end{minipage}%
\begin{minipage}{.5\textwidth}
  \centering
 \includegraphics[width=1\textwidth]{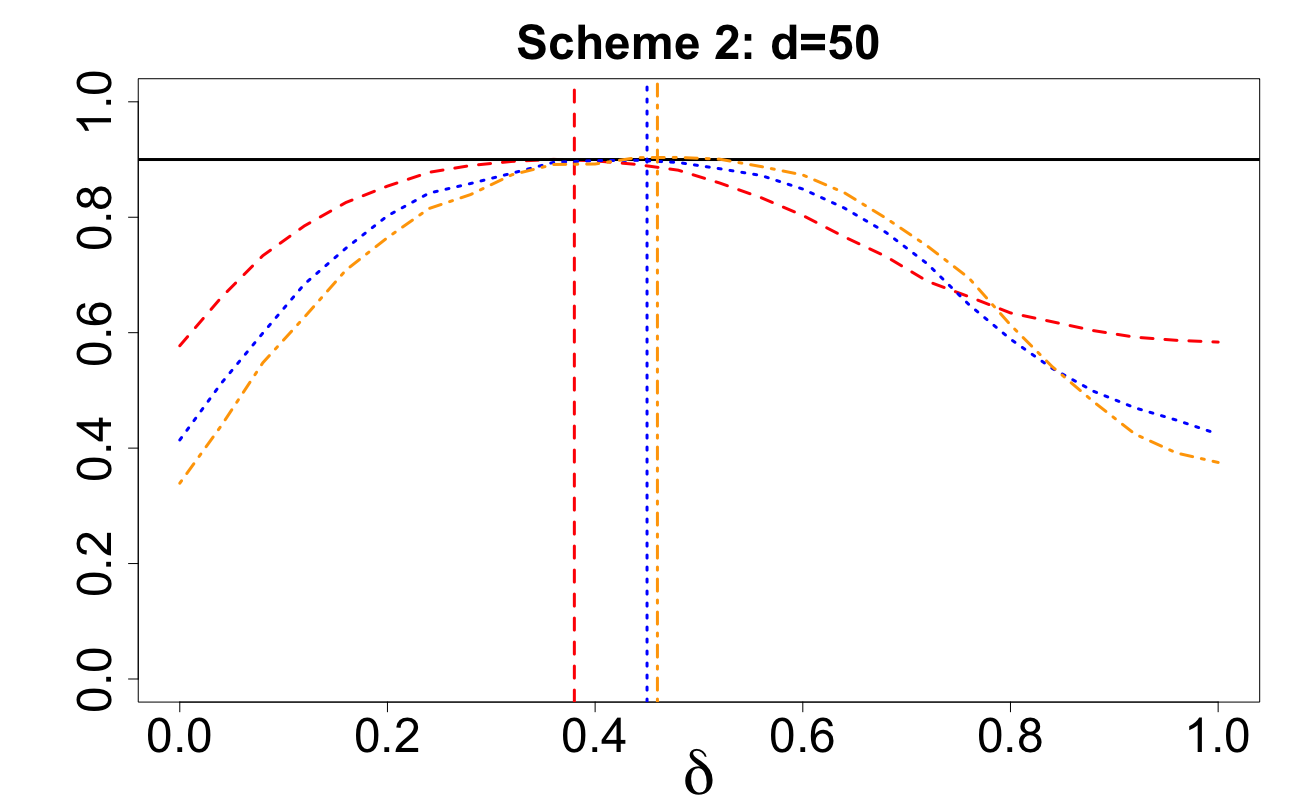}
\caption{Scheme~2: $C_d(\mathbb{Z}_n,r)$ across $\delta$ for $d=50$}
\label{Scheme_2_multiple_n}

\end{minipage}
\end{figure}

\begin{figure}[h]
\centering
\begin{minipage}{.5\textwidth}
  \centering
  \includegraphics[width=1\textwidth]{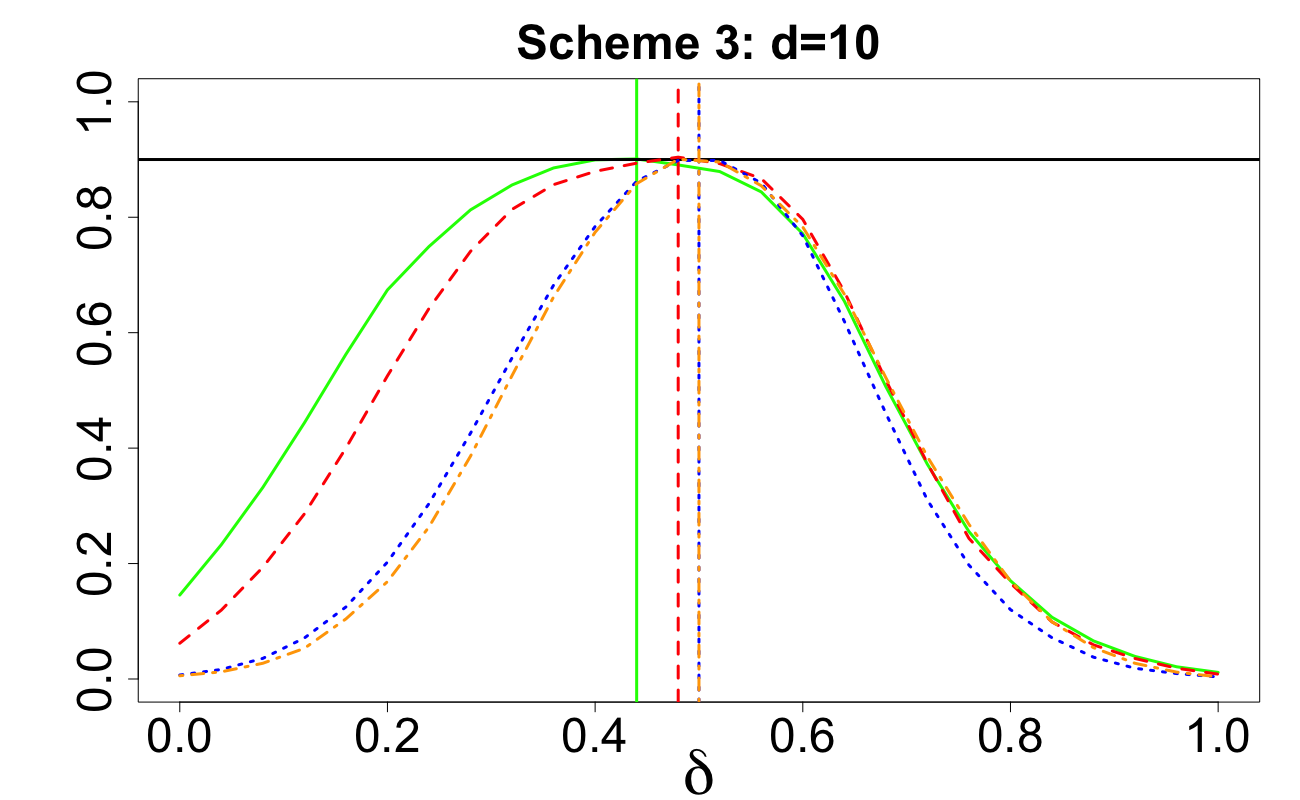}
\caption{Scheme~3: $C_d(\mathbb{Z}_n,r)$ across $\delta$ for $d=10$}
\label{Scheme_3_multiple_n}
\end{minipage}%
\begin{minipage}{.5\textwidth}
  \centering
 \includegraphics[width=1\textwidth]{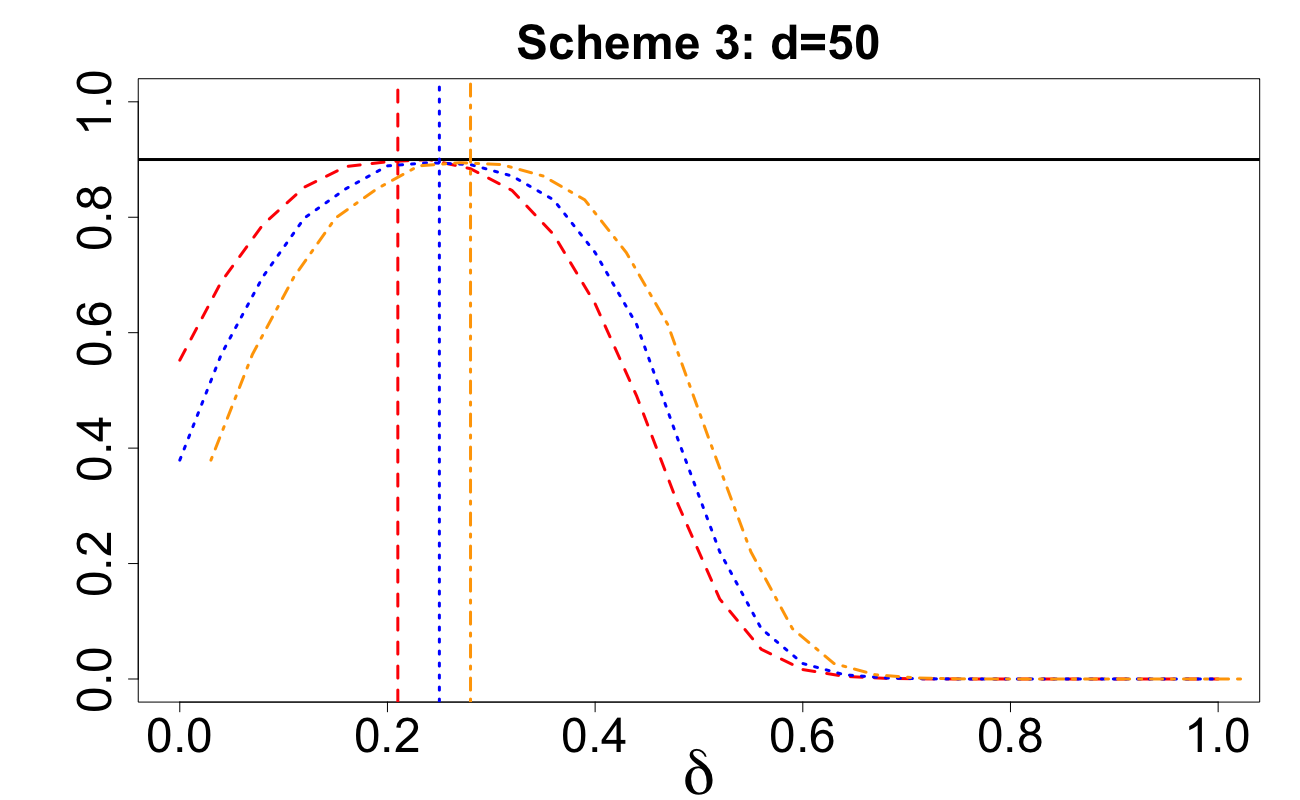}
\caption{Scheme~3: $C_d(\mathbb{Z}_n,r)$ across $\delta$ for $d=50$}
\end{minipage}
\end{figure}

\begin{figure}[h]
\centering
\begin{minipage}{.5\textwidth}
  \centering
  \includegraphics[width=1\textwidth]{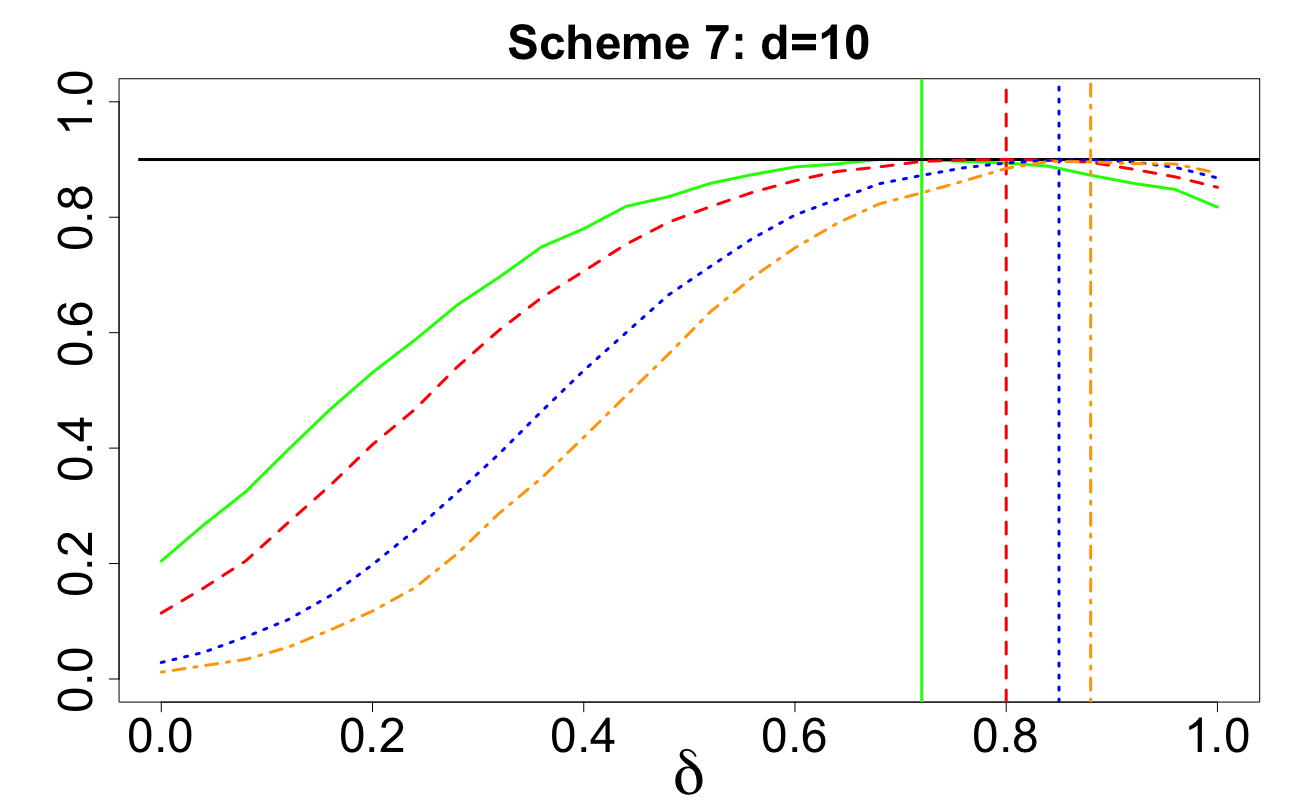}
\caption{Scheme~7: $C_d(\mathbb{Z}_n,r)$ across $\delta$ for $d=10$}
\end{minipage}%
\begin{minipage}{.5\textwidth}
  \centering
 \includegraphics[width=1\textwidth]{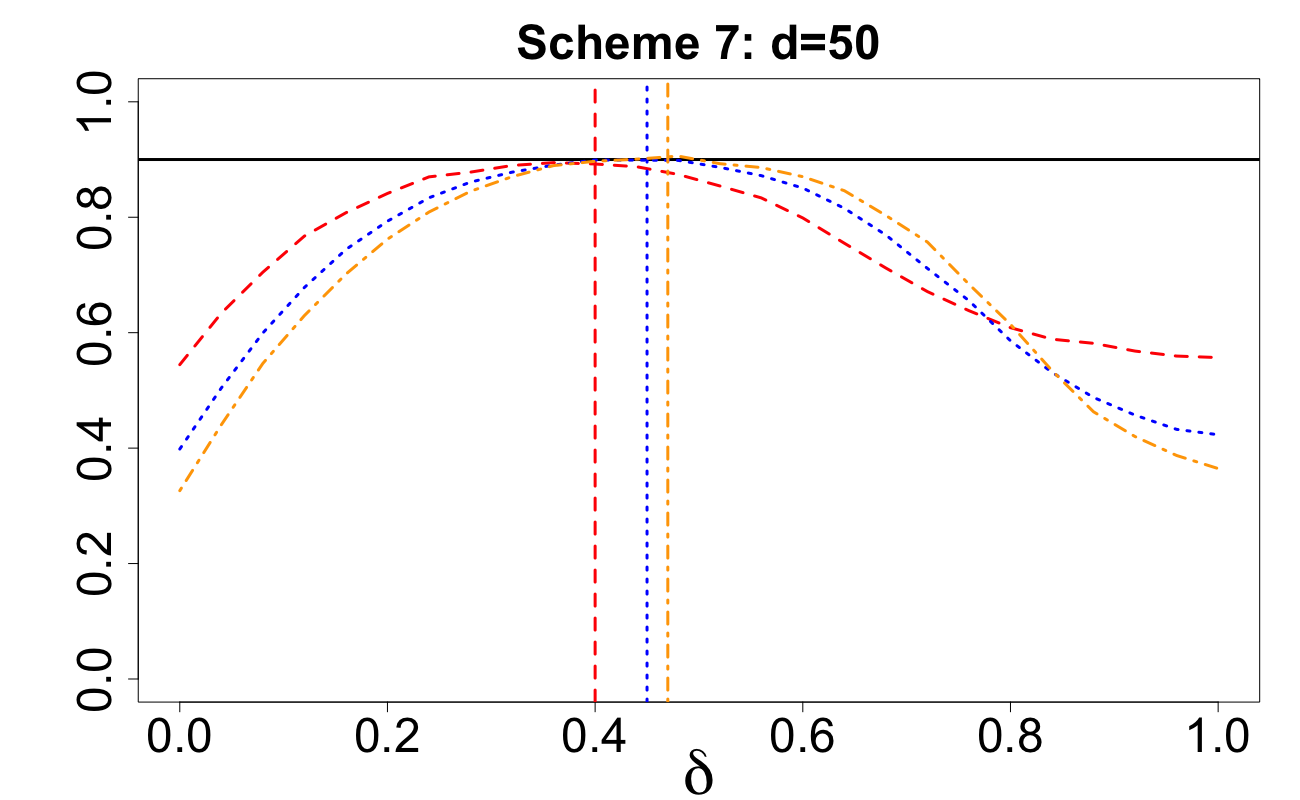}
\caption{Scheme~7: $C_d(\mathbb{Z}_n,r)$ across $\delta$ for $d=50$}
\label{Scheme_7_multiple_n}

\end{minipage}
\end{figure}

From Tables~\ref{table_d_10}--\ref{Table_d_50} and Figures~\ref{Scheme_1_multiple_n}--\ref{Scheme_7_multiple_n} we arrive at the following conclusions:

\begin{itemize}
\item the $\delta$-effect is very important and getting much stronger as $d$ increases;
\item coverage  of unadjusted low-discrepancy sequences is extremely low;
\item properly $\delta$-tuned deterministic Scheme 3, which uses fractional factorial designs of minimum abberation, provides excellent covering;
\item randomized Scheme 4 with suitably chosen parameters of the Beta-distribution, also provides extremely high-quality covering (on average);
\item for all schemes considered, the coverings with the optimal values of $\delta$  fully comply with the result of
Section~\ref{sec:concentration} describing  the area of volume concentration in the cube ${\cal C}_{d}$.
\end{itemize}

%
%

\section{Covering  a cube by cubes}
\label{sec:cubes}

\subsection{Volume of intersection of two cubes}
\label{sec:quantutyA1}

Let us take two cubes:
${\cal C}_{d}=[-1,1]^d$
and
 ${\cal C}_{d}(Z,r)\!=\!\{ Y \in \mathbb{R}^d: \| Y \!-\!Z\|_\infty \! \leq \! { r } \}$,  a  cube of side length $2 r$ centered at a point $Z
 =(z_1, \ldots, z_d)^\top \in {\cal C}_{d}$.
 Denote the fraction of the cube ${\cal C}_{d}$ covered by ${\cal C}_{d}(Z,r)$ by
\be
\label{eq:inters2aaa}
F_{d,Z,{ r }}={{\rm vol}({\cal C}_{d} \cap  {\cal C}_d(Z,{ r }))}/2^d\, .
\ee

Let, like in Section~\ref{sec:quantuty2},  $U=(u_1, \ldots, u_d)^\top$ be a random vector with uniform distribution on ${\cal C}_{d}$ so that $u_1, \ldots, u_d$ are i.i.d.r.v. uniformly distributed on $[-1,1]$. Then
\be
\label{eq:inters2c}
F_{d,Z,{ r }}=  \mathbb{P} \left\{ \| U-Z \|_\infty \leq  { r } \right\}=
 \mathbb{P} \left\{\max_{1\leq j \leq d} |u_j-z_j| \leq r \right\}  \, .\;\;\;\;
\ee
That is,  $F_{d,Z,{ r }}$, as a function of ${ r }$,   is the c.d.f. of the r.v. $\| U-Z \|_\infty = \max_{1\leq j \leq d} |u_j-z_j|  $.

From  Lemma 2 of Section~\ref{sec:appB} the c.d.f. of the r.v. $|u_j-z_j|$ is
\be
\label{eq:inter2b}
G_{d,z_j}(t) = \mathbb{P}\{|u_j\!-\!z_j| \leq t \}= \left\{\begin{array}{cl}
0 & {\rm for\;\;} t\leq 0\\
            t  & {\rm for\;\;} 0<t< 1- |z_j|\\
             \frac12({1\! +\!t\!-\!|z_j|})   & {\rm for\;\;}1- |z_j| \leq t\leq 1+ |z_j| \\
              1 &  1+ |z_j| < t\, .
            \end{array}
            \right.\;\;\;\;\;\;\;
\ee

Since the  c.d.f. of a maximum of independent r.v. is the product of marginal c.d.f.'s, we obtain
$$
F_{d,Z,{ r }}= \prod_{j=1}^d G_{d,z_j}(r)\, .
$$

Two extreme  particular cases of location of $Z$ are:
\begin{itemize}
  \item[(i)] $Z=0$:   $F_{d,0,{ r }}=r^d$, $0 \leq r \leq 1$;
  \item[(ii)] $\|Z\|=\sqrt{d}$, when  $Z$ being a vertex of  the cube ${\cal C}_{d}$: $F_{d,V,{ r }}=(r/2)^d$, $0 \leq r \leq 2$.
%
\end{itemize}

Assume now
that we have the cube ${\cal C}_{d}(\delta)=[-\delta,\delta]^d$ of volume $(2\delta)^d$ and another cube
$
{\cal C}_d(Z^\prime,{ r^\prime })= \{ Y \in \mathbb{R}^d: \| Y-Z^\prime \|_\infty \leq { r^\prime } \}
$ with a center at a point  $Z^\prime=(z_1^\prime, \ldots, z_d^\prime)^\top $.
Denote the fraction of the cube ${\cal C}_{d}(\delta)$ covered by  ${\cal C}_d(Z^\prime,{ r^\prime })$ by
\be
\label{eq:inters2ab1}
F^{(\delta)}_{d,Z^\prime,{ r^\prime }}={{\rm vol}({\cal C}_{d}(\delta) \cap  {\cal C}_d(Z^\prime,{ r^\prime }))}/(2\delta)^d\, .
\ee
Then by  changing  the coordinates and the radius using \eqref{eq:inters2ac}
we get
\be
\label{eq:inters2ad1}
F^{(\delta)}_{d,Z^\prime,{ r^\prime }}=F_{d,Z/\delta,{ r/\delta }} \,.
\ee

\subsection{Proportion  of a cube covered by smaller cubes with random centers}
\label{sec:closed_form_cubes}

Let us take the cube ${\cal C}_{d}=[-1,1]^d$  and $n$ smaller cubes
$
{\cal C}_d(Z_j,{ r })= \{ Y \in \mathbb{R}^d: \| Y-Z_j \|_\infty \leq { r } \}
$ with  centers at  points  $Z_j \in \mathbb{R}^d$.
Denote the fraction of the cube ${\cal C}_{d}$ covered by ${\cal C}_d(\mathbb{Z}_n,{ r })= \cup_{j=1}^n {\cal C}_d(Z_j,{ r })$, the union of these cubes,
 by
\be
\label{eq:inters2aC}
C_{d,\mathbb{Z}_n,{ r }}={{\rm vol}({\cal C}_{d} \cap  {\cal C}_d(\mathbb{Z}_n,{ r }))}/2^d\, .
\ee

Our aim is to obtain a closed form expression for this quantity for arbitrary $d,r$ and~${ n }$ in the case when $Z_1, \ldots, Z_n$ are i.i.d. random vectors uniformly distributed  in the cube  $ {\cal C}_{d}(\delta)=[-\delta,\delta]^d$ with $0<\delta\leq 1$.

Similarly to the combination of \eqref{eq:prod} with \eqref{eq:sss}, for a given $U = (u_1, \ldots, u_d)^\top \in \mathbb{R}^d$,
\be
\mathbb{P} \left\{ U \in {\cal C}_d(\mathbb{Z}_n,r)  \right\}= 1-\bigg(1-F_{d, U/\delta,{  r/\delta }} \bigg)^n\, .
\label{eq:prod7}
\ee
Similarly to
\eqref{eq:prod5},
\be
C_{d,\mathbb{Z}_n,{ r }}  = \mathbb{E}_{_U} \mathbb{P} \left\{ U \in {\cal C}_d(\mathbb{Z}_n,r)  \right\}
= 1- \mathbb{E}_{_U}  \bigg(1-F_{d, U/\delta,{  r/\delta }} \bigg)^n.
\ee

For an integer  $k$, set
\be\label{eq:integral1}
I_k = \frac12 \int_{-1}^{1} \left[G_{d,u/\delta}(r/\delta) \right]^{k}du \,.
\ee
Then, using the binomial theorem, we have
\be\label{eq:cube_approx}
C_{d,\mathbb{Z}_n,{ r }}  = 1 -  \sum_{k=0}^{n}(-1)^{k} \binom{n}{k} I_k^d \,.
\ee
It is possible to evaluate \eqref{eq:integral1} explicity. For $k=0$, we clearly have $I_k=1$. For $k\ge 1$, the integral $I_k $ takes different forms depending on the values of $r$ and $\delta$.
For $k\ge 1$, we have the following:
\begin{itemize}
\item For $r\leq \delta  $:
\bea
I_k = (\delta-r) \left(\frac{r}{\delta}\right)^{k} - \frac{2\delta}{(k+1)} \left\{\left(\frac{\delta+r-1}{2\delta} \right) ^{k+1} - \left(\frac{r}{\delta} \right)^{k+1}  \right \}
\eea
\item For $0 \leq  r - \delta \leq 1$, $r+\delta \geq1$:
\bea
I_k = (r-\delta)  - \frac{2\delta}{(k+1)} \left\{ \left(\frac{\delta+r-1}{2\delta} \right) ^{k+1} - 1  \right \}
\eea
\item For $0 \leq  r - \delta \leq 1$, $r+\delta \leq 1$:
\bea
I_k = (r-\delta)  + \frac{2\delta  }{(k+1)}
\eea
\item For $r - \delta \geq 1$:
\bea
I_k = 1.
\eea
\end{itemize}

In Figures~\ref{Cube_by_cube_d_10}--\ref{Cube_by_cube_d_10_2}, we depict values of $C_{d,\mathbb{Z}_n,{ r }} $ (computed using \eqref{eq:cube_approx}) as a function of $\delta$ for a number of choices of $r$. As in Section~\ref{Numerical_comparison_of_schemes}, we note that the $\delta$-effect holds for the problem of coverage of the cube by smaller cubes.

\begin{figure}[h]
\centering
\begin{minipage}{.5\textwidth}
  \centering
  \includegraphics[width=1\textwidth]{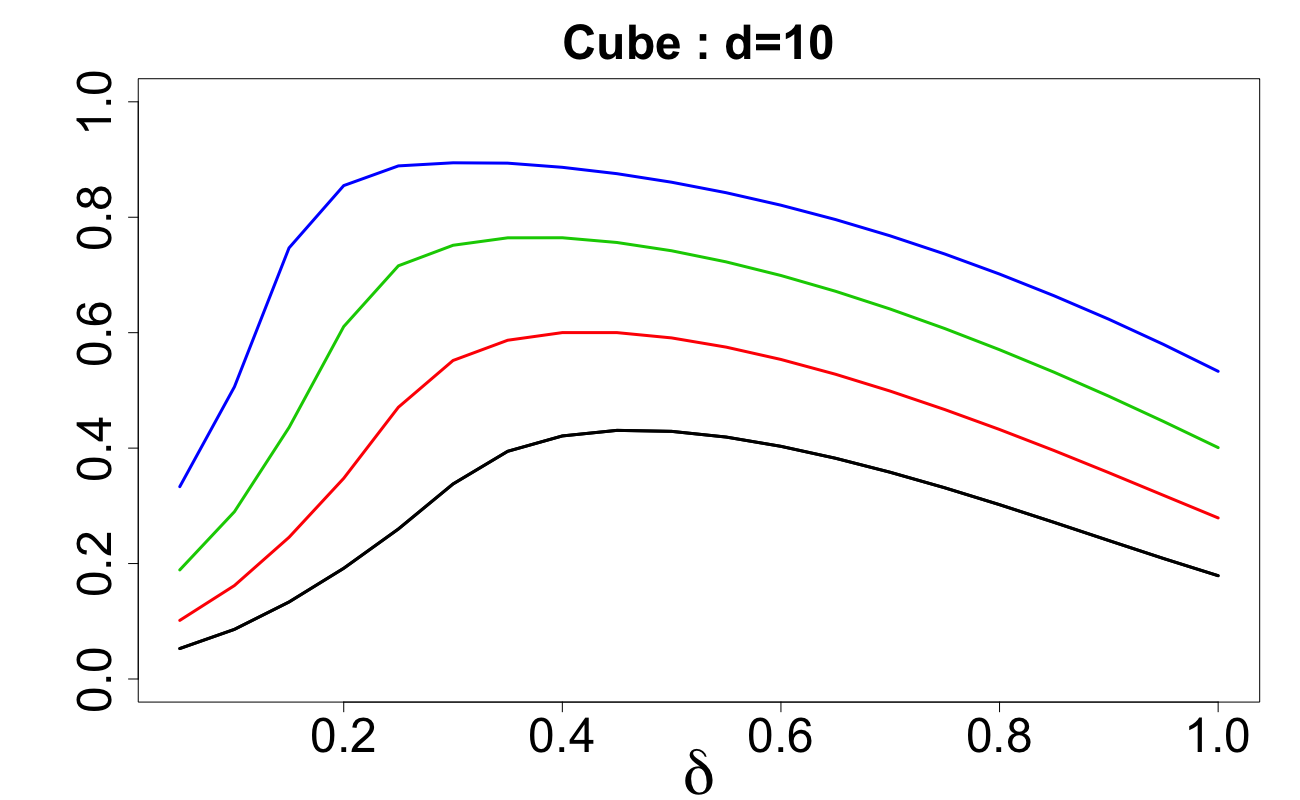}
\caption{$n=50$, $r\in [0.7,0.85]$ increasing by 0.05.}
\label{Cube_by_cube_d_10}
\end{minipage}%
\begin{minipage}{.5\textwidth}
  \centering
 \includegraphics[width=1\textwidth]{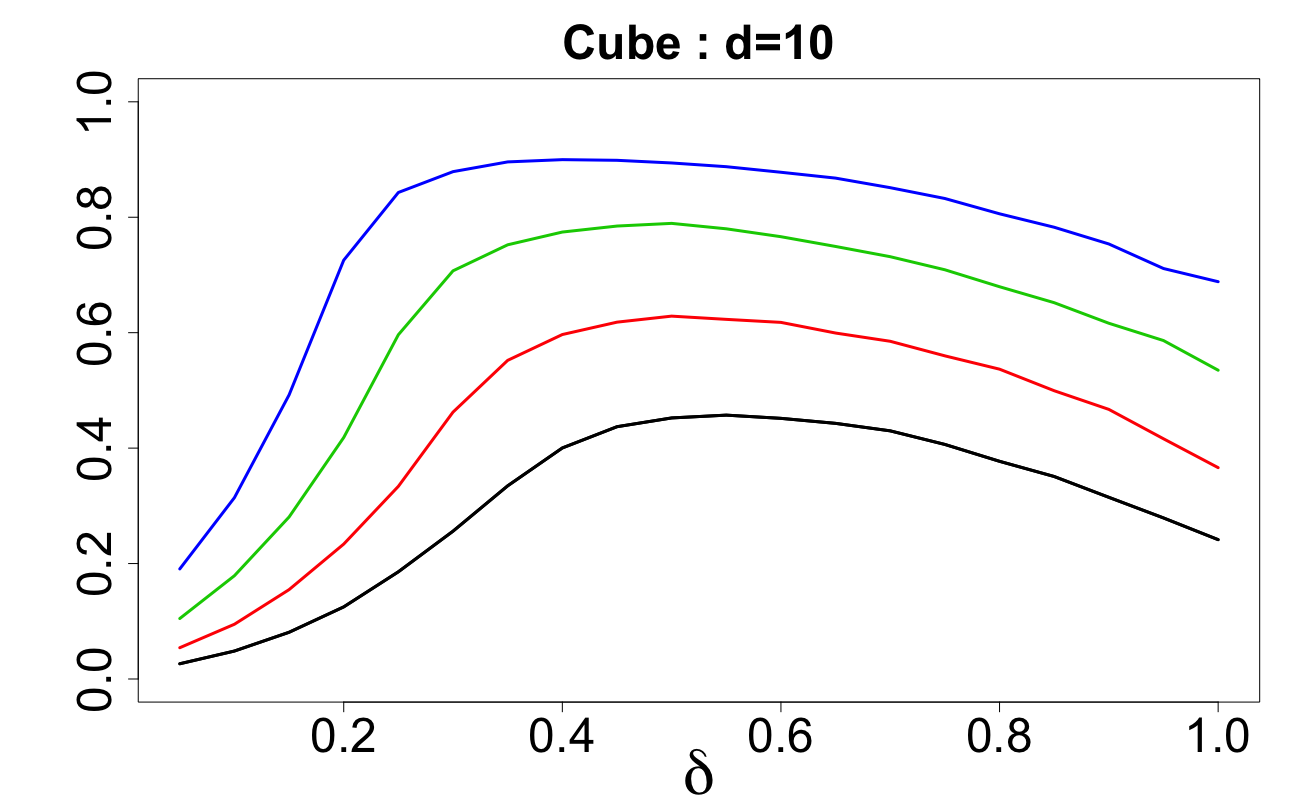}
 \caption{$n=128$, $r\in [0.6,0.8]$ increasing by 0.05.}
\label{Cube_by_cube_d_10_2}
\end{minipage}
\end{figure}

\vspace{-1cm}
\section{Quantization}
\label{sec:quantization}

In this section, we briefly consider the following problem of quantization also known as the problem of facility location.
Let  $X=(x_1, \ldots, x_d)$ be uniform on ${\cal C}_{d}=[-1,1]^d$ and $\mathbb{Z}_n=\{Z_1, \ldots, Z_n\}$ be an $n$-point design.
The mean square quantization error is $\theta_n=\theta(\mathbb{Z}_n)=\mathbb{E}_X\min_{i=1, \ldots, n} \|X-Z_i\|^2$. In the case where $Z_1, \ldots, Z_n$ are i.i.d. uniform on ${\cal C}_{d}(\delta)$, we  will derive a simple approximation for
the expected value of $Q(\mathbb{Z}_n)$ in order to  demonstrate the $\delta$-effect. We shall also notice a strong correlation in design efficiency used for quantization and for  $(1-\gamma)$-covering as studied  in Section~\ref{sec:main_1}.

For deriving an approximation for the quantization mean squared error, we choose Scheme 1 of Section~\ref{sec:main_1}. That is, 
we assume that  $X=(x_1, \ldots, x_d)$ is uniform on ${\cal C}_{d}=[-1,1]^d$ and  $Z_1, \ldots, Z_n$ are i.i.d. uniform on a potentially smaller cube ${\cal C}_{d}(\delta)=[-\delta,\delta]^d$, $0 < \delta \leq 1$. We are interested in finding $\delta=\delta(n,d)$ such that the probability
$\mathbb{P} \{\theta_n \leq r^2\}$ is maximal, where $\theta_n =\min_{i=1, \ldots, n} \|X-Z_i\|^2$.
From Lemma 1, we have for given $X$:
\bea
\mathbb{E}_{_Z} \theta_1=  \mathbb{E}_{_Z}\|X-Z_1\|^2 &=&  \frac{1}{2\delta} \sum_{i=1}^d \int_{-\delta}^{\delta} (x_i-z)^2 dz=\frac{1}{3} \delta^2 d +\|X\|^2 \overset{{\rm def}}{=}\mu_{X,\delta}\, ,
\eea
\bea
{\rm var}( \theta_1)& =&\sum_{i=1}^d {\rm var}_{z}(x_i-z)^2
=\frac{4\delta^2}3 \left[ \frac{\delta^2 d}{15}+ \|X\|^2 \right] \overset{{\rm def}}{=}\sigma^2_{X,\delta}
\eea

We have $\theta_n =\min_{i=1, \ldots, n} \xi_i$, where $\xi_i$ are i.i.d.r.v. with the same distribution as $\theta_1$.
Since $\theta_1$ is a sum of $d$ i.i.d.r.v., for large $d$ we can assume that $\theta_1$ is approximately normal; that is, $\theta_1 \sim N(\mu_{X,\delta},\sigma^2_{X,\delta})$. Under this assumption,
$\xi_i \sim N(\mu_{X,\delta},\sigma^2_{X,\delta})$ and for the conditional expectation of $\theta_n =\min_{i=1, \ldots, n} \xi_i$ we have
\bea
\mathbb{E}_Z \theta_n= \mu_{X,\delta} - \sigma_{X,\delta} E_{n}\, ,
\eea
where $E_n$ is the expectation of the maximum of $n$ i.i.d. $N(0,1)$ r.v.
For any $X \in [0,1]^d$, both $\mu_{X,\delta} $  and $\sigma_{X,\delta}$ increase as $\delta \to 0$ and therefore the behaviour of
$\mathbb{E}\theta_n$ is not obvious when $\delta $ is small.

Since $\mathbb{E} \|X\|^2=d/3$,
\bea
\mathbb{E}_{_X}\mu_{X,\delta}=\frac13{d(1+\delta^2)}\,  \;{\rm and}\;\;
\mathbb{E}_{_X} \sigma^2_{X,\delta}=
\frac{4d\delta^2}9 \left[1+ \frac{\delta^2 }{5} \right] \,.
\eea

%

A rough estimator for $E_n$ is $E_n \simeq \sqrt{2 \log n}$ and (not so rough) estimator for $\mathbb{E}_X \sigma_{X,\eps}$ is
$\mathbb{E}_X \sigma_{X,\eps} \simeq \sqrt{ \mathbb{E}_X \sigma^2_{X,\eps} }.$ This gives
\bea
\mathbb{E} \theta_n &=& \mathbb{E}_{_X} [ \mathbb{E}_Z \theta_n]= \mathbb{E}_X \mu_{X,\delta} - [\mathbb{E}_X \sigma_{X,\delta}] E_{n} \simeq \frac13 { \sqrt{d}} F_{d,n}(\delta)\,,
\eea
where
\be\label{eq:quantization_arrox}
F_{d,n}(\delta)=  \sqrt{d}(1+\delta^2) -  2 \delta \sqrt{ 1+\delta^2/5} \cdot \sqrt{2 \log n} \,.
\ee

We suggest a simple modification to formula \eqref{eq:quantization_arrox} to improve its accuracy for relatively small $d$. That is, we propose using
\be\label{quant_approx_corrected}
\mathbb{E} \theta_n  \simeq \frac13 { \sqrt{d}} \hat{F}_{d,n}(\delta)\,,
\ee
where
\bea
\hat{F}_{d,n}(\delta)=  \sqrt{d}(1+\delta^2) -  \frac85 \delta \sqrt{ 1+\delta^2/5} \cdot \sqrt{2 \log n} \,.
\eea

In Figures~\ref{quant_approx_128}--\ref{quant_approx_512}, we asses the accuracy of the approximation \eqref{quant_approx_corrected}. In these two figures, the solid black line corresponds to $\mathbb{E} \theta_n $ obtained via Monte Carlo methods and the dashed red line depicts the approximation.

\begin{figure}[h]
\centering
\begin{minipage}{.5\textwidth}
  \centering
  \includegraphics[width=1\textwidth]{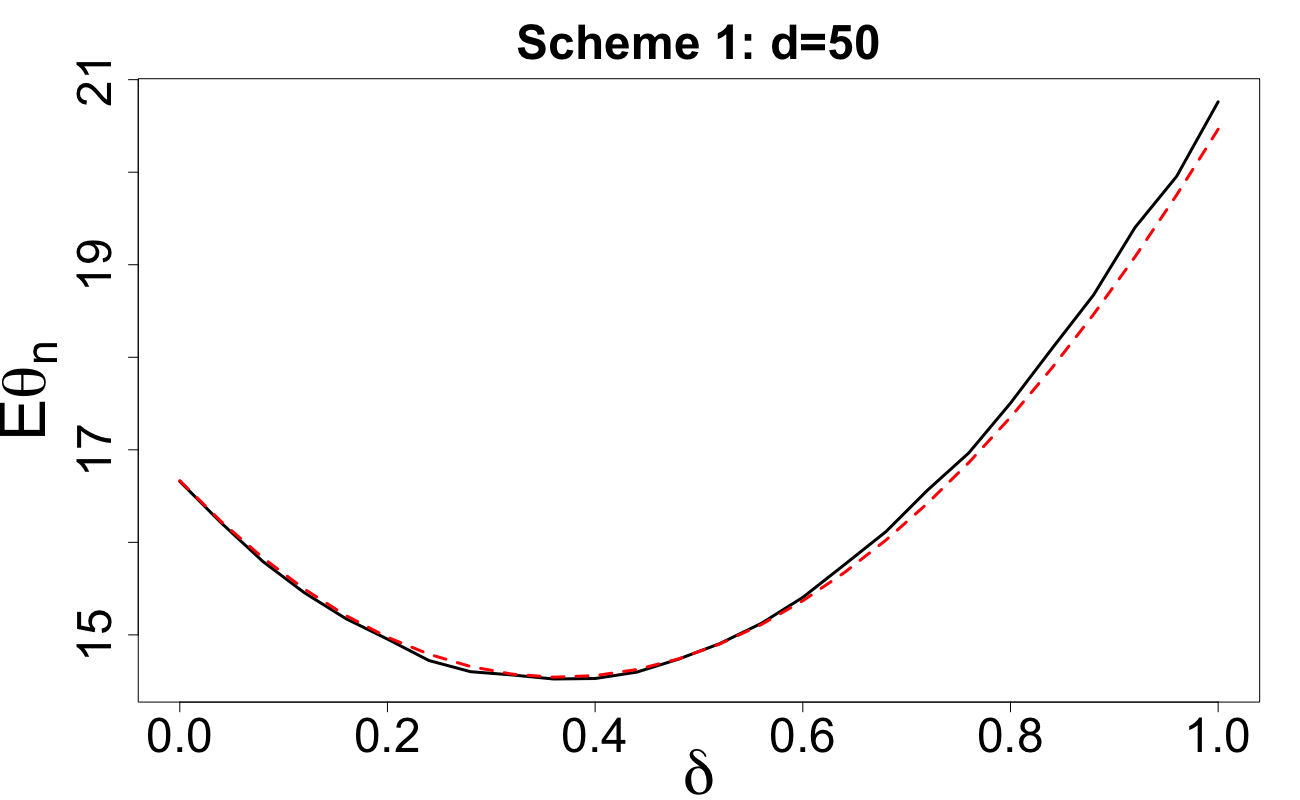}
\caption{$\mathbb{E} \theta_n$ and approximation \eqref{quant_approx_corrected}; $n=128$.}
\label{quant_approx_128}
\end{minipage}%
\begin{minipage}{.5\textwidth}
  \centering
 \includegraphics[width=1\textwidth]{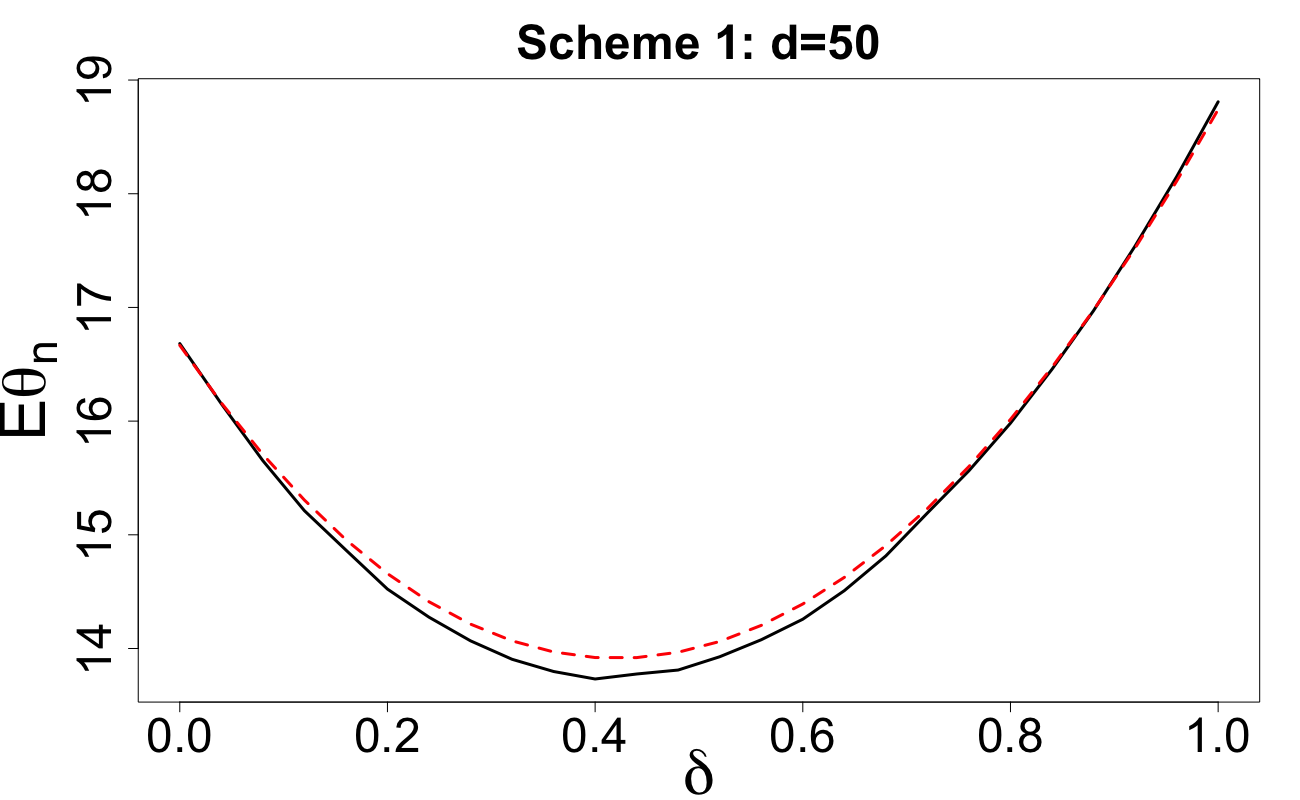}
\caption{$\mathbb{E} \theta_n$ and approximation \eqref{quant_approx_corrected}; $n=512$.}
\label{quant_approx_512}
\end{minipage}
\end{figure}

As follows from results of \cite[Ch.6]{niederreiter1992random}, for efficient covering schemes the order of convergence of the covering radius to 0 as $n \to \infty$ is $n^{-1/d}$. Therefore, for the mean squared distance (which is the quantization error) we should expect the order $n^{-2/d}$ as $n \to \infty$. Therefore, for sake of comparison of quantization errors $\theta_n$ across $n$ we renormalize this error from $\mathbb{E} \theta_n$ to $n^{2/d} \mathbb{E} \theta_n  $.

In Tables~\ref{d_10_quantization}--\ref{d_50_quantization}, we present the minimum value of $n^{2/d} \mathbb{E} \theta_n  $ for a selection of the schemes among those considered in Section~\ref{sec:main_1}. In these tables, the value within the brackets corresponds to the value of $\delta$ where the minimum of $n^{2/d} \mathbb{E} \theta_n$ was obtained. For Scheme~3, typical behaviour of $\mathbb{E} \theta_n$ across $\delta$ for a number and $n$ and $d$  is presented in Figures~\ref{Scheme_3_quant_d_10}--\ref{Scheme_3_quant_d_50}.

\begin{table}[h]
\centering
\begin{tabular}{ |p{2.5cm}||p{2cm}|p{2cm}|p{2cm}|p{2cm}|   }
 \hline
 \multicolumn{5}{|c|}{$d=10$} \\
 \hline
  & $n=64$ &$n=128$& $n=512$ & $n=1024$\\
 \hline
Scheme 1& 4.153 (0.68)   & 4.105 (0.72) &  3.992 (0.80) &3.925 (0.84)  \\
Scheme 3 & 3.663 (0.40) &3.548 (0.44)   & 3.221 (0.48)  & 3.348  (0.52)   \\
Scheme 4,  $\alpha=0.5$& 4.072 (0.56)  &  4.013 (0.60) & 3.839 (0.68) & 3.770   (0.69) \\
Scheme 7&   3.998 (0.68)  & 3.973 (0.76) & 3.936 (0.80) &   3.834 (0.82) \\
Scheme 7, $\delta=1$& 4.569   (1.00)  & 4.425  (1.00) & 4.239 (1.00) & 4.094  (1.00) \\
 \hline
\end{tabular}
\caption{Minimum value of $n^{2/d} \mathbb{E} \theta_n$ and $\delta$ (in brackets) across schemes and $n$ for $d=10$.  }
\label{d_10_quantization}
\end{table}

\begin{table}[h]
\centering
\begin{tabular}{ |p{2.5cm}||p{2cm}|p{2cm}|p{2cm}|p{2cm}|   }
 \hline
 \multicolumn{5}{|c|}{$d=20$} \\
 \hline
  & $n=64$ &$n=128$& $n=512$ & $n=1024$\\
 \hline
Scheme 1& 7.552 (0.52)   & 7.563 (0.56)  &  7.528 (0.64)  & 7.484 (0.68)  \\
Scheme 3 &  7.298 (0.32) &  7.270  (0.33)&  7.133 (0.36) & 7.016 (0.40)  \\
Scheme 4,  $\alpha=0.5$& 7.541 (0.40) &  7.515 (0.44) & 7.457 (0.52)  &  7.421 (0.54) \\
Scheme 7&  7.445 (0.48)   &7.464 (0.56)  &7.487 (0.64)  &  7.453 (0.66)  \\
Scheme 7, $\delta=1$& 9.089   (1.00)  &  9.133 (1.00) &  8.87 (1.00) & 8.681  (1.00) \\
 \hline
\end{tabular}
\caption{Minimum value of $n^{2/d} \mathbb{E} \theta_n$ and $\delta$ (in brackets) across schemes and $n$ for $d=20$.  }
\end{table}

\begin{table}[h]
\centering
\begin{tabular}{ |p{3cm}||p{2cm}|p{2cm}|p{2cm}|  }
 \hline
 \multicolumn{4}{|c|}{$d=50$} \\
 \hline &$n=128$& $n=512$ & $n=1024$\\
 \hline
Scheme 1& 17.608 (0.36)   & 17.634 (0.40) & 17.643 (0.44)    \\
Scheme 3 & 17.483  (0.20) &  17.511 (0.24) & 17.554 (0.27)   \\
Scheme 4,  $\alpha=0.5$& 17.590 (0.28) & 17.670 (0.36)   &  17.620 (0.38)    \\
Scheme 7, $\delta=1$& 20.196 (1.00) & 21.231  (1.00) & 21.711  (1.00) \\
 \hline
\end{tabular}
\caption{Minimum value of $n^{2/d} \mathbb{E} \theta_n$ and $\delta$ (in brackets) across schemes and $n$ for $d=50$.  }
\label{d_50_quantization}
\end{table}

\begin{figure}[h]
\centering
\begin{minipage}{.5\textwidth}
  \centering
  \includegraphics[width=1\textwidth]{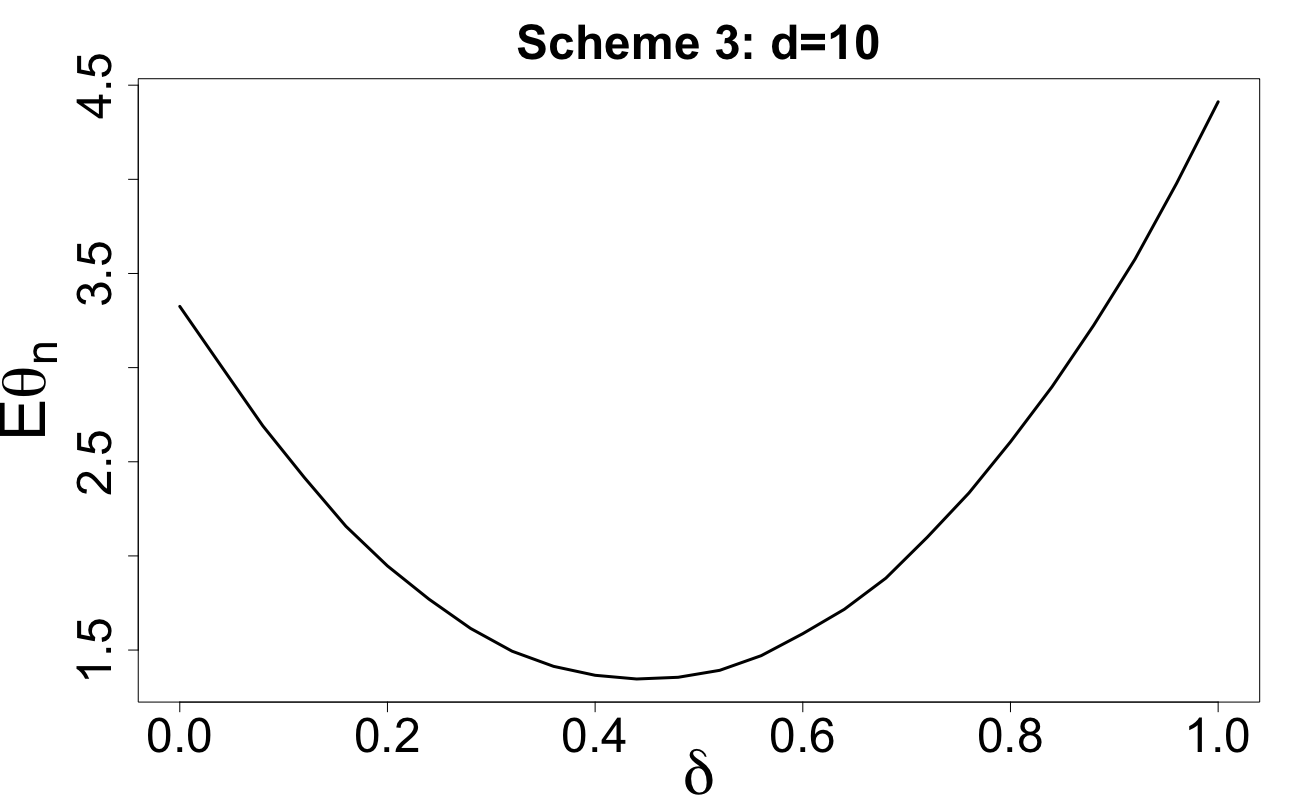}
\caption{$\mathbb{E} \theta_n$ with $n=128$.}
\label{Scheme_3_quant_d_10}
\end{minipage}%
\begin{minipage}{.5\textwidth}
  \centering
 \includegraphics[width=1\textwidth]{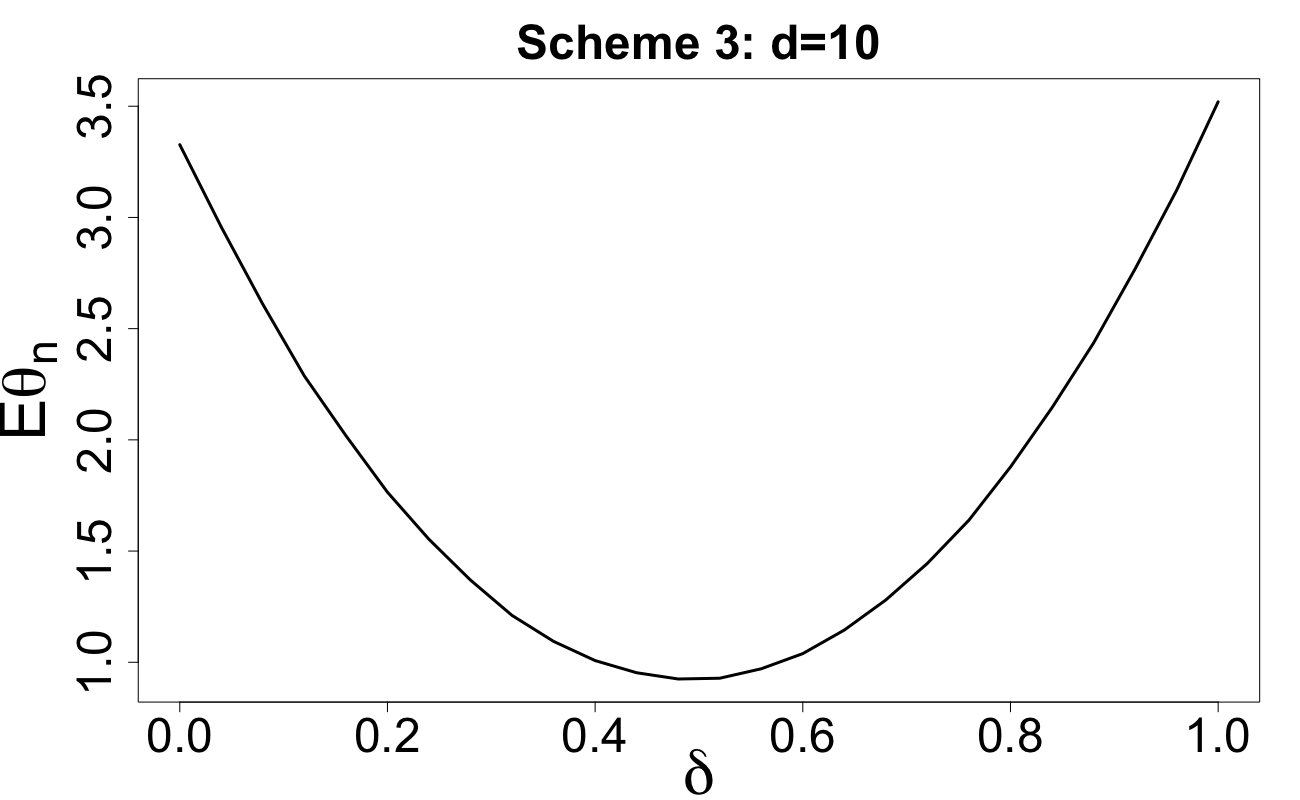}
\caption{$\mathbb{E} \theta_n$ with $n=512$.}
\end{minipage}
\end{figure}

\begin{figure}[h]
\centering
\begin{minipage}{.5\textwidth}
  \centering
  \includegraphics[width=1\textwidth]{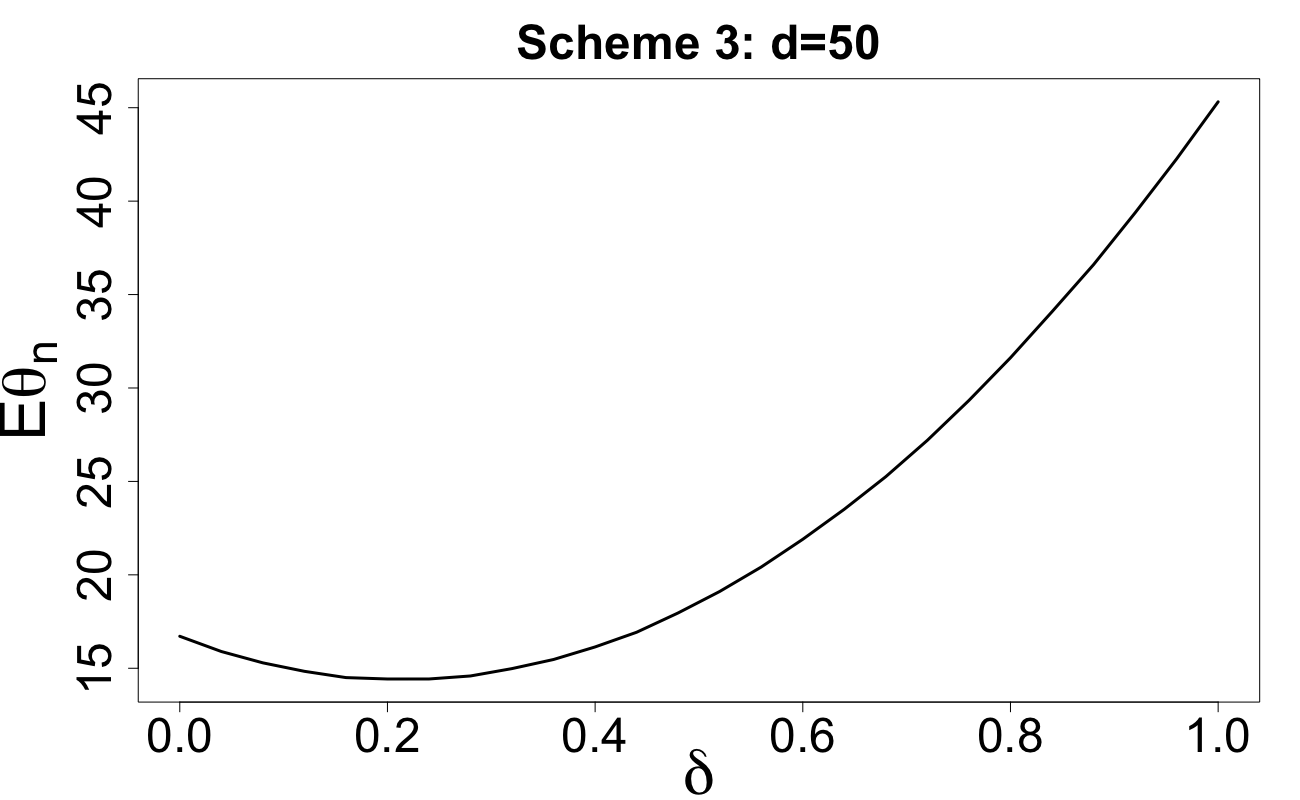}
\caption{$\mathbb{E} \theta_n$ with $n=128$.}
\end{minipage}%
\begin{minipage}{.5\textwidth}
  \centering
 \includegraphics[width=1\textwidth]{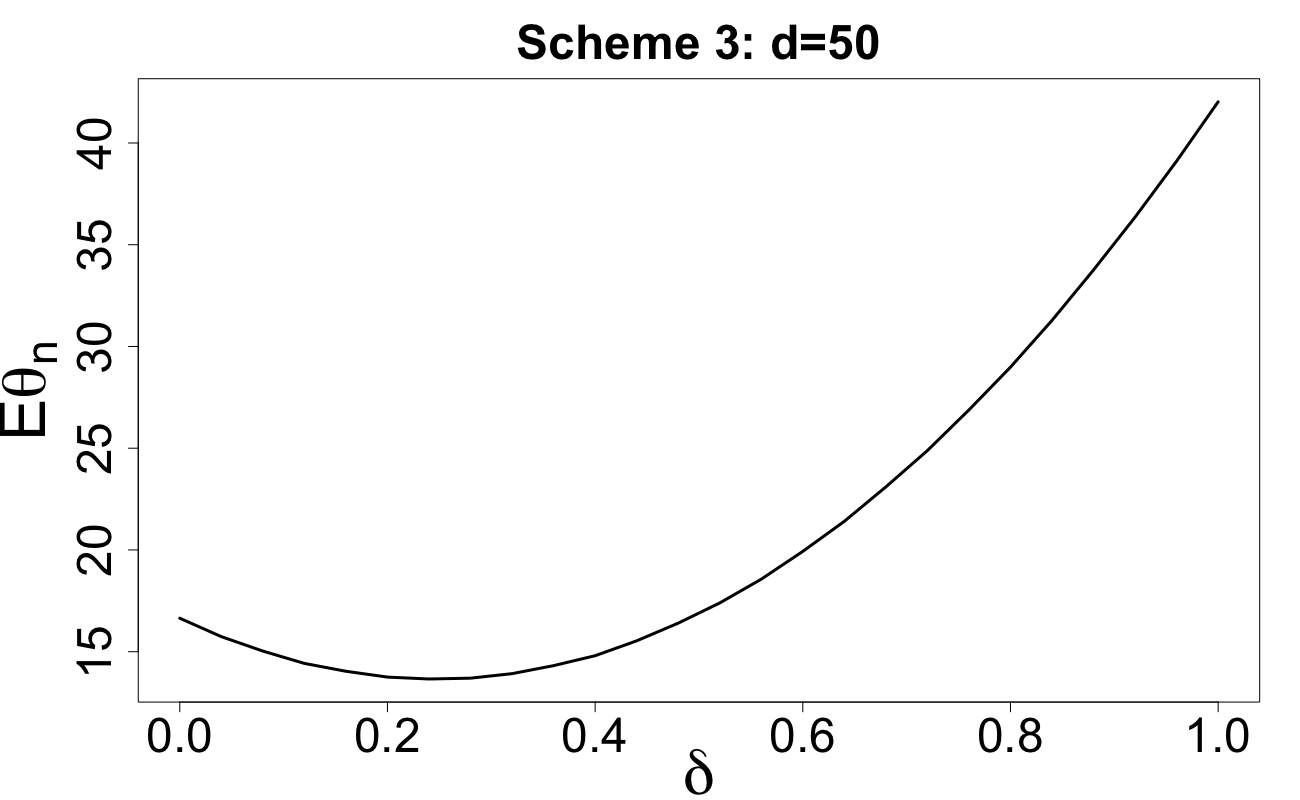}
\caption{$\mathbb{E} \theta_n$ with $n=512$.}
\label{Scheme_3_quant_d_50}
\end{minipage}
\end{figure}

We make the following two main conclusions from analyzing results of this numerical study:

\begin{itemize}
  \item[(a)] the presence of a strong $\delta$-effect, very similar to the effect observed in Section~\ref{sec:main_1}, and
  \item[(b)] for a given design $\mathbb{Z}_n$, there is a very strong correlation between the covering probability as studied   in Section~\ref{sec:main_1} and the normalized quantization
  error $n^{2/d}\mathbb{E} \theta(\mathbb{Z}_n)$.
  \end{itemize}


By comparing the values of $\delta$ in Tables~\ref{d_10_quantization}--\ref{d_50_quantization} with Tables~\ref{table_d_10}--\ref{Table_d_50}, we see a strong similarity between efficient quantization schemes and efficient covering schemes.

\section{Appendix A: Several facts about $d$-dimensional balls and cubes}
\label{sec:appA}

In this appendix,  we  briefly mention several  facts, used in the main part of the paper, related to high-dimensional cubes and balls. Many of these facts are somewhat counter-intuitive and often lead to creation of wrong heuristics in multivariate optimization and misunderstanding of the behaviour of even simple algorithms in high-dimensional spaces. For more details concerning the material of Sections~\ref{sec:Vball}-\ref{sec:concentration}, see~\cite{blum2020foundations}.

\subsection{Volume of the  ball}
\label{sec:Vball}

The volume of the ball $
{\cal B}_d({ r })= \{ x \in \mathbb{R}^d: \| x \| \leq { r }\}
$  can be computed by the formula
\be
\label{eq:unitball}
\mbox{${\rm vol}({\cal B}_d({ r }))= { r }^d  V_d$, where}\;\; V_d= {\rm vol}({\cal B}_d(1))= \frac{  \pi^{d/2} }{ \Gamma (d/2+1)}\,.
 \ee
 The volumes $V_d$ decrease very fast as $d$ grows.
For example, $V_{100}\simeq 2.368 \cdot 10^{-40}$. As $ d \to \infty$,
\be
\label{eq:unitball5}
V_d^{1/d} \simeq  {\sqrt {2\pi e}} \frac 1{\sqrt{d}}  +O \left( \frac {\log d}{{d}^{3/2}} \right)\, .
\ee

\subsection{Radius of the ball of unit volume}

Define $r_d$ by ${\rm vol} ({\cal B}_d({ r }_d))=1$. Table~\ref{Radius_of_ball} gives approximate values of $r_d$.\\

\begin{table}[h]
\centering
\begin{tabular}{|c||c|c|c|c|c|c|c|c|c|c|}
  \hline
  $d$ & 1& 2 & 3 & 4 & 5 & 6& 7& 8& 9 \\ \hline
  $r_d$ &0.5 & 0.564 & 0.62 & 0.671 & 0.717 & 0.761 &0.8 &0.839 &0.876 \\
  \hline  \hline
$d$ & 10 & 20 & 30 & 40 & 50& 100& 200& 500 &1000 \\ \hline
  $r_d$  &0.911 & 1.201 & 1.43 & 1.626 & 1.8 &2.49 &3.477 &5.45&7.682 \\
  \hline
\end{tabular}
\caption{Radius of the ball of unit volume for different dimensions}
\label{Radius_of_ball}
\end{table}

From \eqref{eq:unitball5}, for large $d$
we have
\bea
\label{eq:unitball6}
r_d=\frac {\sqrt{d}}{\sqrt {2\pi e}} +O \left( \frac 1{\sqrt{d}} \right)\, ,
\eea
where  $1/{\sqrt {2\pi e}} \simeq 0.242$. This is only about twice smaller than $\sqrt{d}/2$, the length of the half-diagonal  of the $d$-dimensional unit cube $[0,1]^d$.

For
 $r_{d,2 \delta}$ defined by ${\rm vol} ({\cal B}_d(r_{d,2 \delta}))= {\rm vol} ({\cal C}_{d}(\delta) ) = (2 \delta)^d$, we  have
$ r_{d,2 \delta}= 2 \delta r_d$.

\subsection{Almost all the volume is near the boundary}
\label{sec:boundary}

First, consider  the  cube ${\cal C}_{d}(\delta)=[-\delta,\delta]^d$, with $0<\delta<1$, as interior to   the  cube ${\cal C}_{d}=[-1,1]^d$.
For the ratio of the volumes of these two cubes, we  have ${{\rm vol}({\cal C}_{d}(\delta))} /{{\rm vol}({\cal C}_{d})}= \delta^d $
which  tends to 0 (as $d \to \infty$) exponentially fast for   any $\delta \in (0,1)$.

If, as $d \to \infty$, $\delta $ changes getting closer to 1 but $1-\delta$ tends to 0 slower than $1/d$, then the ratio of the two volumes still tends to 0.
In particular, if $1-\delta=c/d^{1-\delta}$ with $0<\delta<1$ then
\bea
\frac{{\rm vol}({\cal C}_{d}(\delta))} {{\rm vol}({\cal C}_{d})}= \delta^d \simeq \exp\{- c d^{1-\delta} \} \to 0 \, , \;\;\; d \to \infty\, .
\eea

Consider now the balls ${\cal B}_d(1)$ and ${\cal B}_d(1-\epsilon)$. The difference ${\cal B}_d(1)\setminus {\cal B}_d(1-\epsilon)$ is called the annulus.
Using \eqref{eq:unitball} we can compute the ratio of volume of this annulus to the volume of the unit ball:
\bea
\frac{{\rm vol}\left[{\cal B}_d(1)\setminus {\cal B}_d(1-\epsilon)\right]}{{\rm vol}({\cal B}_d(1))}= 1-\eps^d \, .
\eea
This ratio tends to 1 exponentially fast as $d \to \infty$.
The ratio of volume of the ball ${\cal B}_d(1-\epsilon)$ to the volume of the unit ball ${\cal B}_d(1)$ is, similarly to  the case of the cubes above, $(1-\eps)^d$.
This result extends to any measurable set $A \subset \mathbb{R}^d$. Indeed, define the set
$A_{1-\ve}=\{(1-\ve)x \,: x \in A\}$. Then, by splitting $A$ and  $A_{1-\ve}$ into infinitesimal cubes and adding up their volumes, we find
${\rm vol}(A_{1-\ve})=(1-\ve)^d {\rm vol}(A)\, .$

\subsection{The area of volume concentration in a cube }
\label{sec:concentration}

Let $X=(x_1, \ldots, x_d)$ be uniformly distributed on ${\cal C}_{d} =[-1,1]^d$. Then  $x_1^2, \ldots, x_d^2$ are independent r.v. on    $[0,1]$.
The Hoeffding's inequality gives
\bea
\mathbb{P} \left\{ \big| \;\frac{1}{d}(x_1^2+ \ldots+ x_d^2)- \frac{1}{d}\mathbb{E}\left( x^2_1+ \ldots+ x^2_d \right) \;\big| \geq \epsilon  \right\} \leq 2 e^{-2d\epsilon^2 } \, .
\eea
Since $\mathbb{E} x^2_i= \frac13$, we obtain
\bea
\mathbb{P} \left\{ \big| \; \|X\| ^2 -\frac{d}{3} \big|
\geq \epsilon d  \right\} \leq 2 e^{-2d\epsilon^2 } \, .
\eea
Therefore, the main volume in the cube ${\cal C}_{d}$ is concentrated in the annulus around the sphere with radius $\sqrt{d/3}$.

\subsection{Squared norm of a  random  point in a cube}

Let $Z=(z_1, \ldots, z_d)$ be a  random vectors on ${\cal C}_{d}(\delta)=[-\delta,\delta]^d$
consisting of i.i.d. random components $z_i$ having a  distribution with density $p(t)$, $t \in [-\delta,\delta]$, $\delta>0$.

Set $\eta=\sum_{j=1}^d z_j^2$. We have $\mathbb{E} \eta= d \mu_2 $ and ${\rm var}( \eta)=d {\rm var}( z_1^2)=d(\mu_4-\mu_2^2)$,
where $\mu_j$ be the moments of the distribution with density $p(t)$.

For example, when $z_i$  have Beta$(\alpha,\alpha)$ distribution with density
\be
\label{eq:beta1}
p_{\alpha,\delta}(t)= \frac{(2\delta)^{1-2\alpha}}{\mbox{Beta$(\alpha,\alpha)$}} [\delta^2-t^2]^{\alpha-1}\, , \;\;-\delta<t<\delta\, ,\alpha>0,
\ee
where Beta$(\cdot,\cdot)$ is the Beta-function, then
\be
\label{eq:beta2}
 \mu_2=\frac{\delta^2}{2\alpha+1}\,, \;   \mu_4=\frac{3 \delta^4}{(2\alpha+1)(2\alpha+3)}\,
\ee
and therefore
\be
\label{eq:moments_hd0}
\mathbb{E} \eta= \frac{d\delta^2}{2\alpha+1}\,,\;\;  {\rm var}( \eta)= \frac{4d \delta^4 \alpha}{(2\alpha+1)^2(2\alpha+3)}\, .
\ee
If $\alpha=1$, when $Z$ is uniform in the cube ${\cal C}_{d}(\delta)$, then
\be
\label{eq:moments_hd_u0}
\mathbb{E} \eta= \frac13 d\delta^2\,,\;\;  {\rm var}( \theta)= \frac{4}{45} d \delta^4\, .
\ee

\subsection{Distance between two random  points in a cube}

Assume $Z=(z_1, \ldots, z_d)$ and $Z^{\prime}=(z^{\prime}_1, \ldots, z^{\prime}_d)$ are independent random vectors on ${\cal C}_{d}(\delta)=[-\delta,\delta]^d$ consisting of i.i.d. random components $z_i$ and $z^{\prime}_i$ which have some distribution with density $p(t)$, $t \in [-\delta,\delta]$, $\delta>0$. Let $\mu_j$ be the moments of the distribution with density $p(t)$.
Assume that the density $p(t)$ is symmetric around 0 and hence all odd moments are zero: $\mu_{2k+1}=0$ for $k=1,2, \ldots$

The distribution of the squared distances
$$ \theta=\|Z-Z^{\prime}\|^2 = { \sum_{i=1}^d (z_i-z^{\prime}_i)^2 } $$
has the mean and variance that can be easily computed as follows:
\bea
\mathbb{E} \theta= d \mathbb{E}(z_1-z^{\prime}_1)^2 &=& 2d \mu_2\, ,\\
{\rm var}( \theta)= d {\rm var}(z_1-z^{\prime}_1)^2& =& d \left[[\mathbb{E}(z_1-z^{\prime}_1)^4-[\mathbb{E}(z_1-z^{\prime}_1)^2)]^2 \right]= 2d \left[\mu_4+ \mu_2^2  \right]
\eea
For example, when $z_i$ and $z^{\prime}_i$ have Beta$(\alpha,\alpha)$ distribution with density \eqref{eq:beta1} and hence moments \eqref{eq:beta2}, we obtain
\be
\label{eq:moments_hd}
\mathbb{E} \theta= \frac{2d\delta^2}{2\alpha+1}\,,\;\;  {\rm var}( \theta)= \frac{4d \delta^4 (4\alpha+3)}{(2\alpha+1)^2(2\alpha+3)}\, .
\ee
If $\alpha=1$ (that is, when $Z$ and $Z^{\prime}$ are uniform in the cube ${\cal C}_{d}(\delta)$), then
\be
\label{eq:moments_hd_u}
\mathbb{E} \theta= \frac23 d\delta^2\,,\;\;  {\rm var}( \theta)= \frac{28}{45} d \delta^4
\ee

\subsection{Volume of the intersection of two balls of the same radius}
\label{sec:inter}

Let ${\cal B}_d(Z_j,r) $ and $ {\cal B}_d(Z_i,r)$ be two balls in $\mathbb{R}^d$ with same radius and different centers $Z$ and $Z^\prime$.
To compute the volume of the intersection ${\cal B}_d(Z,r) \cap {\cal B}_d(Z^\prime,r)$,  we will use the formula, see , for the volume of the
$d$-dimensional cap (cut in the direction of $Z'$) of height~$h$ from a $d$-dimensional ball ${\cal B}_{d}(Z,r) $:
\be
\label{eq:inters0f}
K_{d,r,h}=\frac12  r^d V_d I_{1-h^2/r^2}\left(\frac{d-1}{2},\frac12 \right) - \frac{h}{d}(r^2-h^2)^{(d-1)/2} V_{d-1}\, ,
\ee
where $V_d$ is defined in \eqref{eq:unitball}, $\Gamma(\cdot)$ is the Gamma-function and
$$
I_t(\alpha,\beta)= \int_{0}^{t} u^{\alpha-1}(1-u)^{\beta-1} du \bigg/ \int_{0}^{1} u^{\alpha-1}(1-u)^{\beta-1} du
$$
is the normalised incomplete Beta-function. In the rhs of \eqref{eq:inters0f}, the first term is the volume of the related $d$-dimensional hyper-sector (this expression  is derived in \cite{li2011concise}) and the second term is the volume of the cone with height $h$ and base
${\cal B}_{d-1}((Z+Z')/2,r') $, where $r'=   \sqrt{r^2-h^2}$.

The volume of the intersection of  the  balls ${\cal B}_d(Z,r)$ and $ {\cal B}_d(Z^\prime,r)$
is therefore
\be
\label{eq:int3}
\mbox{
vol(${\cal B}_d(Z,r) \cap {\cal B}_d(Z^\prime,r))=2K_{d,r,h} $}
\ee where $h=\frac12 \|Z-Z^\prime||$ and $K_{d,r,h}$ is defined in
\eqref{eq:inters0f}.

\subsection{ A  direct  computation of $C_{d,Z,{ r }}$}
\label{sec:quantuty44}

For computing values of $C_{d,Z,{ r }}$, we  can employ the following direct approach based on the use of characteristic functions (c.f.).

\begin{itemize}
  \item[(a)]$\,$  Compute  the c.f.
  $  \psi_z(s)= \int e^{its} \varphi_z(t)dt$
  for $z=z_j$ ($j=1,\ldots, d$), with the density $\varphi_z(t)$  defined either by \eqref{eq:inters2d} or \eqref{eq:inters3d}.

  \item[(b)]$\;$As $u_j$ are independent, the c.f. of $\| U-Z \|^2$ is the product $\psi_{Z}(s)=\prod_{j=1}^{d}\psi_{z_j}(s)$.
  \item[(c)]$\;$The density of $\| U-Z \|^2$ is found using the inversion formula
  $$
  p_{d,Z} (x)= \frac1{2\pi} \int_{-\infty}^\infty  e^{-isx} \psi_{Z}(s) ds\, , \;\;  x \geq 0\, .
  $$
\end{itemize}

For computing the c.f. $  \psi_z(s)= \int e^{its} \varphi_z(t)dt$ we can use the formula
$$
\int_a^b \frac{e^{xt}}{\sqrt{t}}  dt = 2\int_{\sqrt{a}}^{\sqrt{b}} e^{x u^2}  du= \sqrt{\frac{\pi}{x}} \left( {\rm erfi}(\sqrt{bx})-{\rm erfi}(\sqrt{ax})\right)
$$
for any $0 \leq a<b<\infty$ and any complex $x\neq 0$.
Here ${\rm erfi}(x)$ is the imaginary error function
$$
{\rm erfi}(x)=\frac{2}{\sqrt{\pi}} \int_{0}^{x}e^{t^2} dt= \frac{2}{\sqrt{\pi}}  \sum_{j=0}^{\infty} \frac{x^{2j+1}}{j!(2j+1)}\, ;
$$
the series in the right-hand side of this formula converges for all complex $x$.

This approach allows very accurate computation of $C_{d,Z,{ r }}$ but it is very computationally intensive and  can only be performed for given $Z$.

\section{Appendix B: Important auxiliary results}
\label{sec:appB}

{\bf Lemma 1.}
{\it Let $\delta>0$,  $x  \in \mathbb{R} $ and $\eta_{x,\delta}$ be a r.v. $\eta_{x,\delta} = (\xi-x)^2$, where
r.v.~$\xi$ has uniform distribution on $[-\delta,\delta]$.
Then the c.d.f. of the r.v. $\eta_{x,\delta} $ is
\be
\label{eq:inters1b}
F_{x,\delta}(t) = \mathbb{P}\{\eta_{x,\delta} \leq t \}= \left\{\begin{array}{ll}
0 & {\rm for\;\;} t\leq 0\\
             \frac{\sqrt{t}}{\delta}  \cdot 1_{[\,|x|\leq \delta]}& {\rm for\;\;} 0<t< (\delta- |x|)^2\\
             \frac{\delta  -|x|+\sqrt{t}}{2\delta}   & {\rm for\;\;}(\delta- |x|)^2 \leq t\leq (\delta+ |x|)^2 \\
              1 &  (\delta+ |x|)^2 < t\,,
            \end{array}
            \right.\;\;\;\;\;\;\;
\ee
where
$$
1_{[\,|x|\leq \delta]}= \left\{\begin{array}{ll}
1 & {\rm if\;\;} |x | \leq \delta\\
0 & {\rm if\;\;} |x | > \delta\, .
\end{array}
            \right.
$$
 The corresponding  density of  $\eta_{x,\delta} $ is
\be
\label{eq:inters1a}
\varphi_{x,\delta}(t) = \left\{\begin{array}{ll}
             1/(2\delta\sqrt{t}) \cdot 1_{[\,|x|\leq \delta]} & {\rm for\;\;} 0<t< (\delta- |x|)^2\\
             1/(4\delta\sqrt{t}) & {\rm for\;\;} (\delta- |x|)^2 < t\leq (\delta+ |x|)^2 \\
              0 & {\rm otherwise.}
            \end{array}
            \right.
\ee
The first four  central moments of  the r.v. $\eta_{x,\delta} $ are:
\be
\label{eq:inters1c}
\mu_{x,\delta}^{(1)} = E\eta_{x,\delta} =x^2 +\frac{\delta^2}{3}, \;\;\mu_{x,\delta}^{(2)} = {\rm var}(\eta_{x,\delta}) = \frac{4\delta^2}3 \left(x^2 +\frac{\delta^2}{15} \right)\, ,
\ee
\be
\label{eq:inters1c3}
\mu_{x,\delta}^{(3)}=E \left[\eta_{x,\delta} - E\eta_{x,\delta}\right]^3 = \frac{16\delta^4}{15} \left(x^2 +\frac{\delta^2}{63} \right)  \, ,
\ee
\be
\label{eq:inters1c4}
\mu_{x,\delta}^{(4)}=E \left[\eta_{x,\delta} - E\eta_{x,\delta}\right]^4 =
 3 \mu_{x,\delta}^{(1)} \mu_{x,\delta}^{(3)}
\, .
\ee

}

{\bf Proof.} Clearly, if $t\leq 0$ then  $F_{x,\delta}(t)=0$ and so  we only consider the case  $t>0$.
In view of symmetry, for all $x\in \mathbb{R}$, $\delta>0$ and $t \geq 0$, we have
$F_{x,\delta}(t)= F_{-x,\delta}(t)$ and therefore we only need to consider $x \geq 0$. Also, $\eta_{x,\delta} \leq (|x|+\delta)^2$ with probability 1 implying $F_{x,\delta}(t)=1 $ for all $t \geq (|x|+\delta)^2$.

Assume   $0 \leq x \leq \delta$. We then have for all $ t \geq 0$:
\bea
F_{x,\delta}(t)&=& \mathbb{P}\{(\xi-x)^2 \leq t \}= \mathbb{P}\{(\xi-x)^2 \leq t, \; \xi \leq x \}+\mathbb{P}\{(\xi-x)^2 \leq t, \; \xi > x \}\\
&=& \mathbb{P}\{x-\xi \leq \sqrt{ t}, \; \xi \leq x \}+\mathbb{P}\{\xi-x \leq \sqrt{ t}, \; \xi > x \}\\
&=& \mathbb{P}\{x- \sqrt{ t} \leq \xi \leq x \}+\mathbb{P}\{x<\xi \leq x+ \sqrt{ t} \}
\eea
with
\bea
\label{eq:inter0}
\mathbb{P}\{x- \sqrt{ t} \leq \xi \leq x \}= \left\{
                                               \begin{array}{ll}
                                                \sqrt{ t}/(2 \delta)  & {\rm if}\; \sqrt{ t}<x+\delta \\
                                               (x+\delta)/(2 \delta) & {\rm if}\; \sqrt{ t} \geq  x+\delta \, ,
                                               \end{array}
                                             \right.
\eea
\bea
\label{eq:inter01}
\mathbb{P}\{x<\xi \leq x+ \sqrt{ t} \}= \left\{
                                               \begin{array}{ll}
                                                \sqrt{ t}/(2 \delta)  & {\rm if}\; \sqrt{ t}<\delta -x \\
                                               (\delta -x)/(2 \delta) & {\rm if}\; \sqrt{ t} \geq \delta -x\, .
                                               \end{array}
                                             \right.
\eea
This yields the expression \eqref{eq:inters1b} for $F_{x,\delta}(t)$ in the case $|x|\leq \delta$.

If $x > \delta$ then  $\eta_{x,\delta} \geq (x-\delta)^2$ with probability 1 implying  $F_{x,\delta}(t)=0$ for all $t\leq (x-\delta)^2$ and   $\mathbb{P}\{x<\xi \leq x+ \sqrt{ t} \}=0$ for all $t$. Therefore
\bea
F_{x,\delta}(t)=
\mathbb{P}\{x- \sqrt{ t} \leq \xi \leq x \}= \left\{
                                               \begin{array}{ll}
                                               0 & {\rm if}\; \sqrt{ t} \leq  x-\delta \, \\
                                               \frac{\delta-(x- \sqrt{ t})}{2 \delta}  & {\rm if}\; x-\delta<\sqrt{ t}<x+\delta \\
                                               1 & {\rm if}\; \sqrt{ t} \geq  x+\delta \, ,
                                               \end{array}
                                             \right.
\eea
This yields the expression \eqref{eq:inters1b} for $F_{x,\delta}(t)$ in the case $|x|> \delta$.

 Deduction of the formulas \eqref{eq:inters1a} for the density and \eqref{eq:inters1c} for the moments from the expression \eqref{eq:inters1b} for the c.d.f. $F_{x,\delta}(t)$  is an easy exercise.

 \hfill$\Box$

%

{\bf Lemma 2.}
{\it Let $\delta>0$,  $x  \in \mathbb{R} $ and $\eta_{x,\delta}^\prime$ be a r.v. $\eta_{x,\delta}^\prime = |\xi-x|$, where
r.v.~$\xi$ has uniform distribution on $[-\delta,\delta]$.
Then the c.d.f. of the r.v. $\eta_{x,\delta}^\prime $ is
\be
\label{eq:inter2b}
F_{x,\delta}^\prime(t) = \mathbb{P}\{\eta_{x,\delta}^\prime \leq t \}= \left\{\begin{array}{ll}
0 & {\rm for\;\;} t\leq 0\\
             \frac{{t}}{\delta}  \cdot 1_{[\,|x|\leq \delta]}& {\rm for\;\;} 0<t< |\delta- |x|\, |\\
             \frac{\delta  -|x|+{t}}{2\delta}   & {\rm for\;\;}|\delta- |x|\, | \leq t\leq \delta+ |x| \\
              1 &  \delta+ |x| < t\,,
            \end{array}
            \right.\;\;\;\;\;\;\;
\ee

 The corresponding  density of  $\eta_{x,\delta}^\prime $ is
\be
\label{eq:inter2a}
\varphi_{x,\delta}^\prime(t) = \left\{\begin{array}{ll}
             \frac1{\delta} \cdot 1_{[\,|x|\leq \delta]}& {\rm for\;\;} 0<t< |\delta- |x|\, |\\
             \frac1{2\delta} & {\rm for\;\;} |\delta- |x|\,| < t\leq \delta+ |x| \\
              0 & {\rm otherwise.}
            \end{array}
            \right.
\ee
}

 Lemma 2 follows from Lemma 1 by noting that $\eta_{x,\delta}^\prime = \sqrt{\eta_{x,\delta} } $.

Note that $1_{[\,|x|\leq \delta]}=0$ for   $|x|>\delta$ and one of the two non-trivial cases in \eqref{eq:inters1b},   \eqref{eq:inters1a},
\eqref{eq:inter2b}  and \eqref{eq:inter2a}, when
$|x|>\delta$, become trivial as expressions vanish to  zero.

\bibliographystyle{plain}
\bibliography{large_dimension}

\begin{thebibliography}{10}

\bibitem{blum2020foundations}
A.~Blum, J.~Hopcroft, and R.~Kannan.
\newblock {\em Foundations of data science}.
\newblock Cambridge University Press, 2020.

\bibitem{janson1986random}
S.~Janson.
\newblock Random coverings in several dimensions.
\newblock {\em Acta Mathematica}, 156:83--118, 1986.

\bibitem{januszewski1994line}
J.~Januszewski and M.~Lassak.
\newblock On-line covering the unit cube by cubes.
\newblock {\em Discrete \& Computational Geometry}, 12(4):433--438, 1994.

\bibitem{joe2008constructing}
S.~Joe and F.~Y. Kuo.
\newblock Constructing {S}obol sequences with better two-dimensional
  projections.
\newblock {\em SIAM Journal on Scientific Computing}, 30(5):2635--2654, 2008.

\bibitem{johnson1990minimax}
M.~E. Johnson, L.~M. Moore, and D.~Ylvisaker.
\newblock Minimax and maximin distance designs.
\newblock {\em Journal of statistical planning and inference}, 26(2):131--148,
  1990.

\bibitem{kuperberg1994line}
W.~Kuperberg.
\newblock On-line covering a cube by a sequence of cubes.
\newblock {\em Discrete \& Computational Geometry}, 12(1):83--90, 1994.

\bibitem{li2011concise}
S.~Li.
\newblock Concise formulas for the area and volume of a hyperspherical cap.
\newblock {\em Asian Journal of Mathematics and Statistics}, 4(1):66--70, 2011.

\bibitem{niederreiter1992random}
H.~Niederreiter.
\newblock {\em Random number generation and quasi-Monte Carlo methods}.
\newblock Siam, 1992.

\bibitem{second_paper}
J.~Noonan and A.~Zhigljavsky.
\newblock Non-lattice covering and quantization in high dimensions.
\newblock In {\em Black Box Optimization, Machine Learning and No-Free Lunch
  Theorems}, pages 799--860. Springer, 2021.

\bibitem{petrov2012sums}
V.~V. Petrov.
\newblock {\em Sums of independent random variables}.
\newblock Springer-Verlag, 1975.

\bibitem{petrov}
V.~V. Petrov.
\newblock {\em Limit theorems of probability theory: sequences of independent
  random variables}.
\newblock Oxford Science Publications, 1995.

\bibitem{rao1987asymptotic}
B.L.S. Prakasa~Rao.
\newblock {\em Asymptotic theory of statistical inference}.
\newblock Wiley, 1987.

\bibitem{pronzato2012design}
L.~Pronzato and W.~G. M{\"u}ller.
\newblock Design of computer experiments: space filling and beyond.
\newblock {\em Statistics and Computing}, 22(3):681--701, 2012.

\bibitem{sukharev2012minimax}
A~Sukharev.
\newblock {\em Minimax models in the theory of numerical methods}.
\newblock Springer Science \& Business Media, 1992.

\bibitem{SIAM}
B.~Tibken, D.~Constales, and et~al.
\newblock The volume of the intersection of a concentric cube and ball in
  n-dimensional space: collection of approximations.
\newblock In {\em SIAM Review}, volume~39, pages 783--786, 1997.

\bibitem{toth20172}
G.~F. T{\'o}th.
\newblock Packing and covering.
\newblock In {\em Handbook of discrete and computational geometry}, pages
  27--66. Chapman and Hall/CRC, 2017.

\bibitem{toth1993packing}
G.~F. T{\'o}th and W.~Kuperberg.
\newblock Packing and covering with convex sets.
\newblock In {\em Handbook of Convex Geometry}, pages 799--860. Elsevier, 1993.

\bibitem{zhigljavsky2012theory}
A.~Zhigljavsky.
\newblock {\em Theory of global random search}.
\newblock Kluwer Academic Publishers, 1991.

\bibitem{zhigljavsky2007stochastic}
A.~Zhigljavsky and A.~Zilinskas.
\newblock {\em Stochastic global optimization}.
\newblock Springer Science \& Business Media, 2007.

\end{thebibliography}

\end{document}